\documentclass{gtart}
\usepackage{amssymb, amsmath, amscd, latexsym, mathrsfs}
\usepackage[all]{xy}
\newtheorem{thm}{Theorem}[subsection]
\newtheorem{prp}[thm]{Proposition}
\newtheorem{cor}[thm]{Corollary}
\newtheorem{lem}[thm]{Lemma}
\newtheorem{rmk}[thm]{Remark}
\newtheorem{dfn}[thm]{Definition}

\newtheorem{exm}[thm]{Example}

\newtheorem{asp}[thm]{Assumption}

\newcommand{\uli}{\underline}
\newcommand{\ep}{\varepsilon}
\newcommand{\lra}{\longrightarrow}

\newcommand{\lag}{\langle}
\newcommand{\rag}{\rangle}
\newcommand{\pt}{\star}
\newcommand{\smin}{\smallsetminus}

\newcommand{\bbul}{_{\bullet\bullet}}

\newcommand{\saddle}{\mathrm{saddle}}
\newcommand{\trace}{\mathrm{trace}}

\newcommand{\id}{\mathrm{id}}
\newcommand{\im}{\mathrm{im}}

\newcommand{\cone}{\mathrm{cone}}
\newcommand{\emb}{\mathrm{emb}}
\newcommand{\map}{\mathrm{map}}
\newcommand{\intr}{\mathrm{int}}
\newcommand{\op}{^{\mathrm{op}}}
\newcommand{\reg}{\mathrm{rg}}

\newcommand{\supp}{\mathrm{supp}}
\newcommand{\loc}{\mathrm{loc}}
\newcommand{\Th}{\mathrm{Th}\,}
\newcommand{\st}{\mathrm{st}}
\newcommand{\es}{\mathrm{es}}
\newcommand{\colim}{\mathrm{colim}}

\newcommand{\hocolim}{\mathrm{hocolim}}
\newcommand{\hofiber}{\mathrm{hofiber}}
\newcommand{\fiber}{\mathrm{fiber}}
\newcommand{\pr}{\mathrm{pr}}
\newcommand{\mor}{\mathrm{mor}}
\newcommand{\ob}{\mathrm{ob}}

\newcommand{\Or}{\mathrm{O}}
\newcommand{\Un}{\mathrm{U}}
\newcommand{\SOr}{\mathrm{SO}}
\newcommand{\GL}{\mathrm{GL}}
\newcommand{\SL}{\mathrm{SL}}
\newcommand{\Gr}{\mathrm{G}}

\newcommand{\iso}{\mathrm{iso}}
\newcommand{\Diff}{\mathrm{Diff}}

\newcommand{\hocolimsub}[1]{\begin{array}[t]{cc} \mathrm{hocolim} \\
[-1.7mm] \scriptstyle{#1} \end{array}}

\newcommand{\ttag}{\refstepcounter{myeqn}\tag{\themyeqn}}

\newcommand{\incto}{\ar@{^{(}->}}
\newcommand{\into}{\ar@{>->}}
\newcommand{\onto}{\ar@{->>}}

\newcommand{\pd}{^{\rule{0mm}{0mm}}}    

\newcommand{\ssmallint}{\!\smallint\!}

\newcommand{\CC}{\mathbb C}
\newcommand{\FF}{\mathbb F}
\newcommand{\ZZ}{\mathbb Z}

\newcommand{\QQ}{\mathbb Q}
\newcommand{\RR}{\mathbb R}

\newcommand{\sA}{\mathscr A}

\newcommand{\sC}{\mathscr C}
\newcommand{\sD}{\mathscr D}
\newcommand{\sE}{\mathscr E}
\newcommand{\sI}{\mathscr I}
\newcommand{\sH}{\mathscr H}

\newcommand{\sK}{\mathscr K}
\newcommand{\sM}{\mathscr M}
\newcommand{\sP}{\mathscr P}
\newcommand{\sR}{\mathscr R}
\newcommand{\sS}{\mathscr S}

\newcommand{\sT}{\mathscr T}

\newcommand{\sW}{\mathscr W}
\newcommand{\sX}{\mathscr X}
\newcommand{\sY}{\mathscr Y}

\newcommand{\cC}{\mathcal C}
\newcommand{\cD}{\mathcal D}
\newcommand{\cE}{\mathcal E}
\newcommand{\cF}{\mathcal F}
\newcommand{\cG}{\mathcal G}

\newcommand{\cL}{\mathcal L}

\newcommand{\cT}{\mathcal T}

\newcommand{\cV}{\mathcal V}
\newcommand{\cW}{\mathcal W}

\newcommand{\cZ}{\mathcal Z}

\newcommand{\cdW}{{}^{\partial}\mathcal W}
\newcommand{\cddW}{{}^{\partial\partial}\mathcal W}
\newcommand{\cdL}{{}^{\partial}\!\mathcal L}

\newcommand{\bP}{\mathbf P}
\newcommand{\bV}{\mathbf h\mathbf V}
\newcommand{\bW}{\mathbf h\mathbf W}

\newcommand{\U}{U}

\renewcommand{\Theta}{\varTheta}
\renewcommand{\Sigma}{\varSigma}
\renewcommand{\Phi}{\varPhi}
\renewcommand{\Psi}{\varPsi}
\renewcommand{\Lambda}{\varLambda}
\renewcommand{\Gamma}{\varGamma}
\renewcommand{\rho}{\varrho}

\begin{document}

\hspace{4in} {\it version: July 14, 2004}
 
\vspace{1in}
\begin{center}
{\bf \LARGE The stable moduli space of Riemann surfaces: \\ Mumford's 
conjecture}
\end{center}

\begin{center} 
{\large Ib Madsen$^{1,*}$ and Michael Weiss$^{2,**}$}
\end{center}

\vspace{0.5in}

\begin{minipage}[t]{3in}
$^{1)}$ Institute for the Mathematical Sciences
\\ Aarhus University \\ 8000 Aarhus C, Denmark \\ 
{\it email}: imadsen@imf.au.dk
\end{minipage} 
\begin{minipage}[t]{2.7in}
$^{2)}$ Department of Mathematics \\
Aberdeen University \\ Aberdeen AB24 3UE, United Kingdom 
\\ {\it email}: m.weiss@maths.abdn.ac.uk
\end{minipage} 

\vspace{1in} 
\emph{Abstract}: D.Mumford conjectured in \cite{Mumford83} that the 
rational cohomology of the stable moduli space of Riemann surfaces is a 
polynomial algebra generated by certain classes $\kappa_i$ of 
dimension $2i$. For the purpose of calculating rational cohomology, 
one may replace the stable moduli space of Riemann surfaces by 
$B\Gamma_{\infty}$, where $\Gamma_\infty$ is the group of 
isotopy classes of automorphisms of a smooth oriented connected 
surface of ``large'' genus. Tillmann's insight \cite{Tillmann97} 
that the plus construction makes $B\Gamma_{\infty}$ into an infinite 
loop space led to a stable homotopy version of Mumford's conjecture, 
stronger than the original \cite{MadsenTillmann01}. We prove the 
stronger version, relying on Harer's stability theorem~\cite{Harer85},
Vassiliev's theorem concerning spaces of functions with moderate 
singularities \cite{Vassiliev89}, \cite{Vassiliev94} and methods 
from homotopy theory.

\vspace{1in}
$^*)$ partially supported by American Institute of Mathematics
\newline
$^{**})$ partially supported by the Royal Society

\newpage 

\tableofcontents

\newpage

\section{Introduction: Results and methods}
\label{sec-intro}
 
\subsection{Main result}
Let $F=F_{g,b}$ be a smooth, compact, connected and oriented 
surface of genus $g>1$ with $b\ge 0$ boundary circles. Let $\sH(F)$ 
be the space of hyperbolic metrics on $F$ with geodesic boundary 
and such that each boundary circle has unit length. The topological 
group $\Diff(F)$ of orientation preserving diffeomorphisms 
$F\to F$ which restrict to the identity on the boundary 
acts on $\sH(F)$ by pulling back metrics. The orbit space 
\[  \sM(F) = \sH(F)\big/\Diff(F) \]
is the (hyperbolic model of the) moduli space of Riemann surfaces 
of topological type $F$. 

\medskip
The connected component $\Diff_1(F)$ of the identity acts 
freely on $\sH(F)$ with orbit space $\sT(F)$, the Teichm\"{u}ller 
space. The projection from $\sH(F)$ to $\sT(F)$ is a principal 
$\Diff_1$--bundle \cite{EarleEells69}, \cite{EarleSchatz70}. 
Since $\sH(F)$ is contractible and $\sT(F)\cong\RR^{6g-b+2b}$, 
the subgroup $\Diff_1(F)$ must be contractible. Hence the mapping 
class group $\Gamma_{g,b}=\pi_0\Diff(F)$ is homotopy equivalent to the 
full group $\Diff(F)$, and $B\Gamma_{g,b}\simeq B\Diff(F)$. 

When $b>0$ the action of $\Gamma_{g,b}$ on $\sT(F)$ is free so that 
$B\Gamma_{g,b}\simeq \sM(F)$. If $b=0$ the action of $\Gamma_g$ on 
$\sT(F)$ has finite isotropy groups and $\sM(F)$ has singularities. In
this case 
\[  B\Gamma_g \simeq (E\Gamma_g\times \sT(F))\big/\Gamma_g \]
and the projection $B\Gamma_g\to \sM(F)$ is only a rational 
homology equivalence. \newline 
For $b>0$, the standard homomorphisms 
\begin{equation*}
\label{eqn-Harerstability}  
\begin{array}{ccccccc}
\Gamma_{g,b} & \to &\Gamma_{g+1,b}\,, && 
\Gamma_{g,b} & \to &\Gamma_{g,b-1}
\end{array} \ttag 
\end{equation*}
yield maps of classifying spaces that induce isomorphisms in integral 
cohomology in degrees less than $g/2 - 1$ by the stability theorems 
of Harer \cite{Harer85} and Ivanov \cite{Ivanov89}. We let 
$B\Gamma_{\infty,b}$ denote the mapping telescope or homotopy 
colimit of  
\[ B\Gamma_{g,b}\lra B\Gamma_{g+1,b} \lra  B\Gamma_{g+2,b} \lra \cdots\,\,\,.
\]
Then $H^*(B\Gamma_{\infty,b};\ZZ)\cong H^*(B\Gamma_{g,b};\ZZ)$ 
for $*<g/2-1$, and in the same range the cohomology groups are 
independent of $b$. \newline
The mapping class groups $\Gamma_{g,b}$ are perfect for $g>1$,
so we may apply Quillen's plus construction to their classifying
spaces. By the above, the resulting homotopy type  
is independent of $b$ when $g=\infty$; we write
\[ B\Gamma_{\infty}^+ = B\Gamma_{\infty,b}^+\,. \] 
The main result from \cite{Tillmann97} asserts
that $\ZZ\times B\Gamma_{\infty}^+$ is an infinite 
loop space, so that homotopy classes of maps to it 
form the degree 0 part of a generalized cohomology theory. Our main 
theorem identifies this cohomology theory. 

\medskip
Let $G(d,n)$ denote the Grassmann manifold of oriented
$d$--dimensional subspaces of $\RR^{d+n}$, and let $U_{d,n}$ and 
$U_{d,n}^{\bot}$ be the two canonical vector bundles on $G(d,n)$ 
of dimension $d$ and $n$, respectively. The restriction 
\[  U_{d,n+1}^{\bot}|G(d,n) \]
is the direct sum of $U_{d,n}^{\bot}$ 
and a trivialized real line bundle. This yields an inclusion of their 
associated Thom spaces, 
\[  S^1\wedge \Th(U_{d,n}^{\bot}) \lra \Th(U_{d,n+1}^{\bot})\,, \]
and hence a sequence of maps (in fact cofibrations) 
\[ 
\cdots\to \Omega^{n+d}\Th(U_{d,n}^{\bot}) \to 
\Omega^{n+1+d}\Th(U_{d,n+1}^{\bot}) \to \cdots
\]
with colimit 
\begin{equation*}
\label{eqn-hVThomspectrum} 
\Omega^{\infty}\bV\,\,=\,\,\colim_n\, \Omega^{n+d}\Th(U_{d,n}^{\bot}).
\ttag
\end{equation*}
For $d=2$, the spaces $G(d,n)$ approximate the complex 
projective spaces, and 
\[ \Omega^{\infty}\bV \simeq \Omega^{\infty}\CC\bP^{\infty}_{-1}
:= \colim_n\,\, \Omega^{2n+2}\Th(L_n^{\bot}) \]
where $L_n^{\bot}$ is the complex $n$--plane bundle on $\CC
P^n$ which is complementary to the tautological line bundle $L_n$. \newline 
There is a map $\alpha_{\infty}$ from $\ZZ\times B\Gamma_{\infty}^+$
to $\Omega^{\infty}\CC\bP^{\infty}_{-1}$
constructed and examined in considerable detail in 
\cite{MadsenTillmann01}. Our main result is the following theorem 
conjectured in \cite{MadsenTillmann01}:

\begin{thm}
\label{thm-mainresult} The map 
$\,\alpha_{\infty}\co  \ZZ\times B\Gamma_{\infty}^+\lra 
\Omega^{\infty}\CC\bP^{\infty}_{-1}\,$
is a homotopy equivalence. 
\end{thm}

\medskip
Since $\alpha_{\infty}$ is an infinite loop map by
\cite{MadsenTillmann01}, the theorem identifies the generalized 
cohomology theory determined by $\ZZ\times B\Gamma_{\infty}^+$
to be the one associated with the spectrum $\CC\bP^{\infty}_{-1}$. 
To see that theorem~\ref{thm-mainresult} verifies Mumford's 
conjecture we consider the homotopy fibration sequence of 
\cite{Ravenel84}, 
\begin{equation*}
\label{eqn-ravenel}
\CD \Omega^{\infty}\CC\bP^{\infty}_{-1} @>\,\,\omega\,\,>>
\Omega^{\infty} S^{\infty}(\CC P^{\infty}_+) 
@>\,\,\partial\,\, >> \Omega^{\infty+1} S^{\infty} \endCD \ttag
\end{equation*}
where the subscript $+$ denotes an added disjoint base point.
The homotopy groups of $\Omega^{\infty+1} S^{\infty}$ are equal to 
the stable homotopy groups of spheres, up to a shift of one, and 
are therefore finite. Thus $H^*(\omega;\QQ)$ is an isomorphism. 
The canonical complex line bundle over $\CC P^{\infty}$, considered 
as a map from $\CC P^{\infty}$ to $\{1\}\times B\Un$, induces via Bott 
periodicity a map 
\[ L\co \Omega^{\infty} S^{\infty}(\CC\bP^{\infty}_+)
\lra \ZZ\times B\Un, \]
and $H^*(L;\QQ)$ is an isomorphism. Thus we have isomorphisms 
\[ H^*(\ZZ\times B\Gamma^+_{\infty};\QQ)\,\,\cong\,\, 
H^*(\Omega^{\infty}\CC\bP^{\infty}_{-1};\QQ) \,\,\cong\,\, H^*(\ZZ\times
B\Un;\QQ)\,.
\]
Since Quillen's plus construction leaves cohomology undisturbed this
yields Mumford's conjecture:
\[  H^*(B\Gamma_{\infty};\QQ) \,\,\cong \,\,H^*(B\Un;\QQ) \,\,\cong\,\, 
\QQ[\kappa_1,\kappa_2,\dots]\,. \]
Miller, Morita and Mumford \cite{Miller86}, \cite{Morita84}, 
\cite{Morita87}, \cite{Mumford83} defined the classes 
$\kappa_i\in H^{2i}(B\Gamma_{\infty};\QQ)$ 
by integration (Umkehr) of the $(i+1)$--th power of the tangential 
Euler class in the universal smooth $F_{g,b}$--bundles. In the 
above setting
$\kappa_i=\alpha_{\infty}^*L^*(i!\,\textrm{ch}_i)$. \newline 
We finally remark that the cohomology 
$H^*(\Omega^{\infty}\CC\bP^{\infty}_{-1};\FF_p)$ has been calculated 
in \cite{Galatius02} for all primes $p$. The result is 
quite complicated.

\subsection{A geometric formulation}

\medskip
Let us first consider smooth \emph{proper} maps 
$q\co M^{d+n}\to X^n$ of smooth manifolds without boundary, of 
fixed relative dimension $d\ge 0$, and equipped with an orientation 
of the (stable) relative tangent bundle $TM-q^*TX$. Two such maps 
$q_0\co M_0\to X$ and $q_1\co M_1\to X$ are \emph{concordant} (traditionally,
\emph{cobordant}) if there exists a similar map 
$q_{\RR}\co W^{d+n+1}\to X\times\RR$ transverse to $X\times\{0\}$ and 
$X\times\{1\}$, and such that the inverse images of $X\times\{0\}$ and 
$X\times\{1\}$ are isomorphic to $q_0$ and $q_1$ (with all the relevant 
vector bundle data), respectively. 
The Pontryagin--Thom theory, cf.\ particularly \cite{Quillen71}, 
equates the set of concordance classes 
of such maps over fixed $X$ with the set of homotopy classes of maps 
from $X$ into the degree $-d$ term of the universal Thom spectrum, 
\[ \Omega^{\infty+d}{\bf M}\SOr = \colim_n\,\,
\Omega^{n+d}\Th(U_{n,\infty})\,. \] 
The geometric reformulation of 
theorem~\ref{thm-mainresult} is similar in spirit.

We consider smooth proper maps $q\co M^{d+n}\to X^n$ much as before, 
together with a vector bundle epimorphism $\delta q$ from $TM\times\RR^i$ to
$q^*TX\times\RR^i$, where $i\gg 0$, and an orientation of the 
d--dimensional kernel bundle of $\delta q$. (Note that $\delta q$ is 
not required to agree with $dq$, the differential of $q$.) Again, 
the Pontryagin--Thom theory equates the set of concordance classes 
of such pairs $(q,\delta q)$ over fixed $X$ with the set of homotopy 
classes of maps 
\[  X\lra \Omega^{\infty}\bV\,, \]
with $\Omega^{\infty}\bV$ as in~(\ref{eqn-hVThomspectrum}). 
For a pair $(q,\delta q)$ as above which is integrable, 
$\delta q=dq$, the map $q$ is a proper submersion with target $X$ and hence 
a bundle of smooth closed $d$--manifolds on $X$ by Ehresmann's fibration 
lemma \cite[8.12]{BroeckerJaenich73}. Thus the set 
of concordance classes of such integrable pairs over a fixed $X$ is
in natural bijection with the set of homotopy classes of maps
\[  X \lra \coprod B\Diff(F^d) \] 
where the disjoint union runs over a set of representatives of 
the diffeomorphism classes of closed, smooth and oriented
$d$--manifolds. Comparing these two classification results we obtain 
a map 
\begin{equation*}
\label{eqn-elementaryalphamap}
\alpha\co \coprod B\Diff(F^d) \lra \bV 
\end{equation*}
which for $d=2$ is closely related to the map 
$\alpha_{\infty}$ of theorem~\ref{thm-mainresult}. The map 
$\alpha$ is not a homotopy equivalence (which is why we replace 
it by $\alpha_{\infty}$ when $d=2$). However, using submersion 
theory we can refine our geometric understanding of homotopy classes of 
maps to $\bV$ and our understanding of $\alpha$. \newline 
We suppose for simplicity that $X$ is closed. As explained 
above, a homotopy class of maps
from $X$ to $\bV$ can be represented by a pair $(q,\delta q)$ 
with a proper $q\co M\to X$, a vector 
bundle epimorphism $\delta q\co TM\times\RR^i\to
q^*TX\times\RR^i$ and an orientation on $\ker(\delta q)$. We set 
\[ E= M\times\RR \]
and let $\bar q\co E\to X$ be given by $\bar q(x,t)=q(x)$. The epimorphism 
$\delta q$ determines an epimorphism $\delta\bar q\co TE\times\RR^i\to
{\bar q\,}^*TX\times \RR^i$. In fact, obstruction theory shows that 
we can take $i=0$, so we write 
$\delta\bar q\co TE\to {\bar q\,}^*TX$. Since $E$ 
is an open manifold, the submersion 
theorem of Phillips \cite{Phillips67}, \cite{Haefliger71}, 
\cite{Gromov71} applies, showing that
the pair $(\bar q,\delta \bar q)$ is homotopic through vector bundle
surjections to a pair $(\pi,d\pi)$ 
consisting of a submersion $\pi\co E\to X$ and 
its differential $d\pi\co TE\to \pi^*TX$. Let $f\co E\to \RR$ 
be the projection. This is obviously proper; consequently 
$(\pi,f)\co E\to X\times\RR$ is proper. \newline 
The vertical tangent bundle $T^{\pi}E=\ker(d\pi)$ of $\pi$ 
is identified with $\ker(\delta p)\cong \ker(\delta q)\times T\RR$, 
so has a trivial line bundle factor. 
Let $\delta f$ be the projection to that factor. In terms 
of the vertical $1$--jet bundle  
\[  p^1_{\pi}\co  J^1_{\pi}(E,\RR) \lra E \]
whose fiber at $z\in E$ consists of all affine maps from the 
vertical tangent space $(T^{\pi}E)_z$ to $\RR$, the pair $(f,\delta f)$ 
amounts to a section $\hat f$ of 
$p^1_{\pi}$ such that $\hat f(z)\co (T^{\pi}E)_z\to \RR$ is surjective 
for every $z\in E$.  

\medskip
We introduce the notation $h\cV(X)$ for the set of pairs 
$(\pi,\hat f)$, where $\pi\co E\to X$ is a smooth submersion with 
$(d+1)$-dimensional oriented fibers and $\hat f\co E\to
J^1_{\pi}(E,\RR)$ is a section of
$p^1_{\pi}$ with underlying map 
$f\co E\to\RR$, subject to two conditions: for each $z\in E$ 
the affine map $\hat f(z)\co (T^{\pi}E)_z\to\RR$ is surjective, and 
$(\pi,f)\co E\to X\times\RR$ is proper. \newline 
Concordance defines an equivalence relation on $h\cV(X)$. Let 
$h\cV[X]$ be the set of equivalence classes. The 
arguments above lead to a natural bijection
\begin{equation*}
\label{eqn-firstbijection}  
h\cV[X] \,\, \cong \,\,[X, \Omega^{\infty}\bV] \,. 
\ttag
\end{equation*} 
We similarly define $\cV(X)$ as the set of pairs $(\pi,f)$ 
where $\pi\co E\to X$ is a smooth submersion as before and 
$f\co E\to\RR$ is a smooth function, subject to two conditions: 
the restriction of $f$ to any fiber of $\pi$ is regular 
($=$ nonsingular), and $(\pi,f)\co E\to X\times\RR$ is proper.
Let $\cV[X]$ be the correponding set of concordance classes. 
Since elements of $\cV(X)$ are bundles of closed oriented $d$--manifolds over
$X\times\RR$, we have a natural bijection 
\[ \cV[X] \,\,\cong\,\, [X, \coprod B\Diff(F^d)]. \]
On the other hand an element $(\pi,f)\in \cV(X)$ with $\pi\co E\to X$ 
determines a section $j^1_{\pi}f$ of the projection $J^1_{\pi}E\to E$ 
by fiberwise $1$--jet prolongation. The map 
\begin{equation*}
\label{eqn-thejetprolongation}
\begin{array}{cccccc}
\cV(X) & \lra & h\cV(X) \, ; &  
(\pi,f) & \mapsto & (\pi,j^1_{\pi}f)
\end{array}
\ttag
\end{equation*} 
respects the concordance relation and so induces a map 
$\cV[X]\to h\cV[X]$, which corresponds to $\alpha$ 
in~(\ref{eqn-elementaryalphamap}). 

\subsection{Outline of proof} 
The main tool is a special case of the celebrated ``first main 
theorem'' of V.A. Vassiliev \cite{Vassiliev94}, \cite{Vassiliev89}
which can be used to approximate~(\ref{eqn-thejetprolongation}). 
\newline 
We fix $d\ge 0$ as above. For smooth $X$ without boundary we 
enlarge the set $\cV(X)$ to the set $\cW(X)$ consisting of pairs
$(\pi,f)$ with $\pi$ as before but with $f\co E\to \RR$ a 
fiberwise Morse function rather than a fiberwise regular function. 
We keep the condition that the combined map 
$(\pi,f)\co E\to X\times\RR$ be proper. There is a 
similar enlargement of $h\cV(X)$ to a set $h\cW(X)$. An element 
of $h\cW(X)$ is a pair $(\pi,\hat f)$ where $\hat f$ is a section 
of ``Morse type'' of the fiberwise $2$--jet bundle 
$J^2_{\pi}E\to E$ with an underlying map $f$ such that 
$(\pi,f)\co E\to X\times \RR$ is proper. In analogy 
with~(\ref{eqn-thejetprolongation}), we have the $2$--jet prolongation 
map 
\begin{equation*}
\label{eqn-theotherjetprolongation}
\begin{array}{cccccc}
\cW(X) & \lra & h\cW(X) \, ; &  
(\pi,f) & \mapsto & (\pi,j^2_{\pi}f)\,.
\end{array}
\ttag
\end{equation*}  
Dividing out by the concordance relation we get representable functors: 
\begin{equation*}
\label{eqn-repfunctors}
\begin{array}{ccccccc}
\cW[X] & \cong & [X, |\cW|\,]\,, & &  h\cW[X] & \cong & [X, |h\cW| \,]
\end{array} 
\ttag
\end{equation*}
and~(\ref{eqn-theotherjetprolongation}) induces a map 
$j^2_{\pi}\co |\cW|\to |h\cW|$. Vassiliev's first main theorem is
a main ingredient in our proof (in section~\ref{sec-Vassiliev}) of 

\begin{thm} 
\label{thm-middlecolumn} The jet prolongation map $|\cW| \to |h\cW|$ 
is a homotopy equivalence. 
\end{thm} 

There is a commutative square 
\begin{equation*}
\label{eqn-leftstrategic}
\xymatrix{
|\cV| \ar[r] \ar[d] & |\cW| \ar[d] \\
|h\cV| \ar[r] & |h\cW|\,.
}
\ttag
\end{equation*}
We need information about the horizontal maps. This involves 
introducing ``local'' versions $\cW_{\loc}(X)$ and 
$h\cW_{\loc}(X)$ where we focus on the behavior of functions 
$f$ and jet bundle sections $\hat f$ near the fiberwise 
singularity set: 
\[ 
\begin{array}{ccl}  
\Sigma(\pi,f) & = & \{ z\in E\mid df_z=0 \textrm{ on }
(T^{\pi}E)_z\}\,, \\ 
\Sigma(\pi,\hat f) & = & 
\{ z\in E\mid \textrm{ linear part of $\hat f(z)$ vanishes} \}. 
\end{array}
\]
Technically the localization is easiest to achieve as follows. 
Elements of $\cW_{\loc}(X)$ are defined like elements $(\pi,f)$ 
of $\cW(X)$, but we relax the condition that 
$(\pi,f)\co E\to X\times\RR$ be proper to the condition that its 
restriction to $\Sigma(\pi,f)$ be proper. The definition of 
$h\cW_{\loc}(X)$ is similar, and we obtain spaces $|\cW_{\loc}|$ 
and $|h\cW_{\loc}|$ which represent the corresponding 
concordance classes, together with a commutative diagram 
\medskip
\begin{equation*}
\label{eqn-strategicdiagram} 
\xymatrix@M=8pt@W=8pt{
|\cV| \ar[r] \ar[d]_{j^2_{\pi}} & |\cW| \ar[r] \ar[d]_{j^2_{\pi}}
& |\cW_{\loc}| \ar[d]_{j^2_{\pi}} \\
|h\cV| \ar[r] & |h\cW| \ar[r] & |h\cW_{\loc}|.
}  \ttag
\end{equation*}
The next two theorems are proved in section~\ref{sec-bottomrow}. 
They are much easier than theorem~\ref{thm-middlecolumn}. 

\begin{thm} 
\label{thm-twolocs} The jet prolongation map 
$|\cW_{\loc}|\to |h\cW_{\loc}|$ is a homotopy equivalence. 
\end{thm}

\begin{thm} 
\label{thm-bottomrow} The maps 
$|h\cV|\to |h\cW|\to |h\cW_{\loc}|$ define a homotopy 
fibration sequence of infinite loop spaces.
\end{thm}

The spaces $|h\cW|$ and $|h\cW_{\loc}|$ are, like 
$|h\cV|=\Omega^{\infty}\bV$, colimits of certain iterated loop 
spaces of Thom spaces. Their homology can be approached 
by standard methods from algebraic topology. \newline
The three theorems above are valid for any choice of $d\ge 0$. 
This is not the case for the final result that goes into the proof 
of theorem~\ref{thm-mainresult}, although many of the arguments
leading to it are valid in general. 

\begin{thm}
\label{thm-toprow} For $d=2$, the homotopy fiber of $|\cW|\to
|\cW_{\loc}|$ is the space $\ZZ\times B\Gamma_{\infty}^+$. 
\end{thm} 

In conjunction with the previous three theorems this proves 
theorem~\ref{thm-mainresult}: 
\[ 
\begin{array}{ccccccc}
\ZZ\times B\Gamma_{\infty}^+ & \simeq & |h\cV| & \simeq &
\Omega^{\infty}\bV & \simeq & \Omega^{\infty} \CC\bP^{\infty}_{-1}\,. 
\end{array}
\] 

\medskip
The proof of theorem~\ref{thm-toprow} is technically the most 
demanding part of the paper. It rests on compatible stratifications 
of $|\cW|$ and $|h\cW|$, or more precisely on homotopy 
colimit decompositions 
\begin{equation*}
\label{eqn-introdecomp}
\begin{array}{cccccc}
|\cW| & \simeq & \hocolim_R\,\,|\cW_R|\,, & 
|\cW_{\loc}| &  \simeq & \hocolim_R\,\,|\cW_{\loc,R}|
\end{array} \ttag
\end{equation*}
where $R$ runs through the objects of a certain category 
of finite sets. The spaces 
$|\cW_R|$ and $|\cW_{\loc,R}|$ classify certain bundle theories 
$\cW_R(X)$ and $\cW_{\loc,R}(X)$. The proof of~(\ref{eqn-introdecomp})
is given in section~\ref{sec-strat}, and is valid for all $d\ge 0$. 
(Elements of $\cW_R(X)$ are smooth fiber bundles $M^{n+d}\to X^n$ 
equipped with extra fiberwise ``surgery data''. The maps 
$\cW_S(X)\to \cW_R(X)$ induced contravariantly by morphisms $R\to S$ in the 
indexing category involve fiberwise surgeries on some of these data.) \newline
The homotopy fiber of 
$|\cW_R|\to |\cW_{\loc,R}|$ is a classifying 
space for smooth fiber bundles $M^{n+d}\to X^n$ with $d$--dimensional 
oriented fibers $F^d$, each fiber having its boundary identified with a 
disjoint union 
 \[  \coprod_{r\in R} S^{\mu_r}\times S^{d-\mu_r-1}  \]
where $\mu_r$ depends on $r\in R$. The fibers $F^d$ need not 
be connected, but in section~\ref{sec-surgery}
we introduce a modification $\cW_{c,R}(X)$ of $\cW_R(X)$ to 
enforce this additional property, 
keeping~(\ref{eqn-introdecomp}) almost intact. 
Again this works for all $d\ge 0$.

\medskip
When $d=2$ the homotopy fiber of $|\cW_{c,R}|\to |\cW_{\loc,R}|$ 
becomes homotopy equivalent to $\coprod_g B\Gamma_{g,2|R|}$. 
A second modification of~(\ref{eqn-introdecomp}) which we undertake 
in section~\ref{sec-stab} allows us to replace this by 
$\ZZ\times B\Gamma_{\infty,2|R|+1}$, functorially in $R$. 
It follows directly from Harer's theorem that these homotopy 
fibers are ``independent'' of $R$ up to homology equivalences.  
Using an argument from \cite{McDuffSegal76} and \cite{Tillmann97} we  
conclude that the inclusion of any of these homotopy fibers $\ZZ\times
B\Gamma_{\infty,2|R|+1}$ into the homotopy fiber of 
$|\cW|\to |\cW_{\loc}|$ is a homology equivalence. This proves
theorem~\ref{thm-toprow}.

\medskip
The paper is set up in such a way that it proves analogues 
of theorem~\ref{thm-mainresult} for other classes of surfaces, 
provided that Harer type stability results have been 
established. This includes for example spin surfaces. 
See also \cite{Galatius04}.

\bigskip
\section{Families, sheaves and their representing spaces}
\label{sec-bundletheories}
\subsection{Language}
\label{subsec-language}
We will be interested in families of
smooth manifolds, parametrized by other smooth manifolds. In order 
to formalize pullback constructions and gluing properties for such
families, we need the language of sheaves. Let $\sX$ be the category 
of smooth manifolds (without boundary, with a countable base) and 
smooth maps. 

\begin{dfn}
\label{dfn-sheaf}
{\rm A \emph{sheaf} on $\sX$ is a contravariant functor $\cF$ 
from $\sX$ to the category of sets with the following property. 
For every open covering $\{U_i|i\in\Lambda\}$ of some $X$ in $\sX$, 
and every collection $(s_i\in F(U_i))_i$ satisfying 
$s_i|U_i\cap U_j=s_j|U_i\cap U_j$ for all $i,j\in \Lambda$, there is a 
unique $s\in F(X)$ such that $s|U_i=s_i$ for all $i\in\Lambda$. }
\end{dfn}

In definition~\ref{dfn-sheaf},
we do not insist that all of the $U_i$ be nonempty. 
Consequently $\cF(\emptyset)$ must be a singleton. For a 
disjoint union $X=X_1\sqcup X_2$, the restrictions 
give a bijection $\cF(X)\cong \cF(X_1)\times \cF(X_2)$. 
Consequently $\cF$ is determined up to unique natural bijections by 
its behavior on connected nonempty objects $X$ of $\sX$.

\medskip
For the sheaves $\cF$ that we will be considering, an element of 
$\cF(X)$ is typically a family of 
manifolds parametrized by $X$ and with some additional structure. 
In this situation there is usually a sensible concept of 
\emph{isomorphism} between elements of $\cF(X)$, so that 
there might be a temptation to regard $\cF(X)$ as a groupoid. 
We do not include these isomorphisms in our definition of 
$\cF(X)$, however, and we do not suggest that elements of $X$ 
should be confused with the corresponding isomorphism classes
(since this would destroy the sheaf property). This paper is not 
about ``stacks''. All the same, we must 
ensure that our pullback and gluing constructions are well 
defined (and not just up to some sensible notion of isomorphism
which we would rather avoid). This forces us to introduce the following 
purely set--theoretic concept. We fix, once and for all, a set $Z$ 
whose cardinality is at least that of $\RR$.  

\begin{dfn}
\label{dfn-graphicmaps}
 {\rm A map of sets $S\to T$ is \emph{graphic}
if it is a restriction of the projection $Z\times T\to T$. 
In particular, each graphic map with target $T$ is determined 
by its source, which is a subset $S$ of $Z\times T$.}
\end{dfn}

Clearly, a graphic map $f$ with target $T$ is equivalent to a 
map from $T$ to the power set $P(Z)$ of $Z$, which we may 
call the adjoint of $f$. Pullbacks of graphic maps are now easy 
to define: If $g\co T_1\to T_2$ is any map 
and $f\co S\to T_2$ 
is a graphic map with adjoint $f^a\co T\to P(Z)$, then the pullback
$g^*f\co g^*S\to T_1$ is, by definition, the graphic map with adjoint  
equal to the composition 
\begin{equation*}
\label{eqn-basechange}
\xymatrix{
T_1 \ar[r]^g & T_2 \ar[r]^{f^a} & P(Z).
} \ttag
\end{equation*}
If $g$ is an identity, then $g^*S=S$ and $g^*f=f$;  
if $g$ is a composition, $g=g_2g_1$, then
$g^*S={g_1}^*{g_2}^*S$ and $g^*f={g_1}^*{g_2}^*f$. 
Thus, with the above definitions, \emph{base change} is associative. 

\begin{dfn}
\label{dfn-concordantsections}
{\rm Let $\pr\co X\times\RR\to X$ be the 
projection. Two elements $s_0,s_1$ of $\cF(X)$ are \emph{concordant}
if there exist $s\in \cF(X\times\RR)$  
which agrees with $\pr^*s_0$ on an open neighborhood 
of $X\times\,]-\infty,0]$ in $X\times\RR$, and with 
$\pr^*s_1$ on an open neighborhood of $X\times[1,+\infty[\,$
in $X\times\RR$. The element $s$ is then called a concordance 
from $s_0$ to $s_1$.}  
\end{dfn}

It is not hard to show that ``being concordant'' is an 
equivalence relation on the set $\cF(X)$, for every $X$. 
We denote the set of equivalence classes by $\cF[X]$. Then  
$X\mapsto \cF[X]$ is still a contravariant functor on $\sX$. It is 
practically never a sheaf, but it is
\emph{representable} in the following weak sense. There exists a 
space, denoted by $|\cF|$, such that homotopy classes of maps 
from a smooth $X$ to $|\cF|$ are in natural bijection with the elements 
of $\cF[X]$. This follows from very general principles
expressed in Brown's representation theorem \cite{Brown62}. An 
explicit and more functorial construction of $|\cF|$ will be described  
later. To us, $|\cF|$ is more important than $\cF$ itself.
We define $\cF$ in order to pin down $|\cF|$. 

Elements in $\cF(X)$ can usually be regarded as \emph{families} of elements 
in $\cF(\pt)$, parametrized by the manifold $X$. The space $|\cF|$ 
should be thought of as a space which classifies families 
of elements in $\cF(\pt)$.

\subsection{Families with analytic data} 
\label{subsec-analyticsheaves}  
Let $E$ be a smooth manifold, without boundary for now, and $\pi\co
E\to X$ a smooth map to an object of $\sX$. 
The map $\pi$ is a \emph{submersion} if its 
differentials $TE_z \to  TX_{\pi(z)}$ for $z\in E$
are all surjective. In that case, by the implicit function theorem, 
each fiber $E_x=\pi^{-1}(x)$ for $x\in X$ is a smooth submanifold 
of $E$, of codimension equal to $\dim(X)$. We remark that 
a submersion need not be surjective and a surjective submersion 
need not be a bundle.  However, a 
\emph{proper} smooth map $\pi\co E\to X$ which is a submersion 
is automatically a smooth fiber bundle by 
Ehresmann's fibration lemma \cite[thm. 8.12]{BroeckerJaenich73}. 

In this paper, when we informally mention a \emph{family} of 
smooth manifolds parametrized by some $X$ in $\sX$, what we 
typically mean is a submersion $\pi\co E\to X$. The members of the 
family are then the fibers $E_x$ of $\pi$. The \emph{vertical 
tangent bundle} of such a family is the vector bundle 
$T^{\pi}E\to E$ whose fiber at $z\in E$ is the kernel 
of the differential $d\pi\co TE_z\to TX_{\pi(z)}$.

To have a fairly general notion of orientation as well,
we fix a space $\Theta$ with a right action of the 
infinite general linear group over the real numbers:
 $\Theta\times\GL \to \Theta$. 
For an $n$--dimensional vector bundle $W\to B$ let  
$\textrm{Fr}(W)$ be the frame bundle, 
which we regard as a principal $\GL(n)$--bundle on $B$ with $\GL(n)$ acting 
on the right.

\begin{dfn}
\label{dfn-Thetaorientation}
 {\rm By 
a \emph{$\Theta$--orientation} of $W$ we mean 
a section of the associated bundle 
$(\textrm{Fr}(W)\times\Theta)/\GL(n)\lra B$.}
\end{dfn}

This includes a definition of a $\Theta$--orientation on a 
finite dimensional real vector \emph{space}, because a 
vector space is a vector bundle over a point.

\begin{exm}
\label{exm-Theta}
{\rm If $\Theta$ is a single point, then every vector 
bundle has a unique $\Theta$--orienta\-tion. If $\Theta$ is 
$\pi_0(\GL)$ with the action of $\GL$ by translation, then a 
$\Theta$--orientation of a vector bundle is simply an orientation. 
(This choice of $\Theta$ is the one that will be needed in the 
proof of the Mumford conjecture.) If $\Theta$ is $\pi_0(\GL)\times Y$ 
for a fixed space $Y$, where $\GL$ acts by translation on $\pi_0(\GL)$
and trivially on the factor $Y$, then a $\Theta$--orientation on 
a vector bundle $W\to B$ is an orientation on $W$ together 
with a map $B\to Y$. \newline
Let $\widetilde\SL(n)$ be the universal cover of the special 
linear group $\SL(n)$. If $\Theta=\colim_n\Theta_n$ where 
$\Theta_n$ is the pullback of 
\[
\xymatrix{
E\GL(n) \ar[r] & B\GL(n) & \ar[l] B\widetilde\SL(n)\,,
}
\]
then a $\Theta$--orientation on a vector bundle $W$ amounts to a 
spin structure on $W$. Here $E\GL(n)$ can be taken as the total space 
of the frame bundle associated with the universal $n$--dimensional 
vector bundle on $B\GL(n)$.
}
\end{exm}

\medskip
We also fix an integer $d\ge 0$. (For the proof of the 
Mumford conjecture, $d=2$ is the right choice.) The data $\Theta$ 
and $d$ will remain with us, fixed but unspecified, throughout the
paper, except for section~\ref{sec-stab} where we specialize to 
$d=2$ and $\Theta=\pi_0\GL$.

\begin{dfn}
\label{dfn-Vsheaf} {\rm For 
$X$ in $\sX$, let $\cV(X)$ be the set of 
pairs $(\pi,f)$ where $\pi\co E\to X$ is a
graphic submersion of fiber dimension $d+1$, 
with a $\Theta$--orientation of its vertical tangent 
bundle, and $f\co E\to \RR$ is a smooth map, subject to 
the following conditions. 
\begin{description}
\item[{\bf (i)}] The map $(\pi,f)\co E\to X\times\RR$ is proper.
\item[{\bf (ii)}] The map $f$ is fiberwise nonsingular, i.e., 
the restriction of $f$ to any fiber $E_x$ of $\pi$ is a 
nonsingular map.
\end{description} 
}
\end{dfn}

\medskip
For $(\pi,f)\in \cV(X)$ with $\pi\co E\to X$, the 
map $z\mapsto(\pi(z),f(z))$ from $E$ to $X\times\RR$ is
a proper submersion and therefore a smooth bundle with
$d$--dimensional fibers. The $\Theta$--orientation on the 
vertical tangent bundle of $\pi$ is equivalent 
to a $\Theta$--orientation on the vertical tangent bundle of 
$(\pi,f)\co E\to X\times\RR$, since $T^{\pi}E\cong T^{(\pi,f)}E\times\RR$. 
Consequently~\ref{dfn-Vsheaf} is another 
way of saying that an element of $\cV(X)$ is a bundle of smooth 
closed $d$--manifolds on $X\times\RR$ with a 
$\Theta$--orientation of its vertical tangent bundle. 
We prefer the formulation given in 
definition~\ref{dfn-Vsheaf} because it is easier to vary and
generalize, as illustrated by our next definition.  

\begin{dfn}
\label{dfn-Wsheaf} {\rm For 
$X$ in $\sX$, let $\cW(X)$ be the set of 
pairs $(\pi,f)$ as in definition~\ref{dfn-Vsheaf}, subject 
to condition (i) as before, but with condition (ii)
replaced by the weaker condition 
\begin{description}
\item[{\bf (iia)}] the map $f$ is fiberwise Morse.
\end{description}
}
\end{dfn}

\medskip
Recall that a smooth function $N\to \RR$ is a Morse function 
precisely if its differential, viewed as a smooth section 
of the cotangent bundle $TN^*\to N$, is transverse to the 
zero section \cite[II\S6]{GolubitskyGuillemin73}. This 
extends to families of smooth manifolds and Morse functions. 
In other words, if $\pi\co E\to X$ is a smooth 
submersion and $f\co E\to \RR$ is any smooth map, then $f$ 
is fiberwise Morse if and only if the fiberwise differential 
of $f$, a section of the vertical cotangent bundle 
$T^{\pi}E^*$ on $E$, is transverse to the zero section. 
This has the following
consequence for the fiberwise 
singularity set $\Sigma(\pi,f)\subset E$ of $f$. 

\begin{lem} 
\label{lem-singprojetale} Suppose that $f\co E\to \RR$
is fiberwise Morse. Then $\Sigma(\pi,f)$ is a smooth 
submanifold of $E$ and the 
restriction of $\pi$ to $\Sigma(\pi,f)$ is a local 
diffeomorphism, alias \emph{\'{e}tale} map, from 
$\Sigma(\pi,f)$ to $X$.
\end{lem} 

\proof The fiberwise differential 
viewed as a section of the vertical cotangent bundle 
is transverse to the zero section. In 
particular $\Sigma=\Sigma(\pi,f)$ is a submanifold of $E$, of the same
dimension as $X$. But moreover, the fiberwise Morse condition 
implies that for each $z\in\Sigma$, the tangent space $T\Sigma_z$ 
has trivial intersection in $TE_z$ with the vertical tangent space 
$T^{\pi}E_z$. This means that $\Sigma$ is transverse 
to each fiber of $\pi$, and also that the differential 
of $\pi|\Sigma$ at any point $z$ of $\Sigma$ is an invertible 
linear map $T\Sigma_z\to TX_{\pi(z)}$, and consequently that
$\pi|\Sigma$ is a local diffeomorphism. \qed

\begin{dfn} 
\label{dfn-Wlocsheaf} 
{\rm For $X$ in $\sX$ let $\cW_{\loc}(X)$ be the set of
pairs $(\pi,f)$, as in definition~\ref{dfn-Vsheaf}, 
but replacing conditions (i) and (ii) by 
\begin{description}
\item[{\bf (ia)}] the map $\Sigma(\pi,f)\to X\times\RR$ defined by 
$z\mapsto(\pi(z),f(z))$ is proper, 
\item[{\bf (iia)}] $f$ is fiberwise Morse.
\end{description}
}
\end{dfn}

\subsection{Families with formal--analytic data} 
\label{subsec-mocksheaves}

Let $E$ be a smooth manifold and $p^k\co J^k(E,\RR)\to E$ the 
$k$-jet bundle, where $k\ge 0$. Its fiber $J^k(E,\RR)_z$
at $z\in E$ consists of equivalence 
classes of smooth map germs $f\co (E,z)\to \RR$, with $f$ equivalent
to $g$ if the $k$-th Taylor expansions of $f$ and $g$ agree at $z$
(in any local coordinates near $z$). The elements of $J^k(E,\RR)$ are 
called \emph{$k$--jets} of maps from $E$ to $\RR$. 
The $k$-jet bundle $p^k\co J^k(E,\RR)\to E$ is a vector bundle.

Let $u\co TE_z\to E$ be any exponential map at $z$, that is, 
a smooth map such that $u(0)=z$ and the differential at $0$ 
is the identity $TE_z\to TE_z$. Then every jet $t\in J^k(E,\RR)_z$
can be represented by a unique germ $(E,z)\to\RR$ whose 
composition with $u$ is the germ at $0$ of a polynomial function 
$t_u$ of degree $\le k$ on the vector space $TE_z$. The 
constant part (a real number) and the linear part 
(a linear map $TE_z\to \RR$) of $t_u$ do not depend on $u$. We call them
the constant and linear part of $t$, respectively. If 
the linear part of $t$ vanishes, then the quadratic part of $t_u$,
which is a quadratic map $TE_z\to \RR$, is again independent of $u$.
We then call it the quadratic part of $t$. 

\begin{dfn} {\rm A jet $t\in J^k(E,\RR)$ is \emph{nonsingular} 
(assuming $k\ge 1$) if its 
linear part is nonzero. The jet $t$ is 
\emph{Morse} (assuming $k\ge2$) if it has a nonzero linear part or, 
failing that, a nondegenerate quadratic part.}
\end{dfn}

A smooth function $f\co E\to \RR$ induces a smooth section 
$j^kf$ of $p^k$, which we call the 
\emph{$k$-jet prolongation}
of $f$, following e.g. Hirsch \cite{Hirsch76}. (Some  
writers choose to call it the \emph{$k$-jet} of $f$, which can be  
confusing.)  
Not every smooth section of $p^k$ has this form. 
Sections of the form $j^kf$ are called \emph{integrable}. 
Thus a smooth section of $p^k$ is integrable if and only if it 
agrees with the $k$-jet prolongation of its underlying smooth map 
$f\co E\to \RR$.

We need a fiberwise version $J^k_{\pi}(E,\RR)$ of $J^k(E,\RR)$,
fiberwise with respect to a submersion $\pi\co E^{j+r}\to X^j$
with fibers $E_x$ for $x\in X$. 
In a neighborhood of any $z\in E$ 
we may choose local coordinates $\RR^j\times\RR^r$  so that 
$\pi$ becomes the projection onto $\RR^j$ and $z=(0,0)$. Two smooth map
germs $f,g\co (E,z)\to \RR$ define the same element of
$J^k_{\pi}(E,\RR)_z$ if their $k$-th Taylor expansions in the $\RR^r$ 
coordinates agree at $(0,0)$. Thus $J^k_{\pi}(E,\RR)_z$ is a quotient 
of $J^k(E,\RR)_z$ and $J^k_{\pi}(E,\RR)_z$ is identified
with $J^k(E_{\pi(z)},\RR)$.
There is a short exact sequence of vector bundles on $E$, 
\[ \pi^*J^k(X,\RR) \lra J^k(E,\RR) \lra J^k_{\pi}(E,\RR) . \]
Sections of the bundle projection $p^k_{\pi}\co J^k_{\pi}(E,\RR)\to E$
will be denoted $\hat f$, $\hat g$, ..., and their underlying
functions from $E$ to $\RR$ by the corresponding letters $f$, $g$, 
and so on. Such a section 
$\hat f$ is \emph{nonsingular}, resp. \emph{Morse},
if $\hat f(z)$, viewed as an element of $J^k(E_{\pi(z)},\RR)$, 
is nonsingular, resp. Morse, for all $z\in E$.

\begin{dfn} {\rm The fiberwise 
singularity set $\Sigma(\pi,\hat f)$ is the set 
of all $z\in E$ where $\hat f(z)$ is singular (assuming $k\ge1$).
Equivalently, 
\[ \Sigma(\pi,\hat f) = \hat f^{-1}(\Sigma_{\pi}(E,\RR))\,, \]
where $\Sigma_{\pi}(E,\RR)\subset J^2_{\pi}(E,\RR)$ 
is the submanifold consisting of the singular jets, i.e., 
those with vanishing linear part.
} 
\end{dfn}

Again, any smooth function $f\co E\to \RR$ induces a smooth section 
$j^k_{\pi}f$ of $p^k_{\pi}$, which we call the 
\emph{fiberwise $k$-jet prolongation}
of $f$. The sections of the form $j^k_{\pi}f$
are called \emph{integrable}. If $k\ge 1$ and 
$\hat f$ is integrable with 
$\hat f=j^k_{\pi}f$, then 
\[ \Sigma(\pi,\hat f)= \Sigma(\pi,f). \]

\begin{dfn}
\label{dfn-hVsheaf} {\rm For an object 
$X$ in $\sX$, let $h\cV(X)$ be the set of 
pairs $(\pi,\hat f)$ where $\pi\co E\to X$ is a
graphic submersion of fiber dimension $d+1$, 
with a $\Theta$--orientation of its vertical tangent bundle,
and $\hat f$ is a smooth section of $p^2_{\pi}\co J^2_{\pi}(E,\RR)\to E$, 
subject to the following conditions:
\begin{description} 
\item[{\bf (i)}] $(\pi,f)\co E\to X\times\RR$ is proper.
\item[{\bf (ii)}] $\hat f$ is fiberwise nonsingular.
\end{description}
}
\end{dfn}

\begin{dfn} 
\label{dfn-hWsheaf} {\rm For $X$ in $\sX$ let $h\cW(X)$ be the set of
pairs $(\pi,\hat f)$, as in definition~\ref{dfn-hVsheaf}, which 
satisfy condition (i), but where condition (ii) is 
replaced by the weaker condition
\begin{description}
\item[{\bf (iia)}] $\hat f$ is fiberwise Morse.
\end{description}
}
\end{dfn}

\begin{dfn} 
\label{dfn-hWlocsheaf} 
{\rm For $X$ in $\sX$ let $h\cW_{\loc}(X)$ be the set of
pairs $(\pi,\hat f)$, as in definition~\ref{dfn-hVsheaf}, 
but with conditions (i) and (ii)   
replaced by the weaker conditions
\begin{description}
\item[{\bf (ia)}] the map $\Sigma(\pi,\hat f)\to X\times\RR\,;\,
z\mapsto(\pi(z),f(z))$ is proper, 
\item[{\bf (iia)}] $\hat f$ is fiberwise Morse.
\end{description}
}
\end{dfn}

\medskip 
The six sheaves which we have so far defined, together with 
the obvious inclusion and jet prolongation maps, constitute 
a commutative square 

\begin{equation*}
\label{eqn-sheafstrategicdiagram} 
\xymatrix@M=7pt@W=7pt{
\cV \ar[r] \ar[d]_{j^2_{\pi}} &\cW \ar[r] \ar[d]_{j^2_{\pi}}
& \cW_{\loc} \ar[d]_{j^2_{\pi}} \\
h\cV \ar[r] & h\cW \ar[r] & h\cW_{\loc}.
}  \ttag
\end{equation*}

\subsection{Concordance theory of sheaves}
\label{subsec-homotopyofsheaves}

Let $\cF$ be a sheaf on $\sX$ and let $X$ be an object of $\sX$. 
In~\ref{dfn-concordantsections}, we defined the concordance relation 
on $\cF(X)$ and introduced the quotient set $\cF[X]$. 
It is necessary to have a relative version of $\cF[X]$. Suppose that 
$A\subset X$ is a closed subset, where $X$ is in $\sX$. Let 
$s\in \colim_U\cF(U)$ where $U$ ranges over the open 
neighborhoods of $A$ in $\sX$. Note for example that any 
$z\in \cF(\pt)$ gives rise to such an element, namely $s=\{p_U^*(z)\}$ 
where $p_U\co U\to \pt$\,. In this case we often write $z$ instead of
$s$. 

\begin{dfn} 
\label{dfn-relconcordantsections}
{\rm Let $\cF(X,A;s)\subset \cF(X)$ consist of the elements 
$t\in \cF(X)$ whose germ near $A$ is equal to $s$. Two such elements $t_0$ and
$t_1$ are concordant \emph{relative to} $A$ if they are concordant by a 
concordance whose germ near $A$ is the constant concordance from $s$
to $s$. The set of equivalence classes is denoted $\cF[X,A;s]$.} 
\end{dfn}
 
We now construct the representing space $|\cF|$ of $\cF$ and list 
its most important properties. 

\bigskip
Let $\mathbf\Delta$ be the category whose objects 
are the ordered sets $\uli n:=\{0,1,2,\dots,n\}$ for $n\ge 0$, 
with order preserving maps as morphisms. For $n\ge 0$ let 
$\Delta^n_e\subset\RR^{n+1}$ be the extended standard 
$n$-simplex, 
\[  \Delta^n_e:= \{(x_0,x_1,\dots,x_n)\in\RR^{n+1}\mid 
\Sigma x_i=1 \}. \]
An order-preserving map $\uli m\to \uli n$ induces a map of affine spaces 
$\Delta^m_e\to \Delta^n_e$. This makes $\uli n\mapsto \Delta^n_e$ into 
a covariant functor from $\Delta$ to $\sX$. 

\begin{dfn}
\label{dfn-repspaceofsheaf}
{\rm The \emph{representing space} $|\cF|$ of a sheaf $\cF$ on $\sX$ 
is the geometric realization of the simplicial set 
$\uli n\mapsto \cF(\Delta^n_e)$. 
}
\end{dfn}

\medskip
An element $z\in \cF(\pt)$ gives a point $z\in |\cF|$ and $\cF[\pt]
=\pi_0|\cF|$. In appendix~\ref{sec-spacesboxdetails} we prove that 
$|\cF|$ represents the contravariant functor $X\mapsto \cF[X]$. 
Indeed we prove the following slightly more general

\begin{prp}
\label{prp-whatitrepresents}
 For $X$ in $\sX$, let $A\subset X$ be a closed subset and $z\in
\cF(\pt)$. There is a natural bijection $\vartheta$ from the 
set of homotopy classes of maps $(X,A)\to (\,|\cF|,z)$ to 
the set $\cF[X,A;z]$. 
\end{prp}

Taking $X=S^n$ and $A$ equal to the base point, we see that 
the homotopy group $\pi_n(|\cF|,z)$ is identified with 
the set of concordance classes $\cF[S^n,\pt;z]$. We introduce 
the notation 
\[ \pi_n(\cF,z) := \cF[S^n,\pt;z] \,. \] 

\medskip
A map $v\co \cE\to \cF$ of sheaves induces a map 
$|v|\co |\cE|\to |\cF|$ of representing spaces. We call $v$ a \emph{weak 
equivalence} if $|v|$ is a homotopy equivalence. 

\begin{prp} 
\label{prp-relsurjectivity}
Let $v\co\cE\to \cF$ be a map of sheaves on $\sX$. 
Suppose that $v$ induces a surjective map 
\[ \cE[X,A;s] \lra \cF[X,A;v(s)] \]
for every $X$ in $\sX$ with a closed subset $A\subset X$ and any 
germ $s\in \colim_U\cE(U)$, where $U$ ranges over the neighborhoods 
of $A$ in $X$. Then $v$ is a weak equivalence. 
\end{prp}

\proof The hypothesis implies easily that the induced map 
$\pi_0\cE\to \pi_0\cF$ is onto and that, for any choice of base
point $z\in \cE(\pt)$, the map of concordance sets $\pi_n(\cE,z)\to
\pi_n(\cF,v(z))$ induced by $v$ is bijective. Indeed, to see that 
$v$ induces a surjection $\pi_n(\cE,z)\to
\pi_n(\cF,v(z))$, simply take $(X,A,s)=(S^n,\pt,z)$. To see that an 
element $[t]$ in the kernel of this surjection is zero, take $X=\RR^{n+1}$, 
$A=\{z\in \RR^{n+1}\mid \|z\|\ge 1\}$ and $s=p^*t$ where 
$p\co \RR^{n+1}\smin\{0\}\to S^n$ is the radial projection. The hypothesis that
$[t]$ is in the kernel amounts to a null--concordance for $v(t)$ 
which can be reformulated as an element of $\cF[X,A;v(s)]$. 
Our assumption on $v$ gives us a lift of that element to $\cE[X,A;s]$ 
which in turn can be interpreted as a null--concordance of $t$. \qed

\medskip 
Applying the representing space construction to the sheaves 
displayed in diagram~(\ref{eqn-sheafstrategicdiagram}),
we get the commutative diagram~(\ref{eqn-strategicdiagram}) 
from the introduction.

\subsection{Some useful concordances}
\label{subsec-newmodels}
\medskip

\begin{lem} 
\label{lem-basicshrink} {\rm (Shrinking lemma.)} 
Let $(\pi,f)$ be an element of 
$\cV(X)$, $\cW(X)$ or $\cW_{\loc}(X)$, with $\pi\co E\to X$
and $f\co E\to \RR$. Let $e\co X\times\RR\to \RR$ be a smooth 
map such that, for any $x\in X$, the map 
$e_x\co\RR \to\RR$ defined by  
$t\mapsto e(x,t)$ 
is an orientation preserving embedding. Let 
$E^{(1)}=\{z\in E\mid f(z)\in e_{\pi(z)}(\RR)\}$. Let $\pi^{(1)} =
\pi|E^{(1)}$
and 
\[
\begin{array}{ccc}  
f^{(1)}(z) & = & {e_{\pi(z)}}^{-1}f(z)
\end{array}
\]
for $z\in E^{(1)}$. Then $(\pi,f)$ is concordant to $(\pi^{(1)},f^{(1)})$. 
\end{lem} 

\proof Choose an $\ep>0$ and a smooth family of smooth 
embeddings $u_{(x,t)}\co \RR\to \RR$, 
where $t\in\RR$ and $x\in X$, such that $u_{(x,t)}=\id$ whenever 
$t<\ep$ and $u_{(x,1)}=e_x$ whenever $t>1-\ep$. Let 
\[
E^{(\RR)} = \left\{(z,t)\in E\times\RR \,\big|\, f(z)\in 
u_{(\pi(z),t)}(\RR)\right\}. 
\]
Then  $(z,t)\mapsto (\pi(z),t)$ defines a smooth submersion 
$\pi^{(\RR)}$ from $E^{(\RR)}$ to $X\times\RR$, and 
\[
\begin{array}{ccc}
z & \mapsto & {u_{(\pi(z),t)}}^{-1}f(z) 
\end{array}
\]
defines a smooth map $f^{(\RR)}\co E^{(\RR)}\to \RR$. Now
$(\pi^{(\RR)},f^{(\RR)})$ is a concordance  
from $(\pi,f)$ to $(\pi^{(1)},f^{(1)})$, modulo some simple 
re--labelling of the elements of $E{(\RR)}$ to ensure that 
$\pi^{(\RR)}$ is graphic. (As it stands, $E$ is a subset of 
$Z\times X$, compare~\ref{dfn-graphicmaps}, 
and $E^{(\RR)}$ is a subset of $(Z\times X)\times\RR$. But we want 
$E^{(\RR)}$ to be a subset of $Z\times(X\times\RR)$; hence the 
need for relabelling.) \qed

\medskip
Lemma~\ref{lem-basicshrink} has an obvious analogue for the sheaves 
$h\cV$, $h\cW$ and $h\cW_{\loc}$, which we do not state explicitly. 

\begin{lem} 
\label{lem-shrinking2} 
Every class in $\cW[X]$ or $h\cW[X]$ has a representative
$(\pi,f)$, resp.  $(\pi,\hat f)$, in which $f\co E\to\RR$ is a bundle 
projection, so that 
\[ 
E\cong f^{-1}(0)\times\RR\,.
\]
\end{lem}

\proof We concentrate on the first case, starting with an arbitrary  
$(\pi,f)\in \cW[X]$. We do not assume that $f\co E\to \RR$ is a 
bundle projection to begin with. However, by Sard's theorem we can 
find a regular value $c\in \RR$ for $f$. The singularity set of 
$f$ (not to be confused 
with the fiberwise singularity set of $f$) is closed in $E$. 
Therefore its image under the proper map $(\pi,f)\co E\to X\times\RR$
is closed. (Proper maps between locally compact spaces are closed maps). 
The complement of that image is an open
neighborhood $U$ of $X\times\{c\}$ 
in $X\times \RR$ containing no critical points of $f$.
It follows easily that there exists $e\co X\times\RR\to \RR$ 
as in lemma~\ref{lem-basicshrink}, with $e(x,0)=c$ for all $x$ 
and $(x,e(x,t))\in U$ for all $x\in X$ and $t\in\RR$. Apply 
lemma~\ref{lem-basicshrink} with this choice of $e$. 
In the resulting $(\pi^{(1)},f^{(1)})\in \cW(X)$, the map 
$f^{(1)}\co E^{(1)}\to \RR$ is nonsingular and proper, hence a 
bundle projection. (It is not claimed that $f^{(1)}$ is 
fiberwise nonsingular.) \qed 

\medskip
We now introduce two sheaves $\cW^0$ and $h\cW^0$ on $\sX$. They are 
weakly equivalent to $\cW$ and $h\cW$, respectively, but 
better adapted to Vassiliev's integrability theorem, 
as we will explain in section~\ref{sec-Vassiliev}. 

\begin{dfn}
\label{dfn-W0sheaf}
{\rm For $X$ in $\sX$ let $\cW^0(X)$ be the set of 
all pairs $(\pi,f)$ as in definition~\ref{dfn-Wsheaf}, 
replacing however condition (iia) there by the weaker 
\begin{description}
\item[{\bf (iib)}] $f$ is fiberwise Morse in some neighborhood 
of $f^{-1}(0)$. 
\end{description}
}
\end{dfn}

\begin{dfn}
\label{dfn-hW0sheaf}
{\rm For $X$ in $\sX$ let $h\cW^0(X)$ be the set of 
all pairs $(\pi,\hat f)$ as in definition~\ref{dfn-hWsheaf}, 
replacing however condition (iia) by the weaker 
\begin{description}
\item[{\bf (iib)}] $\hat f$ is fiberwise Morse in some neighborhood 
of $f^{-1}(0)$. 
\end{description}
}
\end{dfn} 

\medskip
From the definition, there are inclusions $\cW\to \cW^0$
and $h\cW\to h\cW^0$. There is also a jet prolongation map 
$\cW^0\to h\cW^0$ which we may regard as an inclusion, the 
inclusion of the subsheaf of integrable elements.

\begin{lem}
\label{lem-shrinking1}
The inclusions $\cW\to \cW^0$ and $h\cW\to h\cW^0$ 
are weak equivalences. 
\end{lem}

\proof We will concentrate on the first of the two inclusions, 
$\cW\to \cW^0$. Fix $(\pi,f)$ in $\cW^0(X)$, with $\pi\co E\to X$ 
and $f\co E\to\RR$. We will subject $(\pi,f)$ to a concordance 
ending in $\cW(X)$. Choose an open neighborhood $U$ of 
$f^{-1}(0)$ in $E$ such that, for each $x\in X$, the critical points 
of $f_x =f|E_x$ on $E_x\cap U$ are all nondegenerate. Since $E\smin U$ 
is closed in $E$ and the map $(\pi,f)\co E\to X\times\RR$ is proper, 
the image of $E\smin U$ under that map is a closed subset of
$X\times\RR$ which has empty intersection with $X\times 0$. 
Again it follows that a map $e\co X\times\RR\to \RR$ 
as in~\ref{lem-basicshrink} can be constructed 
such that $e(x,0)=0$ for all $x$ 
and $(x,e(x,t))\in U$ for all $(x,t)\in X\times\RR$. 
As in the proof of lemma~\ref{lem-basicshrink}, use $e$
to construct a concordance from $(\pi,f)$ to some element 
$(\pi^{(1)},f^{(1)})$ which, by inspection, belongs to 
$\cW(X)$. If the restriction of
$(\pi,f)$ to an open neighborhood $Y_1$ of a closed $A\subset X$ 
belongs to $\cW(Y_1)$, then the concordance can be made relative to
$Y_0$, where $Y_0$ is a smaller open neighborhood of $A$ in $X$. \qed

\section{The lower row of diagram~(\ref{eqn-strategicdiagram})}
\label{sec-bottomrow} 

This section describes the homotopy types of the spaces 
in the lower row of~(\ref{eqn-strategicdiagram}) in 
bordism--theoretic terms. One of the conclusions is that the 
lower row is a homotopy fiber sequence, proving
theorem~\ref{thm-bottomrow}. 
We also show that 
the jet prolongation map $|\cW_{\loc}|\to |h\cW_{\loc}|$ is a 
homotopy equivalence (the fact as such does not belong in this
section, but its proof does). In the standard 
case where $d=2$ and $\Theta=\pi_0(\GL)$, the space $|h\cV|$ 
will be identified with $\Omega^{\infty}\CC P^{\infty}_{-1}$.

\subsection{A cofiber sequence of Thom spectra} 
\label{subsec-thomcofiber} 
Let $\Gr\cW(d+1,n)$ be the space of triples 
$(V,\ell,q)$ consisting of 
a $\Theta$--oriented $(d+1)$-dimensional linear 
subspace $V\subset\RR^{d+1+n}$, 
a linear map $\ell\co V\to \RR$ and a quadratic form 
$q\co V\to \RR$, subject to the condition that if $\ell=0$, 
then $q$ is nondegenerate. $\Gr\cW(d+1,n)$ classifies 
$(d+1)$-dimensional $\Theta$--oriented 
vector bundles whose fibers have the above extra structure, 
i.e., each fiber $V$ comes equipped with a Morse type map 
$\ell+q\co V\to \RR$ and with a linear embedding into $\RR^{d+1+n}$. 

\medskip
The tautological $(d+1)$-dimensional vector bundle $\U_{n}$ on 
$\Gr\cW({d+1},n)$ 
is canonically embedded in a trivial bundle 
$\Gr\cW({d+1},n)\times\RR^{{d+1}+n}$. 
Let 
\[ 
\begin{array}{ccc}
\U_{n}^{\bot}& \subset & \Gr\cW({d+1},n)\times\RR^{{d+1}+n} 
\end{array}
\] 
be the orthogonal complement, an $n$-dimensional 
vector bundle on $\Gr\cW({d+1},n)$. 
The tautological bundle $\U_{n}$ comes equipped with the extra
structure consisting of a map from (the total space of) $\U_{n}$ 
to $\RR$ which, on each fiber of $\U_{n}$, is a Morse type map.
(The fiber of $\U_{n}$ over a point $(V,q,\ell)\in \Gr\cW({d+1},n)$ 
is identified with the $(d+1)$-dimensional vector space $V$ and the 
map can then be described as $\ell+q$.) 

\medskip
Let $S(\RR^{d+1})$ be the vector space of quadratic forms on $\RR^{d+1}$ 
(or equivalently, symmetric $(d+1)\times (d+1)$ matrices) 
and $\Delta\subset S(\RR^{d+1})$ the subspace of the degenerate forms 
(not a linear subspace). The complement 
$Q(\RR^{d+1})=S(\RR^{d+1})\smin\Delta$ is the space of non-degenerate
quadratic forms on $\RR^{d+1}$. Since quadratic forms can be diagonalized, 
\[ Q(\RR^{d+1}) = \coprod_{i=0}^{d+1} Q(i,d+1-i) \]
where $Q(i,d+1-i)$ is the connected component containing 
the form $q_i$ given by 
\[  q_i(x_1,x_2,\dots,x_{d+1})= -(x_1^2+\cdots+x_i^2)+(x_{i+1}^2
+\cdots+x_{d+1}^2). \]
The stabilizer $\Or(i,d+1-i)$ of $q_i$ for the (transitive)
action of $\GL(d+1)$ on $Q(i,d+1-i)$ 
has $\Or(i)\times\Or(d+1-i)$ as a maximal compact subgroup and 
$\GL(d+1)$ has $\Or(d+1)$ as a maximal compact subgroup. Hence the 
inclusion 
\[ 
\begin{array}{cccccc}
(\Or(i)\times\Or(d+1-i))\big\backslash \Or(d+1) & \lra & Q(i,d+1-i)\,; & \quad 
\textup{coset of }g & \mapsto &  q_ig 
\end{array}
\]
is a homotopy equivalence, and therefore the subspace 
\begin{equation*}
\label{eqn-Qtopieces}
\begin{array}{rcl} 
Q^0(\RR^{d+1}) & = & \{q_0,q_1,\dots,q_{d+1}\}\cdot\Or(d+1) \\
           & \cong & \rule{0mm}{5mm}\coprod_{i=0}^{d+1}\,
(\Or(i)\times\Or(d+1-i))\big\backslash\Or(d+1)
\end{array} \ttag
\end{equation*}
of $Q(\RR^{d+1})$ is a deformation retract, $Q(\RR^{d+1})\simeq Q^0(\RR^{d})$.

\medskip
For the submanifold $\Sigma(d+1,n)\subset \Gr\cW(d+1,n)$ 
consisting of the triples 
$(V,\ell,q)$ with $\ell=0$ we have  
\begin{equation*}
\label{eqn-UfromQ}
\begin{array}{ccc}
\Sigma(d+1,n)  & \cong & \big(\Or(d+1+n)/\Or(n) \times
Q(\RR^{d+1})\times\Theta\big)\big/\Or(d+1)\,. 
\end{array} \ttag 
\end{equation*}
The restriction of $\U_{n}$ to $\Sigma(d+1,n)$ comes equipped with the 
extra structure of 
a fiberwise nondegenerate quadratic form. There is a canonical 
normal bundle for $\Sigma(d+1,n)$ in $\Gr\cW(d+1,n)$ which
is easily identified with the dual bundle $\U_{n}^*|\Sigma(d+1,n)$. 
Hence there is a homotopy cofiber sequence
\[ 
\xymatrix@M+2pt{
\Gr\cV(d+1,n)\incto[r] & \Gr\cW(d+1,n) \ar[r] & 
\Th(\U_{n}^*|\Sigma(d+1,n)) 
}
\]
where $\Gr\cV(d+1,n)=\Gr\cW(d+1,n)\smin\Sigma(d+1,n)$ 
and $\Th(\dots)$ denotes the Thom space. This leads to a homotopy 
cofiber sequence of Thom spaces 
\[ \Th(\U_{n}^{\bot}|\Gr\cV(d+1,n)) \lra \Th(\U_{n}^{\bot})
\lra \Th(\U_{n}^{\bot}\oplus \U_{n}^*|\Sigma(d+1,n)). \]
(A \emph{homotopy cofiber sequence} is a diagram $A\to B\to C$ of
spaces, where $C$ is pointed, together with a nullhomotopy of the 
composite map $A\to B$ such that 
the resulting map from $\cone(A\to B)$ to $C$ is a weak homotopy 
equivalence.) \newline
We view the space $\Th(U_{n}^{\bot})$ as the $(n+d)$--th space
in a spectrum $\bW$, and similarly for the other two Thom spaces. 
Then as $n$ varies the sequence above becomes a homotopy cofiber sequence 
of spectra
\[ \bV\lra \bW\lra \bW_{\loc}. \]
We then have the corresponding infinite loop spaces 
\[ 
\begin{array}{rcl} 
\Omega^{\infty}\bV & = 
& \colim_n\, \Omega^{d+n}\Th(\U_{n}^{\bot}|\Gr\cV(d+1,n))\,, \\
\Omega^{\infty}\bW & = & \colim_n\, \Omega^{d+n}\Th(\U_{n}^{\bot})\,, \\
\Omega^{\infty}\bW_{\loc} & = &  
\colim_n\, \Omega^{d+n}\Th(\U_{n}^{\bot}\oplus \U_{n}^*|\Sigma(d+1,n)).
\end{array}
\]
(We use CW--models for the spaces involved. For example,
$\Omega^{d+n}\Th(\U_{n}^{\bot})$ can be considered as 
the representing space of the sheaf on $\sX$ 
which to a smooth $X$ associates the set of pointed maps from 
$X_+\wedge S^{d+n}$ to $\Th(\U_{n}^{\bot})$. The representing space 
is a CW--space.) 

The homotopy cofiber sequence of spectra above yields a homotopy 
fiber sequence of infinite loop spaces 
\begin{equation*}
\label{eqn-hofiseinlosp}
\Omega^{\infty}\bV \lra \Omega^{\infty}\bW \lra 
\Omega^{\infty}\bW_{\loc}\,\,, \ttag
\end{equation*}
that is, $\Omega^{\infty}\bV$ is homotopy equivalent to the homotopy 
fiber of the right-hand map. (A \emph{homotopy fiber sequence} 
is a diagram of spaces $A\to B\to C$, where $C$ is pointed, together
with a nullhomotopy of the composite map $A\to C$ such that the 
resulting map from $A$ to $\hofiber(B\to C)$ is a weak homotopy 
equivalence.) 
In particular there 
is a long exact sequence of homotopy groups associated with
diagram~(\ref{eqn-hofiseinlosp}) and a Leray-Serre spectral sequence 
of homology groups.

Suppose that a topological group $G$ acts on a space $Q$ from 
the right. We use the notation $Q_{hG}$ for the 
``Borel construction'' or homotopy orbit space $Q\times_G EG$.
 
\begin{lem}
\label{lem-bWloctype}
There is a homotopy equivalence of infinite loop spaces 
\[ \Omega^{\infty}\bW_{\loc} \simeq  \Omega^{\infty}
S^{1+\infty}(\Sigma(d+1,\infty)_+) \]
where $(\Sigma(d+1,\infty)$ is a disjoint union of homotopy 
orbit spaces,
\[ \Sigma(d+1,\infty) \simeq \coprod_{i=0}^{d+1} 
\Theta_{h\Or(i,d+1-i)}.
\] 
\end{lem}

\proof Since $\U_{n}|\Sigma(d+1,n)$ comes equipped with a fiberwise 
nondegenerate quadratic form, $\U_{n}^*|\Sigma(d+1,n)$ is canonically  
identified with $\U_{n}|\Sigma(d+1,n)$. 
Consequently the restriction 
\[ 
\U_{n}^{\bot}\oplus
\U_{n}^*\,\big|\,\Sigma(d+1,n)
\]
is trivialized, so that
$\Th(\U_{n}^{\bot}\oplus \U_{n}^*\,\big|\,\Sigma(d+1,n)) 
\, \simeq\,  S^{d+1+n}(\Sigma(d+1,n)_+)\,.$
Hence 
\[ \Omega^{\infty}\bW_{\loc}\,\, \simeq  \,\,  
\Omega^{\infty}S^{1+\infty}(\Sigma(d+1,\infty)_+)  \] 
where $\Sigma(d+1,\infty)=\bigcup\Sigma(d+1,n)$. Using the
description~(\ref{eqn-UfromQ}) of $\Sigma(d+1,n)$ 
and the equivariant homotopy equivalence 
$Q(\RR^{d+1})\simeq Q^0(\RR^{d+1})$, see~(\ref{eqn-Qtopieces}), we get 
\[ \Sigma(d+1,n) \,\, \simeq \,\,
\left(\Or(d+1+n)/\Or(n))\times Q^0(\RR^{d+1})\times\Theta\right)
\big/\Or(d+1). \]
The union $\bigcup_n\, \Or(d+1+n)/\Or(n)$ is a contractible 
free $\Or(d+1)$-space, so that $\Sigma(d+1,\infty)$ is homotopy 
equivalent to the homotopy orbit space of the canonical right action of 
$\Or(d+1)$ on the space  
\[  
Q^0(\RR^{d+1})\times \Theta \, \cong\, \left(\coprod_{i=0}^{d+1}\,
(\Or(i)\times\Or(d+1-i))\big\backslash\Or(d+1)\right)\times\Theta\,. 
\]
That in turn is homotopy equivalent to the disjoint 
union over $i$ of the homotopy orbit spaces of
$\Or(i)\times\Or(d+1-i)\simeq \Or(i,d+1-i)$ acting on the left of 
$(O(d+1)\times\Theta)\big/O(d+1)\cong\Theta$. \qed

\medskip 
Let $\Gr(d,n;\Theta)$ be the space of $d$--dimensional 
$\Theta$--oriented linear subspaces in $\RR^{d+n}$.
It can be identified with a subspace of 
$\Gr\cV(d+1,n)=\Gr\cW(d+1,n)\smin\Sigma(d+1,n)$,
consisting of the $(V,\ell+q)$ where $V$ contains the subspace
$\RR\times0\times0$ of $\RR\times \RR^d\times\RR^n$, 
and $\ell+q$ is the linear projection 
to that subspace (so that $q=0$). The injection 
is covered by a fiberwise isomorphism of vector bundles 
\[ T_n^{\bot} \lra \U_{n}^{\bot}\,\big|\,\Gr\cV(d+1,n) \] 
where $T_n^{\bot}$ 
is the standard n-plane bundle on $\Gr(d,n;\Theta)$. 

\begin{lem} 
\label{lem-bVtype} The induced map of Thom spaces 
$\Th(T_n^{\bot}) \lra \Th(\U_{n}^{\bot}\,|\,\Gr\cV(d+1,n))$
is $(d+2n-1)$--connected. Hence 
\[ 
\begin{array}{ccc}
\Omega^{\infty}\bV & \simeq
& \colim_n\, \Omega^{d+n}\Th(T_n^{\bot})\,.
\end{array}
\]
\end{lem} 

\proof It is enough to show that the inclusion of $\Gr(d,n;\Theta)$
in $\Gr\cV(d+1,n)$ is $(d+n-1)$--connected. Viewing both of these 
spaces as total spaces of certain bundles with fiber $\Theta$ reduces 
the claim to the case where $\Theta$ is a single point. Note also that 
$\Gr\cV(d+1,n)$ has a deformation retract consisting of the pairs 
$(V,\ell+q)$ with $q=0$ and $\|\ell\|=1$. 
This deformation retract is homeomorphic 
to the coset space
\[\Or(d)\times\Or(n)\big\backslash\Or(1+d+n),\]
assuming $\Theta=\pt$. We are therefore looking at the 
inclusion of $\left(\Or(d)\times\Or(n)\right)\big\backslash\Or(d+n)$
in $\left(\Or(d)\times\Or(n)\right)\big\backslash\Or(1+d+n)$, which is
indeed $(d+n-1)$--connected. \qed

\medskip
In the standard case where $d=2$ and $\Theta=\pi_0\GL$, we may 
compare the Grassmannian of oriented planes $\Gr(2,2n;\Theta)$
with the complex projective $n$--space. The map 
\[ \CC P^n\lra G(2,2n;\Theta) \]
that forgets the complex structure is $(2n-1)$--connected. The pullback 
of $T_{2n}^{\bot}$ under this map is the realification of the 
tautological complex $n$--plane bundle $L_n^{\bot}$ and the 
associated map of Thom spaces is $(4n-1)$--connected. The spectrum 
$\CC P^{\infty}_{-1}$ with $(2n+2)$--nd space $\Th(L_n^{\perp})$
is therefore weakly equivalent to the Thom spectrum $\bV$. 
We can now collect the main conclusions of this section, 
~\ref{subsec-thomcofiber}, in  

\begin{prp} 
\label{prp-VWexplained} For $d=2$ and $\Theta=\pi_0\GL$, 
the homotopy fiber sequence~{\rm (\ref{eqn-hofiseinlosp})} is 
homotopy equivalent to
\[ \Omega^{\infty}\CC\bP^{\infty}_{-1}
\lra \Omega^{\infty}\bW \lra 
\Omega^{\infty} S^{1+\infty}
\big(\big(\coprod_{i=0}^{3} B\SOr(i,3-i)\big)_+\big). \qed \]
\end{prp}

\subsection{The spaces $|h\cW|$ and $|h\cV|$} 

In section~\ref{subsec-mocksheaves} we described 
the jet bundle $J^2(E,\RR)$ and its fiberwise version as 
certain spaces of smooth map germs $(E,z)\to \RR$, modulo
equivalence. For our use in this section and the next it is 
better to view it as a construction on the tangent bundle. 
For a vector space $V$, let $J^2(V)$ denote the vector space 
of maps 
\[ \begin{array}{cccccc}
 \hat f\co & V\to  & \RR  \,, &\quad 
\hat f(v) & = & c+\ell(v)+q(v) 
\end{array}
\]
where $c\in\RR$ is a constant, $\ell\in V^*$ and $q\co V\to \RR$ 
is a quadratic map. This is a contra\-variant continuous 
functor on vector spaces, so extends to a functor on vector 
bundles with $J^2(F)_z=J^2(F_z)$ when $F$ is a vector bundle over $E$. 

When $F=TE$ is the tangent bundle of a manifold $E$, then there is 
an isomorphism of vector bundles 
\[  J^2(E,\RR) \cong J^2(TE). \]
Indeed after a choice of a connection on $TE$, 
the associated exponential map induces a diffeomorphism germ 
$\exp_z\co (TE_z,0)\to (E,z)$.
Composition with $\exp_z$ is an isomorphism from $J^2(E,\RR)_z$ 
to $J^2(TE_z)$. 

\begin{lem} 
\label{lem-fiberjet} Let $\pi\co E\to X$ be a smooth submersion. 
Any choice of connection on the vertical tangent bundle $T^{\pi}E$ 
induces an isomorphism 
\[ J^2_{\pi}(E,\RR) \lra J^2(T^{\pi}E). \]
This is natural under pullbacks of submersions. 
\end{lem}

\proof In addition to choosing a connection on $T^{\pi}E$, we may choose 
a smooth linear section of the vector bundle surjection 
$d\pi:TE\to \pi^*TX$ and a connection on $TX$. This leads to a
splitting 
\[ TE \cong T^{\pi}E\oplus \pi^*TX \]
and determines a direct sum connection on $TE$. The associated
exponential diffeomorphism germ $\exp\co (TE_z,0) \lra (E,z)$
is fiberwise, i.e., it restricts to a diffeomorphism 
germ 
\begin{equation*}
\label{eqn-diffeogerm}
((T^{\pi}E)_z,0)\to (E_{\pi(z)},z) \ttag
\end{equation*}
for each $z\in E$. 
Indeed, the chosen connection on $T^{\pi}E$ restricts to a connection 
on the tangent bundle of $E_{\pi(z)}$, and any geodesic in
$E_{\pi(z)}$ for that connection is clearly a geodesic in $E$ as well. 
The argument also shows that the diffeomorphism
germ~(\ref{eqn-diffeogerm}), and the isomorphism 
$J^2_{\pi}(E,\RR)_z \lra J^2(T^{\pi}E)_z$ which it induces,
depend only on the choice of a 
connection on $T^{\pi}E$, but not on the choice of a splitting of 
$d\pi:TE\to \pi^*TX$ and a connection on $TX$. (However, 
making use of all the choices, we arrive at a commutative 
diagram of vector bundles 
\medskip
\begin{equation*} 
\xymatrix@M=8pt@W=8pt{
J^2(TE) \onto[r]^{i^*} \ar[d]^{\cong} & J^2(T^{\pi}E) \ar[d]^{\cong}  \\
J^2(E,\RR) \onto[r]^{i^*} & J_{\pi}^2(E,\RR)
} 
\end{equation*}
where the horizontal epimorphisms are induced by inclusions.) 
Finally, if 
\medskip
\begin{equation*} 
\xymatrix@M=8pt@W=8pt{
\varphi^*E \ar[r]^{\bar\varphi} \ar[d]^{\varphi^*\pi} & E \ar[d]^{\pi}  \\
Y  \ar[r]^{\varphi} & X
}
\end{equation*}
is a pullback diagram of submersions, then a choice of connection 
on $T^{\pi}E$ determines a connection on $\bar\varphi^*T^{\pi}E\cong
T^{\varphi^*\pi}\varphi^*E$. The resulting exponential diffeomorphism germs 
are related by a commutative diagram 
\medskip
\begin{equation*} 
\xymatrix@M=8pt@W=8pt{
((T^{\varphi^*\pi}\varphi^*E)_z,0) \ar[r] \ar[d]^{d^{\pi}\bar\varphi} & 
(\varphi^*E_{\varphi^*\pi(z)},z) \ar[d]^{\bar\varphi} \\
((T^{\pi}E)_{\bar\varphi(z)},0) \ar[r] & 
(E_{\pi\bar\varphi(z)},\bar\varphi(z))\,.   
}  
\end{equation*}
This proves the naturality claim. \qed

\medskip
We can re--define $h\cW(X)$ in definition~\ref{dfn-hWsheaf} 
as the set of certain pairs $(\pi,\hat f)$ much as before,
with $\pi\co E\to X$, where $\hat f$ is now a section of
$J^2(T^{\pi}E)$. 
The above lemma tells us that the new definition of $h\cW$ is related to 
the old one by a chain of two weak equivalences. (In the 
middle of that chain is yet another variant of $h\cW(X)$, 
namely the set of triples $(\pi,\hat f,\nabla)$ where $\pi$ and $\hat f$ 
are as in definition~\ref{dfn-hWsheaf}, while $\nabla$ is a connection 
on $T^{\pi}E$.)

\medskip 
Our object now is to construct a natural map 
\begin{equation*}
\label{eqn-Wtobordism}
\tau\co h\cW[X] \lra [X,\Omega^{\infty}\bW]. \ttag 
\end{equation*} 
Here $[\,\,,\,\,]$ in the right-hand side denotes a set 
of homotopy classes of maps. 

We assume familiarity with the 
Pontryagin-Thom relationship between Thom spectra and their 
infinite loop spaces on the one hand, and bordism theory 
on the other. One direction of this relies on transversality theorems, 
the other uses collapse maps to normal bundles of submanifolds
in euclidean spaces.
See \cite{Stong68} and especially \cite{Quillen71}. 
Applied to our situation this identifies $[X,\Omega^{\infty}\bW]$ 
with a group of bordism classes of certain triples $(M,g,\hat g)$. Here 
$M$ is smooth without boundary, $\dim(M)=\dim(X)+d$, and $g, \hat g$
together constitute a vector bundle pullback square 
\begin{equation*}
\label{eqn-cobordismW}
\xymatrix{
TM\times\RR\times\RR^j \ar[r]^-{\hat g} \ar[d] &  
TX\times \U_{\infty}\times\RR^j \ar[d] \\
M \ar[r]^-g  &  X\times \Gr\cW(d+1,\infty)
} \ttag
\end{equation*}
such that the $X$--coordinate of $g$ is 
a proper map $M\to X$. The $\RR^j$ factor in the top row, 
with unspecified $j$, is there 
for stabilization purposes. The map $\hat g$ should be thought 
of as a \emph{stable} vector bundle map 
from $TM\times\RR$ to $TX\times\U_{\infty}$, covering $g$, where 
$U_{\infty}$ is the tautological vector bundle 
of fiber dimension $d+1$ 
on $\Gr\cW(d+1,\infty)$.

Let now $(\pi,\hat f)\in h\cW(X)$, where 
$\hat f$ is a section of $J^2(T^{\pi}E)\to E$ with 
underlying map $f\co E\to \RR$. See definition~\ref{dfn-hWsheaf}. 
After a small deformation which does not affect the concordance 
class of $(\pi,\hat f)$, we may assume that $f$ is transverse 
to $0\in\RR$ (not necessarily fiberwise) and get a manifold 
$M=f^{-1}(0)$ with $\dim(M)=\dim(X)+d$. The 
restriction of $\pi$ to $M$ is a proper map 
$M\to X$, by the definition of $h\cW(X)$. The section $\hat f$ yields
for each $z\in E$ a map 
\[ 
\begin{array}{cccc}
\hat f(z)& = & f(z)+\ell_z+q_z\co & (T^{\pi}E)_z\to \RR 
\end{array}
\]
with the property that the quadratic term $q_z$ is nondegenerate 
when the linear term $\ell_z$ is zero. For $z\in M$ the constant 
$f(z)$ is zero, so the restriction $T^{\pi}E|M$ is a $(d+1)$--dimensional 
vector bundle on $M$ with the extra structure considered 
in section~\ref{subsec-thomcofiber}. Thus $T^{\pi}E|M$ is classified by 
a map from $M$ to the space $\Gr\cW(d+1,\infty)$: 
there is a bundle diagram 
\[ 
\xymatrix{
T^{\pi}E|M \ar[r] \ar[d]  & \U_{\infty} \ar[d] \\
M \ar[r]^-{\kappa} &  \Gr\cW(d+1,\infty).  
}
\]  
Let $g\co M\lra X\times \Gr\cW(d+1,\infty)$ be the map
$z\mapsto(\pi(z),\kappa(z))$. We now have 
a canonical vector bundle map
\[ 
\begin{array}{ccccccc}
\hat g\co TM\times \RR & \cong & TE|M & \cong & \pi^*TX|M\oplus T^{\pi}E|M 
& \lra & TX\times \U_{\infty} 
\end{array}
\]
and we get a triple 
$(M,g,\hat g)$ which represents an element of
$[X,\Omega^{\infty}\bW]$ in the 
bordism-theoretic description. It is easily verified that 
the bordism class of $(M,g,\hat g)$ depends only on the concordance 
class of the pair 
$(\pi,\hat f)$. Thus we have defined the map $\tau$ 
of~(\ref{eqn-Wtobordism}). 

\begin{thm}
\label{thm-Wtobordism} The natural map 
$\tau\co h\cW[X]\to [X,\Omega^{\infty}\bW]$ is a bijection 
when $X$ is a closed manifold.
\end{thm}

\proof We define a map $\sigma$ in the other direction by 
running the construction $\tau$ backwards. We use the 
bordism group description~(\ref{eqn-cobordismW}) of 
$[X,\Omega^{\infty}\bW]$. Let 
$(M,g,\hat g)$ be a representative, with 
$g\co M\to X\times\Gr\cW(d+1,\infty)$ and 
\[ 
\begin{array}{ccc} 
\hat g\co TM\times\RR\times \RR^j & \lra & 
TX\times \U_{\infty}\times \RR^j~. 
\end{array}
\]
By obstruction theory, see lemma~\ref{lem-destabilization} below, we can 
suppose that $j=0$. We write $E=M\times\RR$ and  
$\pi_E\co E\to X$ 
for the composition of the projection $E\to M$ with the first 
component of $g$. The map $\hat g$, now with $j=0$,
has a first component $TM\times\RR\to TX$. We (pre--)compose it with 
the evident vector bundle map from $TE\cong TM\times T\RR$ 
to $TM\times \RR$ which covers the projection 
from $E\cong M\times\RR$ to $M$. The result is a map of 
vector bundles 
\[
\hat \pi_E\co TE \lra TX,
\]
covering $\pi_E$ and surjective in the fibers. Since $E$ is an open 
manifold, Phillips' submersion theorem \cite{Phillips67}, 
\cite{Gromov71}, \cite{Haefliger71} applies to show that 
$(\pi_E,\hat\pi_E)$ is homotopic through 
fiberwise surjective bundle maps to a pair $(\pi,d\pi)$ where 
$\pi\co E\to X$ is a submersion and $d\pi\co TE\to TX$ is its 
differential. \newline 
This homotopy lifts to a homotopy of 
vector bundle maps which are isomorphic on the fibers, starting with  
$\hat g\co TE\to TX\times \U_{\infty}$ and ending with a map 
$TE\to TX\times  \U_{\infty}$ which refines the differential 
$d\pi\co TE\to TX$. 
Its restriction to $T^{\pi}E\subset TE$ is a vector  
bundle map $T^{\pi}E\to \U_{\infty}$, still isomorphic on 
the fibers, which equips each 
fiber $(T^{\pi}E)_z$ of $T^{\pi}E$ with a Morse type map 
\[ \ell_z+q_z\co (T^{\pi}E)_z\to \RR . \]
Let $f\co E\to \RR$ be the 
projection onto the $\RR$ factor, and let
\[ \hat f(z)= f(z)+\ell_z+q_z \in J^2(T^{\pi}E). \] 
The map $f$ is proper, since $X$ and hence $M$
are compact. Consequently the pair $(\pi,\hat f)$ 
represents an element in $h\cW[X]$. Its 
concordance class depends only on the bordism class of 
$(M,g,\hat g)$; the verification uses a relative version 
of lemma~\ref{lem-destabilization}. This describes a map 
\[ \sigma\co [X,\Omega^{\infty}\bW] \lra  h\cW[X]. \]
It is obvious from the constructions that $\tau\circ\sigma=\id$. 
In order to evaluate the composition $\sigma\circ\tau$, it suffices by 
lemma~\ref{lem-shrinking2} to evaluate it on an element $(\pi,\hat f)$ 
where $f\co E\to\RR$ is regular, so that $E\cong M\times\RR$
with $M=f^{-1}(0)$. For $(y,r)\in M\times\RR$, the map 
\[ 
\begin{array}{ccc} 
\hat f(y,r)\co (T^{\pi}(M\times\RR))_{(y,r)} & \lra & \RR 
\end{array}
\]
is a second degree polynomial of Morse type. The homotopy 
\[ \hat f_{t}(y,r) = \hat f(y,tr)+(1-t)r\,,  \]
suitably reparametrized,
shows that $(\pi,\hat f)$ is concordant to 
$(\pi,\hat f_0)$, which represents the image of 
$(\pi,\hat f)$ under $\sigma\circ\tau$. 
Therefore $\sigma\circ\tau=\id$.  \qed

\begin{lem} 
\label{lem-destabilization} Let $T$ and $U$ be $k$-dimensional 
vector bundles over a manifold $M$. Let $\iso(T,U)\to M$ be the 
fiber bundle on $M$ whose fiber at $x\in M$ is the space of 
linear isomorphisms from $T_x$ to $\U_x$.
The stabilization map $\iso(T,U)\to \iso(T\times\RR,U\times\RR)$
induces a map of section spaces which is $(k-\dim(M)-1)$--connected.
\end{lem}

\proof We use the following general principle. Suppose that $Y\to M$
and $Y'\to M$ are fibrations and that $f\co Y\to Y'$ is a map 
over $M$. Suppose that for each $x\in M$, the 
restriction $Y_x\to Y'_x$ of $f$ to the fibers over $x$ is
$c$--connected. Then the induced map of section spaces, 
$\Gamma(Y)\to \Gamma(Y')$, is $(c-m)$--connected where
$m=\dim(M)$. \newline
The proof of this proceeds as follows: Fix $s\in \Gamma(Y')$. 
The homotopy fiber of $\Gamma(Y)\to \Gamma(Y')$ over $s$ is 
easily identified with the section space $\Gamma(Y'')$ of another 
fibration $Y''\to M$, defined by 
\[  Y''_x = \hofiber_{s(x)}\,(Y_x\to Y'_x)\,. \]
By assumption each $Y''_x$ is $(c-1)$--connected. Hence by obstruction 
theory or a simple induction over skeletons, $\Gamma(Y'')$ is 
$(c-1-m)$--connected. Since this holds for arbitrary $s$, all homotopy 
fibers of $\Gamma(Y)\to \Gamma(Y')$ are $(c-1-m)$--connected. 
Consequently $\Gamma(Y)\to \Gamma(Y')$ is $(c-m)$--connected. \newline
Now for the application: The inclusion $\GL(k)\to \GL(k+1)$ 
is $(k-1)$--connected. Hence the stabilization map 
$\iso(T,U)\to \iso(T\times\RR,U\times\RR)$
is $(k-1)$--connected on the fibers, and so induces a
$((k-1)-m)$--connected map of section spaces. \qed

\bigskip 
 
The following concept will be useful in a fiberwise 
version of the Pontryagin-Thom construction which we will 
need in a moment.

\begin{dfn}
\label{dfn-verticaltube}
 {\rm Let $p\co Y\to X$ be a smooth submersion. Let $C$ 
be a smooth submanifold of $Y$ 
and suppose that $p|C$ is still a submersion. A \emph{vertical} tubular 
neighborhood for $C$ in $Y$ consists of a smooth vector bundle 
$q\co N\to C$ with zero section $s$,
and an open embedding $e\co N\to Y$ 
such that $es=\textup{inclusion}\co C\to Y$ and 
$pe=pqe\co N\to X$. 
}
\end{dfn}

Now we give a detailed description of a map
$|h\cW|\to \Omega^{\infty}\bW$ which induces~(\ref{eqn-Wtobordism}). 
It relies entirely on the Pontryagin--Thom collapse construction. 

We begin by describing a variant $h\cW^{(r)}$ of $h\cW$, 
depending on an integer $r>0$. 
Fix $X$ in $\sX$. An element of $h\cW^{(r)}(X)$ is 
a quadruple $(\pi,\hat f,w,N)$ where $\pi\co E\to X$ and 
$\hat f$ are as in definition~\ref{dfn-hWsheaf}. The remaining 
data are a smooth embedding 
\[
\begin{array}{ccc}
 w\co E & \lra & X\times\RR\times\RR^{d+r} 
\end{array}
\] 
which covers $(\pi,f)\co E\to X\times\RR$, and a \emph{vertical} tubular 
neighborhood $N$ for the submanifold $w(E)$
of $X\times\RR\times\RR^{d+r}$, so that the projection 
$N\to w(E)$ is a map over $X\times\RR$.  
The forgetful map taking an element $(\pi,\hat f,w,N)$ to
$(\pi,\hat f)$ is a map of sheaves 
\[  h\cW^{(r)}\lra h\cW \]
on $\sX$. This is highly connected if $r$ is large, by Whitney's 
embedding theorem and the tubular neighborhood theorem, so that
the resulting map from $\colim_r\,h\cW^{(r)}$ to $h\cW$ is a weak equivalence 
of sheaves. (The sequential direct limit is formed by sheafifying the 
``naive'' direct limit, which is a presheaf on $\sX$. It is easy 
to verify that passage to representing spaces commutes with sequential
direct limits up to homotopy equivalence.)  
\newline 
Let $\cZ^{(r)}$ be the sheaf taking an $X$ in $\sX$ to the 
set of maps
\[   X\times\RR   \lra \Omega^{d+r}\Th(\U_{r}^{\bot}). \]
Then the representing space of
$\cZ^{(r)}$ approximates $\Omega^{\infty}\bW$, 
that is, $\colim_r\,|\cZ^{(r)}| \simeq\Omega^{\infty}\bW$.
The Pontryagin-Thom collapse construction 
gives us a map of sheaves 
\begin{equation}
\label{eqn-taumap}
\tau^{(r)}\co h\cW^{(r)} \lra \cZ^{(r)}. \ttag
\end{equation} 
In detail: let $(\pi,\hat f,w,N)$ be an element of
$h\cW^{(r)}(X)$, where $\hat f$ is a section of 
$J^2(T^{\pi}E)\to E$, see lemma~\ref{lem-fiberjet}. 
The differential $dw$ determines, for each $z\in E$,
a triple 
$(V_z,\ell_z,q_z)\in\Gr\cW(d+1,r)$. Here $V_z$ is $dw((T^{\pi}E)_z)$, 
viewed as a subspace of the vertical tangent 
space at $w(z)$ of the projection  
\[ 
\begin{array}{ccc}
X\times\RR\times\RR^{d+r} & \lra & X\,, 
\end{array}
\]
which we in turn may identify with $\RR^{d+1+r}$, and $\ell_z+q_z$ 
is the non-constant part of $\hat f(z)$. In particular 
$z\mapsto (V_z,\ell_z,q_z)$ defines a map $\kappa\co 
E\to \Gr\cW(d+1,r)$.
This extends canonically to a pointed map 
\[  \Th(N) \lra  \Th(\U_{r}^{\bot}) \]
because $N$ is identified with $\kappa^*\U_{r}^{\bot}$. But $\Th(N)$ is a 
quotient of $X\times\RR\times S^{d+r}$ where we regard 
$S^{d+r}$ as the one-point compactification of $\RR^{d+r}$.  
Thus we have constructed a map 
\[ X\times\RR\times S^{d+r} \lra \Th(\U_{r}^{\bot}) \]
or equivalently, a map $X\times\RR\lra 
\Omega^{d+r}\Th(\U_{r}^{\bot})$.  
Viewed as an element of $\cZ^{(r)}(X)$, that map is
the image of $(\pi,\hat f,w,N)$ under $\tau^{(r)}$ 
in~(\ref{eqn-taumap}). Taking colimits over $r$, we therefore have a diagram
\[ 
\xymatrix@M=8pt@W=8pt{
|h\cW| & \ar[l]_-{\simeq}   
\rule[-2mm]{0mm}{0mm}\colim_r\,|h\cW^{(r)}| \ar[r] &
\rule[-2mm]{0mm}{0mm} \colim_r\,|\cZ^{(r)}| 
\ar[r]^-{\simeq} & \Omega^{\infty}\bW  
}
\]
which we informally describe as a map $\tau\co |h\cW|\to
\Omega^{\infty}\bW$.  

\begin{thm} 
\label{thm-spaceWtobordism}
The map $\tau\co|h\cW|\to \Omega^{\infty}\bW$ is a homotopy
equivalence. 
\end{thm}

\proof This follows from theorem~\ref{thm-Wtobordism} and 
the fact that $\tau$ can be taken to be a map between spaces with a
monoid structure up to homotopy. Informally, the monoid structure 
on $|h\cW|$ is induced by a monoid structure on $\cW$ itself 
given by ``disjoint union'':
\[  \cW(X)\times \cW(X) \stackrel{\mu}{\lra} \cW(X) \,; 
((\pi,\hat f),(\psi,\hat g))\mapsto 
(\pi\sqcup\psi,\hat f\sqcup \hat g) \]
where the source of $\pi\sqcup\psi$ is the disjoint union of the 
sources of $\pi$ and $\psi$. (See the remark just below.) \newline
To make the monoid structure 
explicit in the case of the target, we introduce 
$\bW\vee \bW$ and the corresponding 
infinite loop space 
\[ \Omega^{\infty}(\bW\vee\bW) = \colim_n\, \Omega^{d+n}
\big(\Th(\U_{n}^{\bot})\vee \Th(\U_{n}^{\bot})\big). \]
The two maps from $\bW\vee \bW$ to $\bW$ which collapse 
one of the two wedge summands lead to a weak equivalence
$\Omega^{\infty}(\bW\vee\bW) \simeq \Omega^{\infty}(\bW)
\times \Omega^{\infty}(\bW)$
and the fold map $\bW\vee\bW\to \bW$ induces an addition map 
from $\Omega^{\infty}(\bW)
\times \Omega^{\infty}(\bW) \simeq 
\Omega^{\infty}(\bW\vee\bW)$ to $\Omega^{\infty}(\bW)$. \newline
It is clear that $\tau$ can be upgraded to respect the
additions. Now theorem~\ref{thm-Wtobordism}
with $X=\pt$ implies that $\tau$ induces a bijection 
\[ \pi_0|h\cW|  \lra \pi_0(\Omega^{\infty}\bW) \]
and consequently that $\pi_0|h\cW|$ is a group, since 
$\pi_0(\Omega^{\infty}\bW)$ is. Next, we use theorem~\ref{thm-Wtobordism}
with $X= S^n$. The monoid structures
imply the isomorphisms
\[
\begin{array}{cccccc} 
\pi_n|h\cW| & \cong & [ S^n,|h\cW|\,]\big/[\pt,|h\cW|\,]\,,\qquad &
\pi_n(\Omega^{\infty}\bW) & \cong & 
[ S^n,\Omega^{\infty}\bW]\big/[\pt,\Omega^{\infty}\bW]   
\end{array}
\]
for arbitrary choices of base points. Thus the map $\tau$ 
induces an isomorphism of homotopy groups, and Whitehead's theorem 
implies that it is a homotopy equivalence, since we are in a 
CW--situation. \qed

\emph{Remark.} To avoid set--theoretical problems related to 
disjoint unions, one should regard $\mu$ in the above proof as a map from 
a certain subsheaf $\cW\bar\times\cW$ of $\cW\times \cW$ to $\cW$. 
An element $((\pi,\hat f),(\psi,\hat g))$ of $(\cW\times\cW)(X)$ 
belongs to $(\cW\bar\times\cW)(X)$ if the sources $E(\pi)$ 
and $E(\psi)$ of $\pi$ and $\psi$, respectively, are disjoint. 
Let $\mu$ take $((\pi,\hat f),(\psi,\hat g))$ to 
$(\pi\cup\psi, \hat f\cup \hat g)$ with 
\[ \pi\cup\psi\co E(\pi)\cup E(\psi) \lra X\,. \] 
Note that the inclusion $\cW\bar\times\cW\lra \cW\times \cW$
is a weak equivalence.

\medskip
The arguments above work in a completely similar fashion to 
identify $|h\cV|$. In fact the map $\tau$ in 
theorem~\ref{thm-spaceWtobordism} restricts to a map 
from $|h\cV|$ to $\Omega^{\infty}\bV$ and the analogue of 
theorem~\ref{thm-Wtobordism} holds. Keeping the letter $\tau$ for 
this restriction, we therefore have 

\begin{thm} 
\label{thm-spaceVtobordism}
The map $\tau\co|h\cV|\to \Omega^{\infty}\bV$ is a homotopy
equivalence.  \qed
\end{thm}

\subsection{The space $|h\cW_{\loc}|$}

We start with a description of $[X,\Omega^{\infty}\bW_{\loc}]$ 
as a bordism group. This is very similar 
to the description of $[X,\Omega^{\infty}\bW]$ used 
in the construction of the map~(\ref{eqn-Wtobordism}).

\begin{lem}
\label{lem-unusualbordism}
 For $X$ in $\sX$, the group 
$[X,\Omega^{\infty}\bW_{\loc}]$ can be identified with the 
group of bordism classes of triples $(M,g,\hat g)$
consisting of a smooth $M$ without 
boundary, $\dim(M)=\dim(X)+d$, and a vector bundle pullback square 
\[
\xymatrix{
TM\times\RR\times\RR^j \ar[r]^-{\hat g} \ar[d] & TX\times
\U_{\infty}\times\RR^j \ar[d] \\
M \ar[r]^-g &  X\times \Gr\cW(d+1,\infty)
}
\]
with $j\gg0$, such that the map 
$g^{-1}(X\times\Sigma(d+1,\infty))\to X$ induced by $g$ is proper. 
\end{lem}

\proof The standard bordism group description of 
the homotopy set $[X,\Omega^{\infty}\bW_{\loc}]$ has representatives which are 
vector bundle pullback squares
\begin{equation*}
\label{eqn-standardbordism}
\xymatrix@+3pt{
TY\times\RR\times\RR^k \ar[r]^-{\hat g^{\rule{0mm}{0mm}}_{Y}} \ar[d] & 
TX\times\Sigma(d+1,\infty)\times\RR^k \ar[d]  \\
Y \ar[r]^-{g^{\rule{0mm}{0mm}}_{Y}}  & X\times \Sigma(d+1,\infty)
} \ttag
\end{equation*}
for some $k\gg0$, where the map $Y\to X$ determined by $g_Y$ is proper, 
$\partial Y=\emptyset$ and $\dim(Y)=\dim(X)-1$. See 
lemma~\ref{lem-bWloctype}. We produce reciprocal maps relating 
this bordism group to the one in lemma~\ref{lem-unusualbordism}.
\newline
We first identify 
$\U_{\infty}|\Sigma(d+1,\infty)$ with its 
dual using the canonical 
quadratic form $q$, 
and then with the normal bundle $N$ of $\Sigma(d+1,\infty)$ in
$\Gr\cW(d+1,\infty)$. Let 
$(M,g,\hat g)$ be a triple as above, lemma~\ref{lem-unusualbordism}.
We may assume that $g$ is transverse to $X\times\Sigma(d+1,\infty)$.
Then $Y=g^{-1}(X\times\Sigma(d+1,\infty))$ is a smooth 
submanifold of $M$, of codimension $d+1$, with normal bundle $N_{Y}$. 
Restriction of $g$ and $\hat g$ yields a
vector bundle pullback square 
\[
\xymatrix@+4pt{
(TY\oplus N_{Y})\times\RR\times\RR^j \ar[r] \ar[d]
& TX\times N\times\RR^j \ar[d] \\
Y \ar[r] &   X\times \Sigma(d+1,\infty)\,.
}
\] 
But since $N_{Y}$ is also identified with the pullback of 
$N$, this amounts to a vector bundle pullback square as 
in~(\ref{eqn-standardbordism}). \newline
Conversely, given data $Y$, $g^{\rule{0mm}{0mm}}_{Y}$ and 
$\hat g^{\rule{0mm}{0mm}}_{Y}$ as 
in~(\ref{eqn-standardbordism}), let $M$ be the (total space of the) 
pullback of $N$ to $Y$. There is a canonical map
from $M$ to $N\subset \Gr\cW(d+1,\infty)$, 
and another from $M$ to $X$, hence a map 
$g\co M\to X\times\Gr\cW(d+1,\infty)$. Moreover $\hat
g^{\rule{0mm}{0mm}}_{Y}$ determines 
the $\hat g$ in a triple $(M,g,\hat g)$ as above. 
It is easy to verify that the two maps of bordism groups
so constructed are well defined and that they are 
reciprocal isomorphisms. 
\qed

\medskip 
We now turn to the construction of a localized version 
of~(\ref{eqn-Wtobordism}), namely, a natural map  
\begin{equation*}
\label{eqn-Wloctobordism}
\tau_{\,\loc}\co h\cW_{\loc}[X] \lra [X,\Omega^{\infty}\bW_{\loc}]. \ttag
\end{equation*}
Let $(\pi,\hat f)\in h\cW_{\loc}(X)$, where $\pi\co E\to X$ 
is a submersion with $(d+1)$-dimensional fibers 
and $\hat f$ is a section of $J^2(T^{\pi}E)\to E$ with 
underlying map $f\co E\to \RR$. See definition~\ref{dfn-hWlocsheaf}
and~\ref{lem-fiberjet}. 
We may assume that $f$ is transverse to $0$ and get a manifold 
$M=f^{-1}(0)$. Proceeding exactly 
as in the construction of the map~(\ref{eqn-Wtobordism}), 
we can promote this to a triple $(M,g,\hat g)$ where $(g,\hat g)$ 
is a vector bundle pullback square  
\[
\xymatrix{
TM\times\RR\times\RR^{j} 
\ar[r]^-{\hat g} \ar[d] & TX\times \U_{\infty}\times\RR^j
\ar[d] \\
M \ar[r]^-g  &  X\times \Gr\cW(d+1,\infty)\,.
}
\]
This time, however, we cannot expect that the $X$-component 
of $g$, which is $\pi|M$, is proper. But its restriction to
\[ g^{-1}(X\times\Sigma(d+1,\infty))=\Sigma(\pi,\hat f)\cap M \]
is proper, thanks to condition $\textup{(ia)}$
in definition~\ref{dfn-hWlocsheaf}. Therefore 
$(M,g,\hat g)$ represents an element in
$[X,\Omega^{\infty}\bW_{\loc}]$. This is the image of 
$(\pi,\hat f)$ under $\tau_{\,\loc}$.

\begin{thm}
\label{thm-Wloctobordism} For compact $X$ in $\sX$, the natural map 
$\tau_{\,\loc}\co h\cW_{\loc}[X]\to [X,\Omega^{\infty}\bW_{\loc}]$ 
is a bijection.
\end{thm}

\proof There is a map $\sigma_{\,\loc}$ in the other direction. 
The construction of $\sigma_{\,\loc}$ is analogous to 
that of $\sigma$ in the proof of theorem~\ref{thm-Wtobordism}.
It is clear that $\tau_{\,\loc}\circ\sigma_{\,\loc}$ is the identity. 
The verification of $\sigma_{\,\loc}\circ\tau_{\,\loc}=\id$ uses 
lemma~\ref{lem-hsingloc} below. \qed 

\begin{lem} 
\label{lem-hsingloc} Let $(\pi,\hat f)\in h\cW_{\loc}(X)$, with 
$\pi\co E\to X$. Let $U$ be an open neighborhood of 
$\Sigma(\pi,\hat f)$ in $E$. Then 
$(\pi|U,\hat f|U)\in h\cW_{\loc}(X)$ is concordant to $(\pi,\hat f)$. 
\end{lem} 

\proof The concordance that we need is an element 
$(\pi^{\sharp},\hat f^{\,\sharp})$ in 
$h\cW_{\loc}(X\times\RR)$. Let 
$E^{\sharp}\subset E\times\RR$ be the 
union of $E\times\,]-\infty,1/2[\,$ and $U\times\RR$.
Let $\pi^{\sharp}(z,t)=(\pi(z),t)$
and $\hat f^{\,\sharp}(z,t)=(\hat f(z),t)$ for $(z,t)\in E^{\sharp}$. 
Some renaming of the 
elements of $E^{\sharp}$ is required to ensure that $\pi^{\sharp}$ 
is graphic. \qed

\bigskip

Next we give a short description of a map
$|h\cW_{\loc}|\to \Omega^{\infty}\bW_{\loc}$ which 
induces~(\ref{eqn-Wloctobordism}). This is analogous 
to the construction of the map named $\tau$ in 
theorem~\ref{thm-spaceWtobordism}. 

Fix an integer $r>0$ and $X$ in $\sX$. To the data $(\pi,\hat f)$ in 
definition~\ref{dfn-hWlocsheaf}, with $\pi\co E\to X$ and $f\co E\to\RR$,
we add the following: a smooth embedding 
\[
\begin{array}{ccc}
 w\co E & \lra & X\times\RR\times\RR^{d+r} 
\end{array}
\] 
which covers $(\pi,f)\co E\to X\times\RR$, a vertical tubular 
neighborhood $N$ for the submanifold $w(E)$
of $X\times\RR\times\RR^{d+r}$, and a smooth function 
$\psi\co E\to [0,1]$ such that $\psi(z)=1$ for all $z\in
\Sigma(\pi,\hat f)$. We require that the restriction of 
$(\pi,f)\co E\to X\times\RR$ to the support of $\psi$ be a 
proper map. \newline
Making $X$ into a variable
now, we can interpret the forgetful map taking $(\pi,\hat f,w,N,\psi)$ to
$(\pi,\hat f)$ as a map of sheaves 
\[  h\cW^{(r)}_{\loc}\lra h\cW_{\loc} \]
on $\sX$. This map is highly connected if $r$ is large.
Let $\cZ^{(r)}_{\loc}$ be the sheaf taking an $X$ in $\sX$ to the 
set of maps 
\[   X\times\RR   \lra 
\Omega^{d+r}\cone\left(\Th(\U_{r}^{\bot}|\Gr\cV(d+1,r))
\hookrightarrow
\Th(\U_{r}^{\bot})\right). \]
Here the cone is a reduced mapping cone, regarded
as a quotient of a subspace of 
\[ \Th(\U_{r}^{\bot})\times[0,1] \,, \]
with
$\Th(\U_{r}^{\bot})\times\{1\}$ corresponding to the base of the cone.  
The Pontryagin-Thom collapse construction 
gives us a map of sheaves 
\begin{equation}
\label{eqn-spaceWloctobordism}
\tau^{(r)}_{\loc}\co h\cW^{(r)}_{\loc} \lra \cZ^{(r)}_{\loc}. \ttag
\end{equation} 
In detail: let $(\pi,\hat f,w,N,\psi)$ be an element of
$h\cW^{(r)}_{\loc}(X)$. We assume that $\hat f$ is a section of 
$J^2(T^{\pi}E)\to E$, see~\ref{lem-fiberjet}. 
The differential $dw$ determines, for each $z\in E$, a triple 
$(V_z,\ell_z,q_z)\in\Gr\cW(d+1,r)$, as in the proof of 
theorem~(\ref{thm-spaceWtobordism}). This gives us a map  
\[\kappa\co E\to \Gr\cW(d+1,r)\times[0,1]\,,\] 
with first coordinate 
determined by $dw$ and second coordinate equal to $\psi$. The map
$\kappa$ fits into a vector bundle pullback square 
\[
\xymatrix{ 
N \ar[r]^-{\hat\kappa} \ar[d] & \U_{r}^{\bot}\times[0,1] \ar[d] \\ 
E \ar[r]^-{\kappa} & \Gr\cW(d+1,r)\times[0,1].
}
\] 
Now we obtain a map 
from $X\times\RR\times S^{d+r}$ to the mapping cone 
\[  \cone\left(\Th(\U_{r}^{\bot}|\Gr\cV(d+1,r))\hookrightarrow
\Th(\U_{r}^{\bot})\right),   \]
viewed as a subquotient of $\Th(\U_{r}^{\bot})\times[0,1]$, 
by $z\mapsto \hat\kappa(z)$ for $z\in N$ and $z\mapsto \pt$ 
for $z\notin N$. It can also be written in the form
\[  X\times\RR  \lra 
\Omega^{d+r}\cone\big(\Th(\U_{r}^{\bot}|\Gr\cV(d+1,r))
\hookrightarrow\Th(\U_{r}^{\bot})\big) \]
so that it is an element of $\cZ^{(r)}_{\loc}(X)$. 
This defines the map $\tau^{(r)}_{\loc}$. Taking colimits over $r$, we 
therefore have a diagram
\[ 
\xymatrix@M=8pt@W=8pt{
|h\cW_{\loc}| & \ar[l]_-{\simeq} 
\rule[-2mm]{0mm}{0mm}\colim_r\,|h\cW^{(r)}_{\loc}| \ar[r] & 
\rule[-2mm]{0mm}{0mm}\colim_r\,|\cZ^{(r)}_{\loc}| 
\ar[r]^-{\simeq} & \Omega^{\infty}\bW_{\loc}  
}
\]
which we informally describe as a map 
$\tau_{\,\loc}\co |h\cW_{\loc}|\to
\Omega^{\infty}\bW_{\loc}$. The following is a straightforward
consequence of theorem~\ref{thm-Wloctobordism} (cf. the proof 
of theorem~\ref{thm-spaceWtobordism}):

\medskip
\begin{thm} 
\label{thm-spaceWloctobordism}
The map $\tau_{\loc}\co|h\cW_{\loc}|\to \Omega^{\infty}\bW_{\loc}$ is
a homotopy equivalence. \qed
\end{thm}

\medskip 
The combination of theorems~\ref{thm-spaceWloctobordism}, 
\ref{thm-spaceWtobordism}, \ref{thm-spaceVtobordism} and 
proposition~\ref{prp-VWexplained} amounts to a proof of 
theorem~\ref{thm-bottomrow} from the introduction. 

\begin{rmk}
\label{rmk-unloved}
{\rm We are left with the task of saying exactly how 
the lower row of diagram~(\ref{eqn-strategicdiagram}) 
should be regarded as a homotopy fiber sequence. 
Define a sheaf $h\cV_{\loc}$ on $\sX$ by copying
definition~\ref{dfn-hVsheaf}, the 
definition of $h\cV$, but leaving out condition (i). 
Then $|h\cV_{\loc}|$ is contractible by an application of 
proposition~\ref{prp-whatitrepresents}. Any choice 
of nullhomotopy for the inclusion $|h\cV|\to |h\cV_{\loc}|$
determines a nullhomotopy for 
$|h\cV|\to |h\cW_{\loc}|$, since $|h\cV_{\loc}|\subset |h\cW_{\loc}|$.
A nullhomotopy for $|h\cV|\to |h\cW_{\loc}|$
constructed like that is understood in
theorem~\ref{thm-bottomrow}.
}
\end{rmk}

\subsection{The space $|\cW_{\loc}|$} 

The goal is to prove theorem~\ref{thm-middlecolumn}, 
i.e., to show that the inclusion of 
$\cW_{\loc}$ in $h\cW_{\loc}$ is a weak equivalence. 
We begin with the observation that the analogue of 
lemma~\ref{lem-hsingloc} holds for $\cW_{\loc}$:  

\begin{lem} 
\label{lem-singloc} Let $(\pi,f)\in \cW_{\loc}(X)$, with 
$\pi\co E\to X$. Let $U$ be an open neighborhood of 
$\Sigma(\pi,f)$ in $E$. Then 
$(\pi|U,f|U)\in \cW_{\loc}(X)$ is concordant to $(\pi,f)$. \qed
\end{lem}

\begin{cor} 
\label{cor-Wloctobordism}
For $X$ in $\sX$, there are natural bijections between 
$\cW_{\loc}[X]$ and either of the two sets below:
\begin{description}
\item[(i)] the set of bordism classes 
of triples $(\Sigma,p,g)$, where $\Sigma$ is a smooth 
manifold without boundary,  
$p\co \Sigma \to X\times\RR$ is a proper smooth map whose 
$X$-coordinate $\Sigma\to X$ is an \'{e}tale map ($=$ local 
diffeomorphism), and 
$g$ is a map from $\Sigma$ to $\Sigma(d+1,\infty)$;
\item[(ii)] the set of bordism classes 
of triples $(\Sigma_0,v,c)$ where $\Sigma_0$ is a smooth manifold 
without boundary, $v\co \Sigma_0\to X$ 
is a proper smooth codimension 1 immersion with oriented normal bundle 
and $c$ is a map from $\Sigma_0$ to $\Sigma(d+1,\infty)$. 
\end{description}
\end{cor}

The bordism relation in both 
cases involves certain maps to $X\times[0,1]$: \'{e}tale maps 
in the case of (i), and codimension one immersions in the case of (ii). 

\proof 
An element $(\pi,f)$ of $\cW_{\loc}(X)$ determines by 
lemma~\ref{lem-singprojetale} 
a triple $(\Sigma,p,g)$ as in (i), where $\Sigma$ is $\Sigma(\pi,f)$ and 
$p(z)=(\pi(z),f(z))$ for $z\in \Sigma\subset E$.
The map $g$ 
classifies the  
vector bundle $T^{\pi}E|\Sigma$ together with the nondegenerate 
quadratic form determined by (one-half) the fiberwise Hessian of $f$.  
Conversely, given a triple $(\Sigma,p,g)$ we can make an element 
$(\pi,f)$ in $\cW_{\loc}(X)$. Namely, let $E\to \Sigma$ be the 
$(d+1)$-dimensional vector bundle classified by $g$, 
with the canonical quadratic form $q\co E\to \RR$.
Let $(\pi,f)\co E\to X\times\RR$ agree with $q+\bar p$,
where $\bar p$ denotes the composition of the vector bundle 
projection $E\to \Sigma$ with $p\co\Sigma\to X\times\RR$.
The resulting maps from $\cW_{\loc}[X]$ to the bordism set in (i), 
and from the bordism set in (i) to $\cW_{\loc}[X]$, are 
inverses of one another: One of the compositions is obviously 
an identity, the other is an identity by lemma~\ref{lem-singloc}.
\newline 
Next we relate the bordism set in (i) to that in (ii). A triple 
$(\Sigma,p,g)$ as in (i) gives rise to a triple $(\Sigma_0,v,c)$ 
as in (ii) provided $p$ is transverse to $X\times 0$. In that 
case we set $\Sigma_0=p^{-1}(X\times 0)$ and define $v$ and $c$ 
as the restrictions of $p$ and $g$, respectively.  
Conversely, a triple $(\Sigma_0,v,c)$ as in (ii) does of course 
determine a triple $(\Sigma,p,g)$ as in (i) with 
$\Sigma=\Sigma_0\times\RR$. The resulting maps from the bordism set in (i) 
to that in (ii), and vice versa, are inverses of one another: 
One of the compositions is obviously 
an identity, the other is an identity by a shrinking lemma 
analogous to (but easier than) lemma~\ref{lem-shrinking2}. 
\qed

\bigskip
It is well-known that the bordism set (ii) in 
corollary~\ref{cor-Wloctobordism} is in natural bijection with 
\[
\begin{array}{ccc}
[X,\Omega^{\infty} S^{1+\infty}(\Sigma(d+1,\infty)_+)] & \cong &
[X,\Omega^{\infty}\bW_{\loc}]. 
\end{array}
\]
Indeed, Pontryagin-Thom theory allows us to represent 
elements of the homotopy set 
$[X,\Omega^{\infty} S^{1+\infty}(\Sigma(d+1,\infty)_+)]$ 
by quadruples $(\Sigma_0,v,\hat v,c)$ where $\Sigma_0$ is smooth 
without boundary, $\dim(\Sigma_0)=\dim(X)-1$, the maps 
$v$ and $\hat v$ constitute a vector bundle pullback square
\[
\xymatrix{
T\Sigma_0\times\RR\times\RR^{j} 
\ar[r]^-{\hat v} \ar[d] & TX\times\RR^j \ar[d] \\
\Sigma_0 \ar[r]^-v &  X 
}
\]
(for some $j\gg0$) 
with proper $v$, and $c$ is any map from $\Sigma_0$ to $\Sigma(d+1,\infty)$. 
By lemma~\ref{lem-destabilization} we can take $j=0$ and by 
immersion theory \cite{Smale59}, \cite{Hirsch59}, 
\cite{Haefliger71}
we can assume $\hat v=dv$, that is, $v$ 
is an immersion and $\hat v$ is its (total) differential.

Consequently $\cW_{\loc}[X]$ is in natural bijection with 
$[X,\Omega^{\infty}\bW_{\loc}]$. It is easy to verify that this 
natural bijection is induced by the composition 
\[ 
\xymatrix@M=6pt{ 
|\cW_{\loc}| \incto[r] &  |h\cW_{\loc}| \ar[r]^-{\tau_{\,\loc}} &  
\Omega^{\infty}\bW_{\loc} 
}
\]
where $\tau_{\,\loc}$ is the map of~(\ref{eqn-spaceWloctobordism}),
(\ref{eqn-Wloctobordism}) and theorem~\ref{thm-spaceWloctobordism}. 
We conclude that the composition is a homotopy equivalence
(cf. the proof of theorem~\ref{thm-spaceWtobordism}).
Since $\tau_{\loc}$ itself is a homotopy equivalence, it follows 
that the inclusion $|\cW_{\loc}|\hookrightarrow |h\cW_{\loc}|$ is a 
homotopy equivalence. This is theorem~\ref{thm-twolocs} from the 
introduction.

\section{Application of Vassiliev's $h$-principle} 
\label{sec-Vassiliev}

This section contains the proof of theorem~\ref{thm-middlecolumn}. 
It is based upon a special case of Vassiliev's \emph{first 
main theorem}, \cite[ch.III]{Vassiliev94} and 
\cite{Vassiliev89}. \newline  
Let $\mathfrak A\subset J^2(\RR^r,\RR)$ denote the space of $2$-jets 
represented by $f\co (\RR^r,z)\to \RR$ with $f(z)=0$, $df(z)=0$ 
and $\det(d^2f(z))=0$, where $d^2f(z)$ denotes the Hessian. 
This set has codimension $r+2$ and is invariant under 
diffeomorphisms $\RR^r\to \RR^r$. \newline
Let $N^r$ be a smooth compact manifold with boundary and let 
$\psi\co N\to \RR$ be a fixed smooth function with 
$j^2\psi(z)\notin \mathfrak A$ for $z$ in a neighborhood of 
the boundary. (Use local coordinates near $z$. The condition means 
that near $\partial N$, all singularities of $\psi$ with 
value 0 are of Morse type, i.e., nondegenerate.) Define spaces  
\[
\begin{array}{rcl}
\Phi(N,\mathfrak A,\psi)&  = & \{ f\in C^{\infty}(N,\RR)  
\mid f=\psi \textup{ near }\partial N, \quad j^2f(z)\notin \mathfrak
A\textup{ for }z\in N\} , \\
h\Phi(N,\mathfrak A,\psi)
&  = & \{\hat f\in \Gamma J^2(N,\RR) 
\mid \hat f=j^2\psi \textup{ near }\partial N, 
\quad \hat f(z)\notin \mathfrak A\textup{ for }z\in N\} , 
\end{array}
\]
where $\Gamma J^2(N,\RR)$ denotes the space of smooth sections 
of the jet bundle $J^2(N,\RR)\to N$. Both are equipped with 
the standard $C^{\infty}$ topology. The special case of 
Vassiliev's theorem that we need is the statement that the map 
\begin{equation*}
\label{eqn-Vassiliev}
j^2\co \Phi(N,\mathfrak A,\psi) \lra h\Phi(N,\mathfrak A,\psi) \ttag
\end{equation*}
induces an isomorphism in cohomology with arbitrary untwisted
coefficients. (Equivalently by the universal coefficient theorem, 
it induces an isomorphism in integral homology.)

\bigskip 
We briefly indicate how~(\ref{eqn-Vassiliev}) relates to the 
jet prolongation map 
$|\cW|\to |h\cW|$ or equivalently (by lemma~\ref{lem-shrinking1})
to the map $|\cW^0|\to |h\cW^0|$. Let $(N,\psi)$ be as above 
with $\dim(N)=d+1$. We assume in addition that $\psi(N)\subset A$ 
and $\psi(\partial N)\subset \partial A$, where $A\subset \RR$ is a 
compact interval with $0\in \intr(A)$, and that $\psi$ is nonsingular 
near $\partial A$. For $X$ in $\sX$, 
let $\cW^0_{\psi}(X)\subset \cW^0(X)$ consist of the pairs 
$(\pi,f)$ as in~\ref{dfn-W0sheaf}, with $\pi\co E\to X$, such that 
$E$ contains an embedded copy of $N\times X$, the map $f$ agrees with 
$(z,x)\mapsto \psi(z)$ on a neighborhood of $\partial N\times X$ in 
$N\times X$, and 
$f^{-1}(0)\subset N\times X$. Restricting $f$ to $N\times X$ 
defines a map from $\cW^0_{\psi}(X)$ to the set of smooth maps 
$X\to \Phi(N,\psi,\mathfrak A)$. Making $X$ into a variable, 
we have a map of sheaves which easily leads to a 
weak homotopy equivalence 
\[ |\cW^0_{\psi}|\simeq 
\Phi(N,\psi,\mathfrak A). \]
Analogous definitions, 
with $\cW^0$ replaced by $h\cW^0$ and $\psi$ by its jet prolongation 
$j^2_{\pi}\psi$, lead to a weak homotopy 
equivalence 
\[ |h\cW^0_{\psi}|\simeq 
h\Phi(N,\psi,\mathfrak A). \] 
Arranging these two homotopy equivalences in a commutative square,
we deduce from~(\ref{eqn-Vassiliev})
that the jet prolongation map 
$|\cW^0_{\psi}|\to |h\cW^0_{\psi}|$ is a homology
equivalence. \newline
Given an element $(\pi,f)\in \cW^0(X)$ with $\pi\co E\to X$, it 
is of course not always possible to find a pair $(N,\psi)$
and an embedding $N\times X\to E$ over $X$ with the good properties
above. However, the problem can 
always be solved locally. Namely, each $x\in X$ has an open 
neighborhood $U$ in $X$ such that $\pi^{-1}(U)$ admits 
such an embedding, $N\times U\to \pi^{-1}(U)$, for suitable $(N,\psi)$.
This fact, its analogue for the sheaf $h\cW^0$ and a general gluing 
technique, developed in section~\ref{subsec-sheaveswithcats}
below, allow us then to conclude that $|\cW^0|\to |h\cW^0|$ 
induces an isomorphism in homology.

\subsection{Sheaves with category structure}
\label{subsec-sheaveswithcats}
Our goal here is to develop an abstract 
gluing principle, summarized in 
proposition~\ref{prp-transportprojection} and 
relying on definition~\ref{dfn-sheafnerve}.
It is a translation into the language of sheaves of something 
which homotopy theorists are very familar with:
the homotopy invariance property of homotopy colimits.
See section~\ref{subsec-hocolim} for background and motivation. 
Since it is relatively easy to reduce the homotopy colimit concept 
to the classifying space construction for categories, our translation 
effort begins with a discussion of sheaves taking values 
in the category of small categories, and a ``classifying sheaf'' 
construction for such sheaves.

\medskip 
Let $\cF\co \sX\to \sC\mathit{at}$ be a sheaf with values in small 
categories. Taking nerves defines a sheaf with values in the category of 
simplicial sets, 
\[ N_{\bullet}\cF\co \sX\to \sS\!\mathit{ets}_{\bullet} \]
with $N_0\cF=\ob(\cF)$ the sheaf of objects 
and $N_1\cF=\mor(\cF)$ the sheaf of morphisms. We have the associated 
bisimplicial set $N_{\bullet}\cF(\Delta^{\bullet}_e)$ and recall
\cite{Quillen73} that the realization of its diagonal 
is homeomorphic to either of its double realizations, 
\begin{equation*}
\label{eqn-simplicialdiagonal}
\begin{array}{ccccc}
|\,\uli k\mapsto N_k\cF(\Delta^k_e)\,| 
& \cong  & \big|\,\uli\ell\mapsto|\,\uli k\mapsto 
   N_k\cF(\Delta^{\ell}_e)\,|\big| & = & 
\big|\,\uli\ell\mapsto B(\cF(\Delta^{\ell}_e))\,\big| \\
\rule{0mm}{5mm} & \cong  & \big|\,\uli k\mapsto|\,\uli\ell\mapsto 
N_k\cF(\Delta^{\ell}_e)\,|\big| & = &
\big|\,\uli k\mapsto |N_k\cF|\big|\,. 
\end{array}
\ttag
\end{equation*}
There is a topological category $|\cF|$ with object space $|N_0\cF|$ 
and morphism space $|N_1\cF|$. (To be quite precise, $|\cF|$ is a 
category object in the category of compactly generated Hausdorff 
spaces.) Since $|N_k\cF|=N_k|\cF|$ by~~\ref{prp-productsheaf}, 
the last of the five expressions in~(\ref{eqn-simplicialdiagonal}) is the 
classifying space $B|\cF|$ of the topological category 
$|\cF|$. 

We next give another construction of $B|\cF|$ related to 
Steenrod's coordinate bundles (i.e., bundles viewed as 1-cocycles). We shall 
consider locally finite open covers $\sY=(Y_j)_{j\in J}$ of spaces 
$X$ in $\sX$, indexed by a \emph{fixed} infinite set $J$. 
The local finiteness condition means that each $x\in X$ 
has a neighborhood $U$ such that $\{j\in J\mid Y_j\cap U\ne \emptyset\}$
is a finite subset of $J$. We use a fixed indexing set 
$J$, independent of $X$ in $\sX$, to ensure good gluing 
properties: suppose that $X$ is the union of two open subsets, 
$X=X'\cup X''$, with intersection $A=X'\cap X''$, and that $(Y'_j)_{j\in J}$
and $(Y''_j)_{j\in J}$ are open coverings of $X'$ and $X''$, 
respectively. The coverings agree on $A$ if $Y'_j\cap A=Y''_j\cap A$ 
for all $j\in J$. In that case, $(Y'_j\cup Y''_j)_{j\in J}$ is an open covering
of $X$ which \emph{induces} the open coverings 
$(Y'_j)_{j\in J}$ and $(Y''_j)_{j\in J}$ of $X'$ and $X''$, respectively. 

For each 
finite nonempty subset $S\subset J$ we write 
\[  Y_S = \bigcap_{j\in S} Y_j\, . \]
Associated to the cover $\sY$  there is a topological category, 
denoted $X_{\sY}$ in \cite[\S4]{Segal68}, with 
\[ \ob(X_{\sY})= \coprod_SY_S\,, \qquad 
\mor(X_{\sY}) = \coprod_R\coprod_{S\supset R} Y_S\,, \]
source map given by the identities $Y_S\to Y_S$ and target map 
given by the inclusions $Y_S\to Y_R$ for $S\supset R$. 
A continuous functor from $X_{\sY}$ to a topological group $G$, 
viewed as a topological category with one object, is equivalent to 
a collection of maps 
\[   \varphi_{RS}\co Y_S\lra G~, \]
one for each pair $R\subset S$ of finite subsets of $J$, subject 
to certain cocycle conditions. The cocycle conditions are 
listed in definition~\ref{dfn-sheafnerve}
below, but in the more general setting where the group of maps 
from $Y_S$ to $G$ is replaced by the category $\cF(Y_S)$.   

\medskip
\begin{dfn}
\label{dfn-sheafnerve}
{\rm 
For $X$ in $\sX$ an element of $\beta\cF(X)$ is a pair
$(\sY,\varphi\bbul)$ where $\sY$ is a locally finite open 
cover of $X$, indexed by $J$, and $\varphi\bbul$ associates to each
pair of finite, nonempty subsets $R\subset S$ of $J$ a morphism 
$\varphi_{RS} \in N_1\cF(Y_S)$ subject to the following 
cocycle conditions: 
\begin{description}
\item[(i)] every $\varphi_{RR}$ is an identity morphism;  
\item[(ii)] for $R\subset S\subset T$, we have 
$\varphi_{RT}=(\varphi_{RS}|Y_T)\circ\varphi_{ST}$. 
\end{description}
Condition (ii) includes the condition that the right-hand 
composition is \emph{defined}; in particular, taking $S=T$ one finds that 
the source of $\varphi_{RS}$ is the object $\varphi_{SS}$, and taking 
$R=S$ one finds that the target of $\varphi_{ST}$ is $\varphi_{SS}|Y_T$. 
}
\end{dfn}

\medskip
The sets $\beta\cF(X)$ define a sheaf 
$\beta\cF\co \sX\to \sS\!\mathit{ets}$ and hence 
a space $|\beta\cF|$. The following key theorem is one of our main
tools used in the proof of both theorem~\ref{thm-middlecolumn} and 
theorem~\ref{thm-toprow}. Its proof is deferred to 
appendix~\ref{sec-spacesboxdetails}. 

\begin{thm}
\label{thm-cocyclesheaf} The spaces $|\beta\cF|$ and $B|\cF|$ are 
homotopy equivalent. 
\end{thm} 

Consider the example where $\cF(X)$ is the set of continuous maps 
from $X$ to a topological group $G$, made into a group by pointwise 
multiplication. An element $(\sY,\varphi\bbul)$ of $\beta\cF(X)$ 
is a collection of gluing data for a principal $G$--bundle 
$P\to X$ with chosen trivializations over each $Y_R$. Namely, 
\[ P= \coprod_R \{R\}\times Y_R\times
G\bigg/\sim \]
where $R$ runs through the finite nonempty subsets of $J$, 
and the equivalence relation identifies $(R,x,g_1)$ with 
$(S,x,g_2)$ if $R\subset S$
and $\varphi_{RS}(x)g_2=g_1$. \newline
The topological category $|\cF|$ is a topological group 
and comes with a continuous homomorphism $|\cF|\to G$ which 
is clearly a weak homotopy equivalence. So $B|\cF|\simeq BG$.
Thus theorem~\ref{thm-cocyclesheaf} reduces to the well--known
statement that concordance classes of principal Steenrod 
$G$--bundles are classified by $BG$. 

\medskip
Consider next the case where $\cF(X)=\map(X,\sC)$ for a small
topological category $\sC$. That is, $\ob(\cF(X))$ and $\mor(\cF(X))$
are the sets of continuous maps from $X$ to $\ob(\sC)$ and 
$\mor(\sC)$, respectively. Then an element of $\beta(\cF(X))$ is   
a covering $\sY$ of $X$ together with a continuous functor from $X_{\sY}$ 
to $\sC$. If $k\mapsto N_k\sC$ is a 
good simplicial space in the sense of \cite{Segal74}, then 
the canonical map $B|\cF|\to B\sC$ is a weak equivalence 
since it is induced by weak equivalences 
$N_k|\cF|\cong|N_k\cF|\to N_k\sC$. Therefore
theorem~\ref{thm-cocyclesheaf} applied to this situation implies that 
homotopy classes of maps $X\to B\sC$ are in natural bijection 
with concordance classes of pairs consisting of a covering $\sY$ 
and a continuous functor from $X_{\sY}$ to $\sC$. This statement
may have folklore status. It appears explicitly in 
lectures given by tom Dieck in 1972, but it seems that tom Dieck
attributes it to Segal. (We are indebted to R. Vogt who kindly sent us
copies of a few pages of lecture notes taken by himself at the time.)

\medskip 
In our applications of theorem~\ref{thm-cocyclesheaf}, the categories 
$\cF(X)$ will typically be partially ordered sets or will have been 
obtained from a functor 
\[  \cF_{\bullet}\co \sC\op \lra \textrm{ sheaves on }\sX\,, \]
where $\sC$ is a small (discrete) category. Given such a functor one 
can define a category  
valued sheaf $\sC\op\ssmallint\cF_{\bullet}$ on $\sX$. Its value 
on a connected manifold $X$ 
is the category whose objects are pairs $(c,\omega)$ with 
$c\in \ob(\sC)$, $\omega\in \cF_c(X)$ and where a morphism 
$(b,\tau)\to (c,\omega)$ is a morphism $f\co b\to c$ in $\sC$ 
with $f^*(\omega)=\tau$. Then 
\[ |\beta(\sC\op\ssmallint\cF_{\bullet})\,|\,\simeq\,\, 
B|\sC\op\ssmallint\cF_{\bullet}|\,\simeq\, \hocolimsub{c\in
  \sC}\,|\cF_c|  \]
(see section~\ref{subsec-hocolim} for details). 

\begin{dfn} 
\label{dfn-earlysheafhocolim}
{\rm The sheaf $\beta(\sC\op\ssmallint\cF_{\bullet})\co \sX
\lra \sC\!\mathit{at}$ will be written $\hocolimsub{c\in
\sC}\,\cF_c$\,. }
\end{dfn}

\medskip
Spelled out, an element of $(\hocolim_c\,
\cF_c)(X)$ consists of 
\begin{itemize}
\item[(i)] a covering $\sY$ of $X$ indexed by $J$, 
\item[(ii)] a functor $\theta$ from the  
poset of pairs $(S,z)$, where $S\subset J$ is finite nonempty 
and $z\in\pi_0(Y_S)$, to $\sC$, 
\item[(iii)] and finally elements 
$\omega_{S,z}\in \cF_{\theta(S,z)}(Y_{S,z})$, where $Y_{S,z}$ denotes 
the connected component of $Y_S$ corresponding to $z\in \pi_0(Y_S)$. 
The elements $\omega_{S,z}$ are related to each other via the maps  
\[  \cF_{\theta(T,z)}(Y_{T,z}) \lra \cF_{\theta(S,\bar z)}(Y_{T,z})
\longleftarrow \cF_{\theta(S,\bar z)}(Y_{S,\bar z}) \]
for each $S\subset T$ and $z\in \pi_0(Y_T)$ with image 
$\bar z\in \pi_0(Y_S)$.
\end{itemize}

\medskip
We close with an application of theorem~\ref{thm-cocyclesheaf} 
which will be used below to extend the special case of Vassiliev's 
theorem mentioned earlier. 

\medskip
\begin{dfn}
\label{dfn-transportprojection}
{\rm Let $\cE,\cF\co \sX\to \sC\!\mathit{at}$ be sheaves and $g\co
\cE\to \cF$ a map between them.
We say that $g$ is a \emph{transport projection}, or that it has 
the unique lifting property for morphisms, if the 
following square is a pullback square of sheaves on $\sX$:
\[ 
\xymatrix{
N_1\cE \ar[r]^{d_0} \ar[d]^g & N_0\cE \ar[d]^g \\
N_1\cF \ar[r]^{d_0} & N_0\cF 
}
\]
where $d_0$ is the source operator. 
}
\end{dfn}

\begin{dfn}
\label{dfn-lifting} 
{\rm A natural transformation $u\co \cF\to\cG$ of 
sheaves on $\sX$ has the \emph{concordance lifting property} if,  
for $X$ in $\sX$ and $s\in\cF(X)$, any concordance 
$h\in \cG(X\times\RR)$ starting at $u(s)$ lifts to a 
concordance $H\in \cF(X\times\RR)$ starting at $s$.
}
\end{dfn} 

Let $g\co \cE\to \cF$ be a map of set-valued sheaves on $\sX$. 
An element $a\in \cF(\pt)$ gives rise to an element again 
denoted $a\in\cF(X)$ for each $X\in \sX$. The \emph{fiber} 
of $g$ over $a$ is the sheaf $\cE_a$ defined by 
\[ \cE_a(X)= \{s\in\cE(X)\mid g(s)=a\}. \]

\begin{prp} 
\label{prp-transportprojection}   
Let $g\co \cE\to \cF$ and $g'\co \cE'\to \cF$ be transport
projections and let 
$u\co \cE\to \cE'$ be a map of sheaves over $\cF$ which 
respects the category structures. Suppose that 
the maps $N_0\cE\to N_0\cF$ and $N_0\cE'\to N_0\cF$ 
obtained from $g$ and $g'$ have the concordance 
lifting property and that, for 
each object $a$ of $\cF(\pt )$, the restriction 
$N_0\cE_a \to N_0\cE'_a$ of $u$ to the fibers 
over $a$ is a weak equivalence (resp. induces an integral homology 
equivalence of the representing spaces). 
Then $\beta u\co \beta\cE\to \beta\cE'$ is a 
weak equivalence (resp. induces an integral homology 
equivalence of the representing spaces). 
\end{prp} 

\proof According to theorem~\ref{thm-cocyclesheaf} it suffices to 
prove that $u$ induces a homotopy (homology) equivalence 
from $B|\cE|$ to $B|\cE'|$. By~(\ref{eqn-simplicialdiagonal}) 
and lemma~\ref{lem-simplicialfibration} it is then also enough to show that 
\[   N_k(u) \co N_k\cE \lra N_k\cE' \]
becomes a homotopy equivalence (homology equivalence) after passage to 
representing spaces, for each $k\ge0$. We note that the simplicial
spaces obtained from a bisimplicial set by realizing in either direction
are good in the sense of~\cite{Segal74}. \newline 
Since $g$ and $g'$ are transport projections, an obvious 
inductive argument shows that, for each $k$, the diagrams 
\[
\xymatrix{
N_k\cE \ar[r] \ar[d]^g & N_0\cE \ar[d]^g \\
N_k\cF \ar[r] & N_0\cF 
}\qquad\qquad
\xymatrix{
N_k\cE' \ar[r] \ar[d]^{g'} & N_0\cE' \ar[d]^{g'} \\
N_k\cF \ar[r] & N_0\cF 
}
\]
are pullback squares. Passage to representing 
spaces therefore turns them into homotopy cartesian squares 
by~\ref{prp-sheavespullback}, since the maps
$N_0\cE\to N_0\cF$ and $N_0\cE'\to N_0\cF$
have the concordance lifting property.
Hence it suffices to consider the case $k=0$, 
\[   N_0u\co N_0\cE\lra N_0\cE'. \]
Again, $N_0\cE\to N_0\cF$ and $N_0\cE'\to N_0\cF$
have the concordance lifting property and $N_0u$ induces a 
weak equivalence (homology equivalence) of the fibers. 
By~\ref{prp-sheavespullback}, 
the fibers turn into homotopy fibers upon passage to 
representing spaces. Consequently $N_0u\co N_0\cE\to N_0\cE'$
is a homotopy equivalence (homology equivalence). \qed

\bigskip
\subsection{Armlets} 
We begin by defining sheaves $\cW^{\sA}$ and $h\cW^{\sA}$ on $\sX$ 
with values in partially ordered sets, and natural 
transformations 
\[
\xymatrix@R=6pt{
 &  \sP\!\mathit{osets} \ar[dddd]^{\textup{forget}} \\
 &  \ar@<-5ex>@{=>}[dd] \\
\sX \ar@/^1pc/[uur]^{\cW^{\sA}} \ar@/_1pc/[ddr]_{\cW^0} & \\
 &  \\
 &  \sS\!\mathit{ets}
}
\qquad\qquad
\xymatrix@R=6pt{
 &  \sP\!\mathit{osets} \ar[dddd]^{\textup{forget}} \\
 &  \ar@<-4ex>@{=>}[dd] \\
\sX \ar@/^1pc/[uur]^{h\cW^{\sA}} \ar@/_1pc/[ddr]_{h\cW^0} & \\
 &  \\
 &  \sS\!\mathit{ets}
}
\]
where $\cW^0$ and $h\cW^0$ are the sheaves introduced in 
section~\ref{subsec-newmodels}, weakly equivalent to 
$\cW$ and $h\cW$, respectively. 

\begin{dfn}
\label{dfn-armlet} {\rm An \emph{armlet} for an element $(\pi,f)\in
\cW^0(X)$ is a compact 
interval $A\subset\RR$ such that $0\in\intr(A)$
and $f$ is fiberwise transverse to the endpoints of $A$.  
}
\end{dfn}

\begin{dfn}
\label{dfn-harmlet} {\rm An \emph{armlet} for an element $(\pi,\hat f)\in
h\cW^0(X)$ is a compact interval $A\subset\RR$ such that $0\in\intr(A)$ and 
\begin{description}
\item[(i)] $f$ is fiberwise transverse to the endpoints of $A$; 
\item[(ii)] $\hat f$ is integrable on an open neighborhood of  
$f^{-1}(\RR\smin\intr(A))$.
\end{description} 
}
\end{dfn}

\medskip
We introduce a partial ordering on elements of $\cW^0(X)$
or $h\cW^0(X)$ equipped with armlets, namely for elements of $\cW^0(X)$:
\[ 
\begin{array}{ccccc}
(\pi,f,A)\le (\pi',f',A') & \textup{ if } & (\pi,f)=(\pi',f') & 
\textup{ and } &  A\subset A'
\end{array}
\]
and similarly for elements of $h\cW^0(X)$.

\medskip
\begin{dfn} {\rm For a connected $X$ in $\sX$ we let $\cW^{\sA}(X)$ 
denote the partially ordered set of elements $(\pi,f,A)$ 
with $A$ an armlet for $(\pi,f)\in \cW^0(X)$. 
Similarly, $h\cW^{\sA}(X)$ is the partially ordered set of 
elements $(\pi,\hat f,A)$ where $(\pi,\hat f)\in 
h\cW^0(X)$ and $A$ is an armlet for $(\pi,\hat f)$. If 
$X$ is not connected we (must) define 
\[
\begin{array}{ccc} 
\cW^{\sA}(X) = \prod_i \cW^{\sA}(X_i)\,, & & 
h\cW^{\sA}(X) = \prod_i h\cW^{\sA}(X_i)  
\end{array} 
\]
where the $X_i$ are the path components of $X$. 
}
\end{dfn}

\bigskip
Any sheaf $\cF\co \sX\to \sS\!\mathit{ets}$ can be considered 
to be a sheaf with a trivial category structure, so
that each $\cF(X)$ is the 
object set of a category which has only identity morphisms. 
In this case an element $(\sY,\varphi\bbul)$ of $\beta\cF(X)$ reduces 
to a pair consisting of a locally finite open covering of $X$, 
indexed by $J$, and a single element $\varphi\in \cF(X)$, namely, 
the unique element restricting to $\varphi_{SS}\in \cF(Y_S)$ for every 
finite nonempty subset $S$ of $J$.
Thus $\beta\cF\cong \beta\!\pt\times \cF$ where $\pt$ denotes 
the terminal sheaf, again viewed as a sheaf with category values. 
In particular there is a forgetful projection 
$\beta\cF \lra \cF$ which is a weak equivalence, since $|\beta\pt|$
is contractible by theorem~\ref{thm-cocyclesheaf}.

\medskip
\begin{prp}
\label{prp-localarmlets}
The forgetful maps $\beta\cW^{\sA}\to \cW^0$ and 
$\beta h\cW^{\sA}\to h\cW^0$ are weak equivalences of sheaves. 
\end{prp}

The proof of proposition~\ref{prp-localarmlets} will be 
broken up into the following three lemmas.  

\medskip
\begin{lem}
\label{lem-neighborhoodlatch}
Let $X$ be in $\sX$ and $(\pi,f)\in \cW^0(X)$.
Every $x\in X$ has an open neighborhood $U$ in $X$ such that 
the image of $(\pi,f)$ in $\cW^0(U)$ admits an armlet. 
\end{lem}

\proof Write $\pi\co E\to X$ and $E_x=\pi^{-1}(x)$. 
By Sard's theorem, we can find numbers $a<0$ and $b>0$ such
that $f_x\co E_x\to\RR$ is transverse to $a$ and $b$ (in other words, 
$a$ and $b$ are regular values of $f_x$). Let $A=[a,b]$. Let 
$C\subset E$ be the closed subset consisting of all $z\in E$ 
where $f$ has a fiberwise singularity and $f(z)=a$ or $f(z)=b$. 
Then $\pi|C$ is proper and so $\pi(C)$
is a closed subset of $X$. Let $U=X\smin \pi(C)$. \qed

\begin{lem}
\label{lem-coverarmlet} With the assumptions of
lemma~\ref{lem-neighborhoodlatch},
there exists an element of $\beta \cW^{\sA}(X)$ 
mapping to $(\pi,f)$ under the forgetful transformation 
$\beta \cW^{\sA}\to \cW^0$. 
\end{lem}

\proof Choose a locally finite covering of $X$ by open subsets 
$Y_j$, where $j\in J$, such that the restriction of $(\pi,f)$ 
to each $Y_j$ admits an armlet $A_j\subset\RR$. For a finite nonempty 
subset $S\subset J$ with nonempty $Y_S$
let $A_S= \bigcap_{j\in S}A_j$. Then $A_S$ is an armlet for the 
restriction of $(\pi,f)$ to $Y_S$. Therefore, given 
nonempty finite $R,S\subset J$ with $R\subset S$ and $Y_S\ne
\emptyset$, we can define $\varphi_{RS}\in N_1\cW^{\sA}(Y_S)$ to be 
the relation 
\[ (\pi,f,A_S)|Y_S \le (\pi,f,A_R)|Y_S . \]  
The data $\varphi_{RS}$ then constitute an element of $\beta \cW^{\sA}(X)$ 
which clearly projects to $(\pi,f)\in \cW^0(X)$. 
\qed

\medskip
It follows from the two previous lemmas 
that the forgetful map 
$\beta \cW^{\sA}[X]\to \cW^0[X]$ is surjective for any $X$ in $\sX$. 
What we really need in order to prove the first half of 
proposition~\ref{prp-localarmlets} is the relative surjectivity as in 
proposition~\ref{prp-relsurjectivity}. This comes from the next lemma, 
in which we assume that our fixed indexing set $J$ is 
\emph{uncountable}. (The assumption is not needed in 
proposition~\ref{prp-localarmlets} because the homotopy type of 
$|\beta...|$ is independent of the cardinality of $J$ as long as $J$
is infinite.)  

\begin{lem}
\label{lem-relativecoverarmlet} 
For $X$ in $\sX$, let $(\pi,f)\in \cW^0(X)$. Let $C$ be a closed subset of $X$ 
and suppose that a
germ of lifts of $(\pi,f)$ across $\beta \cW^{\sA} \lra \cW^0$
has been specified near $C$. Then there exists an 
element in $\beta \cW^{\sA}(X)$ which lifts $(\pi,f)\in \cW(X)$ and  
extends the prescribed germ of lifts near $C$. 
\end{lem}

\proof Let $U$ be a sufficiently small open neighborhood of $C$ in $X$ 
so that the prescribed germ of lifts is represented by an actual lift of 
$(\pi,f)|U$ across $\beta \cW^{\sA}(U) \lra \cW^0(U)$. This gives us 
a locally finite covering $\sY'$ of $U$, and for
each nonempty finite $S\subset J$ and each $z\in\pi_0(Y'_S)$, 
a compact interval $A'_{S,z}\subset \RR$ such that $0\in
\intr(A'_{S,z})$. We have $A'_{S,z}\subset A'_{R,\bar z}$ if $R\subset S$ 
and $\bar z$ is the image of $z$ under $\pi_0(Y'_S)\to \pi_0(Y'_R)$. 
Making $U$ smaller if necessary, we can assume that the covering
$\sY'$ is locally finite in the strong sense that every $x\in X$
has a neighborhood in $X$ which intersects only finitely many 
of the $Y'_j$.  \newline
Now we make a locally finite 
covering of $X$ by open subsets $Y_j$ as follows. 
For $j\in J$ such that $Y'_j$ is nonempty, let $Y_j=Y'_j$.
For all other $j\in J$ (and there are many such since $J$ is uncountable) 
define $Y_j$ in such a way that $Y_j$ avoids a 
fixed neighborhood of $C$ and the restriction of
$(\pi,f)$ to each path component $z\in\pi_0(Y_j)$ admits an armlet
$A_{j,z}$. 
\newline 
It remains to find enough armlets. We need one armlet $A_{S,z}\subset\RR$
for each nonempty finite $S\subset J$ and every component $z\in
\pi_0(Y_S)$. These armlets must satisfy 
$A_{S,z}\subset A_{R,\bar z}$ if $R\subset S$ 
and $\bar z$ is the image of $z$ under $\pi_0(Y_S)\to \pi_0(Y_R)$.
But, reasoning as in the proof of lemma~\ref{lem-coverarmlet}, we find 
that it is enough to say what $A_{S,z}= A_{j,z}$ should be when 
$S$ is a singleton $\{j\}$. We have already said it in the cases 
where $Y_j\ne Y'_j$~; in the other cases we say
$A_{j,z}:=A'_{j,z}$. \qed

\medskip
The proof of the second half of proposition~\ref{prp-localarmlets}
goes like the proof of the first half, except for one additional 
observation which is related to condition (ii) in 
definition~\ref{dfn-harmlet}. For $X$ in $\sX$ let 
$h_c\cW^0(X)$ consist of all $(\pi,\hat f)\in h\cW^0(X)$, 
with $\pi\co E\to X$ etc., such that $\hat f$ is integrable 
on some open $U\subset E$ and $\pi$ restricted to $E\smin U$ is 
proper. 

\begin{lem} The inclusion of sheaves $h_c\cW^0\hookrightarrow h\cW^0$ 
is a weak equivalence. 
\end{lem} 

\proof  Let $(\pi,\hat f)\in h\cW^0(X)$, 
with $\pi\co E\to X$. Choose an open $U\subset E$ such that $\pi$ 
restricted to $E\smin U$ is proper and such that the closure 
of $U$ has empty intersection with $f^{-1}(0)$. Using the convexity 
of the fibers of $J^2_{\pi}(E,\RR)\to E$, especially 
over points $z\in U$, one may deform $\hat f$ (leaving $f$ unchanged) in 
such a way that it becomes integrable on $U$. This shows that
$h_c\cW^0[X]\to h\cW^0[X]$ is surjective. The argument can easily 
be refined to prove a relative statement as in the hypothesis
of proposition~\ref{prp-relsurjectivity}. \qed

\medskip
\subsection{Proof of theorem~\ref{thm-middlecolumn}}
According to lemma~\ref{lem-shrinking1} and 
proposition~\ref{prp-localarmlets} it remains to show that 
\[  j^2_{\pi}\co \beta\cW^{\sA}\to \beta h\cW^{\sA} \]
is a weak equivalence. To this end we introduce a new 
sheaf 
\[  \cT^{\sA} \co \sX \lra \sP\!\mathit{osets}. \]
Suppose given a smooth submersion $\pi\co E\to X$ with
$(d+1)$--dimensional fibers and a
$\Theta$--orientation on $T^{\pi}E$, as in 
definition~\ref{dfn-Vsheaf} and~\ref{dfn-Wsheaf}.
We consider pairs $(\psi, A)$ where $\psi\co E\to \RR$ is a 
smooth function such that $(\pi,\psi)\co E\to X\times\RR$ is proper,
$A\subset \RR$ is a compact interval with $0\in \intr(A)$,
and $\psi$ is fiberwise transverse to $\partial A$. There is no 
restriction on the fiberwise singularities that $\psi$ might have. 

\medskip
\begin{dfn}
\label{dfn-forgetfulsheaf} {\rm For connected $X$ in $\sX$, the set 
$\cT^{\sA}(X)$ consists of triples 
$(\pi,\psi,A)$ as above, modulo the equivalence relation 
which has $(\pi,\psi,A)\sim (\pi,\zeta,A)$ if 
$\psi^{-1}(A)=\zeta^{-1}(A)$ and the support of $\psi-\zeta$ 
is contained in the interior of $\psi^{-1}(A)$.
}
\end{dfn} 

\medskip
As for $\cW^{\sA}$, we get 
$\cT^{\sA}\co \sX\to \sP\mathit{osets}$. Moreover there is 
an obvious commutative diagram of sheaves
\begin{equation*}
\label{eqn-boxcarrier}
\xymatrix@C=12pt{
{\cW^{\sA}} \ar[rr]^{j^2_{\pi}} \ar[dr]_p 
& & {h\cW^{\sA}} \ar[dl]^q \\
& \cT^{\sA} & 
} \ttag
\end{equation*} 
where $p(\pi,f,A)$ and $q(\pi,\hat f,A)$ are
the equivalence classes of $(\pi,f,A)$; in the second case $f$ 
is the underlying function of $\hat f$.

\medskip
Let $(\pi,\psi,A)$ be a representative of an element of 
$\cT^{\sA}(X)$ with $\pi\co E\to X$, $\psi\co E\to \RR$ 
and $A\subset \RR$. The manifold $\psi^{-1}(A)$ 
is independent of the choice of representative for the equivalence 
class, and $\pi|\psi^{-1}(A)$ is a proper submersion, hence 
a smooth fiber bundle by Ehresmann's fibration 
lemma~\cite{BroeckerJaenich73}. Moreover,  
near the boundary $\partial\psi^{-1}(A)=\psi^{-1}(\partial A)$, 
the function $\psi$ is independent of the choice of representative.

\medskip
\begin{lem}
\label{lem-armletfibration}
The maps $p$ and $q$ in~(\ref{eqn-boxcarrier}) have the 
concordance lifting property. 
\end{lem} 

\proof We only give the proof for $p$, since the proof for $q$ is 
much the same. Suppose given a concordance 
$[\pi,\psi,A]\in \cT^{\sA}(X\times\RR)$  
and a lift to $\cW^{\sA}(X\times 0)$ of its restriction to 
$X\times 0$. The projection 
\begin{equation*}
\label{eqn-boxconcordance}
\xymatrix{
\psi^{-1}(A) \ar[r]^{\pi} & X\times \RR
} \ttag
\end{equation*}
is a smooth manifold bundle. Hence there exists a diffeomorphism 
$N\times \RR\cong \psi^{-1}(A)$ over $X\times\RR$, where $N=\psi^{-1}(A)
\cap \pi^{-1}(X\times0)$. But what we need here is a diffeomorphism 
\[  u\co N\times\RR \lra \psi^{-1}(A) \]
over $X\times \RR$ such that $\psi(u(z,t))=\psi(u(z,0))$ for all 
$(z,t)$ near $\partial N\times\RR$, and of course $u(z,0)=z$ for all $z\in N$. 
Constructing such a diffeomorphism $u$ is equivalent to 
constructing a smooth vector field $\xi=du/dt$  
on $\psi^{-1}(A)$ which 
\begin{itemize}
\item[(i)] covers the vector field 
$(x,t)\mapsto (0,1)\in TX_x\times T\RR_t\,$
on $X\times I$, 
\item[(ii)] satisfies $\lag d\psi,\xi\rag\equiv 0$ near
$\psi^{-1}(\partial A)$.
\end{itemize}
(Actually $\xi$ is also 
prescribed on a neighborhood 
of $\psi^{-1}(X\times C)$ where $C=\RR\smin\,]0,1[\,$,   
due to the details in definition~\ref{dfn-concordantsections}.)
This problem has local solutions which can be pieced 
together by means of a partition of unity on
$\psi^{-1}(A)$. Hence $u$ with the required properties 
exists. \newline
Now we define the lifted concordance 
$(\pi,f,A)\in \cW^{\sA}(X\times I)$ in such a way 
that $f(u(z,t))=f(u(z,0))$ for $(z,t)\in N\times \RR$, 
bearing in mind that $f(u(z,0))=f(z)$ is prescribed for all 
$z\in N$ and $f$ must 
equal $\psi$ outside $u(N\times I)=\psi^{-1}(A)$.    \qed

\medskip
\begin{prp}
\label{prp-homologyvassiliev} The fiberwise jet prolongation map 
\[  
j^2_{\pi}\co |\beta\cW^{\sA}| \lra |\beta h\cW^{\sA}|
\]
induces an isomorphism on integral homology. 
\end{prp}

\proof 
This will be deduced from proposition~\ref{prp-transportprojection}
and diagram~(\ref{eqn-boxcarrier}). Both maps $p$ and $q$ 
in~(\ref{eqn-boxcarrier}) are transport projections
in the sense of~\ref{dfn-transportprojection}. 
We must determine the fibers 
of $p$ and $q$ and check that $j^2_{\pi}$ 
induces a homology equivalence between fibers over the same
point. \newline
We first determine the fiber $p^{-1}(\tau)$ of 
\[ p\co \cW^{\sA} \lra \cT^{\sA} \]
over an element $\tau=[F,\psi,A]\in \cT^{\sA}(\pt)$. That is, 
for each $X$ in $\sX$ we are interested in the subset of
$\cW^{\sA}(X)$ which maps to the element 
$[\pi,\psi\circ\pr_F,A]\in \cT^{\sA}(X)$ 
where $\pi$ and $\pr_F$ are the projections $F\times X\to X$ 
and $F\times X\to F$, respectively. This subset consists of 
$(\pi,f,A)\in \cW^{\sA}(X)$ with $\pi$ and $A$ as above, 
where $f\co F\times X\to \RR$ satisfies the conditions 
\begin{itemize}
\item[(i)] 
$\supp(f-\psi\circ\pr_F)\, \subset\, \intr(\psi^{-1}(A))\times X$, 
\item[(ii)] $f(\psi^{-1}(A)\times X)\subset A$.
\end{itemize} 
Because of (i), we can identify 
the fiber of $p$ over 
$\tau$ with a subsheaf of the sheaf taking $X$ in $\sX$ to
the set of smooth maps from $X$ to 
\[ \Phi(\psi^{-1}(A),\mathfrak A,\psi), \] 
using the notation of~(\ref{eqn-Vassiliev}). 
Similarly, the fiber $q^{-1}(\tau)$ of $q$ in ~(\ref{eqn-boxcarrier}) over the 
same element $\tau\in \cT^{\sA}(\pt)$ can be identified with a
subsheaf of the sheaf
taking $X$ in $\sX$ to the set of smooth maps from $X$ to 
\[  h\Phi(\psi^{-1}(A),\mathfrak A,\psi). \]
The inclusions of these subsheaves are weak equivalences by inspection.
(That is to say, condition (ii) means nothing after 
passage to concordance classes.)    
Thus the representing spaces $|p^{-1}(\tau)|$ and 
$|q^{-1}(\tau)|$ have canonical comparison maps to 
$\Phi(\psi^{-1}(A),\mathfrak A,\psi)$ 
and $h\Phi(\psi^{-1}(A),\mathfrak A,\psi)$, respectively, 
which are homotopy equivalences. With these as identifications,
the jet prolongation map from $|p^{-1}(\tau)|$ to $|q^{-1}(\tau)|$
turns into a special case of~(\ref{eqn-Vassiliev}), and so is a 
homology equivalence by Vassiliev's first main theorem.  \qed

\bigskip
Combining 
lemma~\ref{lem-shrinking1}, proposition~\ref{prp-localarmlets}
and proposition~\ref{prp-homologyvassiliev}, we get that 
\[  j^2_{\pi}\co |\cW|\lra |h\cW| \]
induces an isomorphism in homology. Both $|\cW|$ and $|h\cW|$ are
spaces with a monoid structure up to homotopy (cf. the proof 
of theorem~\ref{thm-spaceWtobordism}) and $j^2_{\pi}$ respects 
this additional structure. 
The target $|h\cW|$ is an infinite loop space by 
theorem~\ref{thm-spaceWtobordism}, hence it is group complete. 
(That is, the monoid $\pi_0|h\cW|$ is a group.) 
Since $H_*(j^2_{\pi};\ZZ)$ is an isomorphism,
especially when $*=0$, the source $|\cW|$ is also group complete. 
It is well known that the connected components of a space 
with a group complete 
monoid structure up to homotopy are  
\emph{simple}, and that a map between simple spaces  
is a homology equivalence if and only if it is a homotopy
equivalence. This completes the proof of 
theorem~\ref{thm-middlecolumn}. \qed

\newpage 

\section{Some homotopy colimit decompositions}
\label{sec-strat}

The organization and the main results of this section can be summarized in a 
commutative diagram of sheaves on $\sX$ and maps of sheaves 
\begin{equation*}
\label{eqn-stratificdiagram}
\xymatrix{ 
\cW  \ar[r]  & {\cW_{\loc}} \\
\cW^{\mu}  \ar[u]_{\simeq} \ar[r]  & 
{\cW^{\mu}_{\loc}} \ar[u]_{\simeq} \ar[d]^{\simeq} \\
\cL \ar[r] \ar[u]_{\simeq} & {\cL_{\loc}}    \\
{\hocolimsub{T\textup{ in }\sK}\,\cL_T \ar[u]_{\simeq}}
\ar[d]^{\simeq} \ar@<.8ex>[r] &
{\hocolimsub{T\textup{ in }\sK}\,\cL_{\loc,T}} \ar[u]_{\simeq} 
\ar[d]^{\simeq} \\
{\hocolimsub{T\textup{ in }\sK}\,\cW_T} 
\ar@<.8ex>[r]  & 
{\hocolimsub{T\textup{ in }\sK}\,\cW_{\loc,T}\,.} 
} \ttag
\end{equation*}
The symbol $\simeq$ indicates weak 
equivalences. The homotopy colimits in the 
diagram are homotopy colimits in 
the category of sheaves on $\sX$, as in 
definition~\ref{dfn-earlysheafhocolim}. 
But their representing spaces can be regarded as homotopy colimits in the 
category of spaces according to lemma~\ref{lem-internalhocolim}. 
The top row of diagram~(\ref{eqn-stratificdiagram})
is the inclusion map $\cW \to \cW_{\loc}$. 
The bottom row is what we eventually want to 
substitute for the top row in order to prove 
theorem~\ref{thm-toprow}. 

\medskip 
The following preliminary remarks about~(\ref{eqn-stratificdiagram})
might help the reader through this rather demanding section. 

The elements of $\cW(X)$ and $\cW_{\loc}(X)$ are  
families parametrized by $X$ of $(d+1)$--manifolds $E_x$ equipped
with, among other things, Morse functions $f_x\co E_x\to\RR$.
The same description applies to $\cW^{\mu}(X)$, $\cW^{\mu}_{\loc}(X)$,
$\cL(X)$ and $\cL_{\loc}(X)$ in the second and third row 
of~(\ref{eqn-stratificdiagram}), except that we ask for more structure 
around the critical points. 
In particular, in the important case of $\cL(X)$ we
insist on proper Morse functions $f_x$ whose critical 
points $z\in E_x$ are separately enclosed in certain standard neighborhoods.
Each of these standard neighborhoods 
$N_z\subset E_x\smin \partial E_x$ is a $(d+1)$--manifold with 
boundary; the restricted map $f_x|N_z$ is proper, regular 
on $\partial N_z$ and has no critical 
points in the interior of $N_z$ except $z$. 
These data will enable us later on to 
move the critical values of $f$ up or down, independently 
of each other. \newline
In going from the third row of~(\ref{eqn-stratificdiagram})
to the fourth row, we are adding ``local'' decisions 
which, for each critical point in sight, specify whether the 
corresponding critical value should eventually be moved towards $-\infty$, 
$+\infty$ or $0$. For more precision, suppose that we are dealing 
with a family $(\pi,f)\co E\to X\times\RR$ of $(d+1)$--manifolds 
and proper Morse functions, plus standard neighborhoods for 
the critical points, i.e., an element of $\cL(X)$.
Let $\Sigma(\pi,f)$ be the fiberwise singularity set. Recall
that the projection $\Sigma(\pi,f)\to X$ is \'{e}tale. 
On some connected components (alias sheets) of $\Sigma(\pi,f)$, the map 
$f$ might neither be bounded above nor below. This 
makes a reasonable partition of $\Sigma(\pi,f)$ into a positive, 
a negative and a neutral part \emph{globally} impossible.
But the problem can be solved \emph{locally} in $X$.
Namely, for any $x\in X$ there exist an open neighborhood $U_x$ 
of $x$ in $X$, and a partition of
$\Sigma(\pi,f)\cap \pi^{-1}(U_x)$ into three closed parts: a ``positive'' 
part where $f$ is bounded below, a ``negative'' part where $f$ 
is bounded above, and a ``neutral'' part where $f$ is bounded below 
and above. (The partition is usually not unique.) The neutral part 
will always be a finite covering space of $U_x$ and, by making 
$U_x$ smaller, we can 
assume that it is trivialized, i.e. identified with $T\times U_x$
for a finite set $T$ with some extra structure.
By fixing $T$ and adding these trivialization and partition data to 
the definition of $\cL(X)$ or $\cL_{\loc}(X)$,
we obtain the definitions of $\cL_T(X)$ and $\cL_{\loc,T}(X)$. 
The local existence statement just described translates into 
homotopy colimit decompositions, i.e., the
equivalence between third and fourth row of~(\ref{eqn-stratificdiagram}). 
This should not come as a surprise, since our definition of 
the sheaf-theoretic homotopy colimit, 
definition~\ref{dfn-earlysheafhocolim}, 
involves open coverings and therefore obviously has a ``local'' flavor. 
\newline 
Finally to pass from the fourth row in~(\ref{eqn-stratificdiagram})
to the fifth, we produce concordances which remove critical 
point sheets labelled positive or negative and which 
move the remaining critical values towards $0$. 
By considering a regular level, we are led to 
weak equivalences $\cL_T\simeq \cW_T$ and 
$\cL_{\loc,T}\simeq \cW_{\loc,T}$
where $\cW_T(X)$ and $\cW_{\loc,T}(X)$ are defined in terms of 
bundles of closed $d$--manifolds on $X$ and fiberwise surgery data.

\subsection{Description of main results}
We now give a 
description of the lower row in diagram~(\ref{eqn-stratificdiagram}). 
This begins with a definition of the category $\sK$ by which the homotopy 
colimits are indexed. 

\medskip
\begin{dfn}
\label{dfn-kcategory} {\rm An object of $\sK$ is a finite set 
$S$ equipped with a map to the set $\{0,1,2,\dots,d+1\}$. A morphism 
from $S$ to $T$ is a pair $(k,\ep)$ where $k$ is an injective map, 
over $\{0,1,\dots,d+1\}$, from $S$ to $T$ and $\ep$ is a function 
$T\smin k(S)\to \{-1,+1\}$. The composition of two morphisms
$(k_1,\ep_1)\co S\to T$ and 
$(k_2,\ep_2)\co T\to U$ is 
$(k_2k_1,\ep_3)\co S\to U$ where $\ep_3$ agrees with 
$\ep_2$ outside $k_2(T)$ and with $\ep_1\circ{k_2}^{-1}$ on 
$k_2(T\smin k_1(S))$. 
}
\end{dfn}

\bigskip
\begin{dfn}
\label{dfn-WlocTsheaf} 
{\rm Let $T$ be an object of $\sK$.  
For $X$ in $\sX$, let $\cW_{\loc,T}(X)$ be the set of  
smooth, riemannian $(d+1)$-dimensional vector
bundles $\omega\co V\lra T\times X$ equipped with 
a fiberwise linear isometric involution 
$\rho$ and a $\Theta$--orientation,
all subject to the following conditions. 
\begin{description}
\item[(i)] For $(t,x)\in T\times X$,
the dimension of the fixed point space of $-\rho$ acting on 
the fiber $V_{(t,x)}$  
is equal to the label of $t$ in $\{0,1,\dots,d+1\}$;
\item[(ii)] The composition 
$V\to T\times X\to X$ is a graphic map.
\end{description}
A smooth map $g\co X\to Y$ induces a map 
$\cW_{\loc,T}(Y)\to \cW_{\loc,T}(X)$, given by 
pullback of vector bundles along 
$\id\times g\co T\times X \to T\times Y$. This makes 
$\cW_{\loc,T}$ into a sheaf on $\sX$.
}
\end{dfn} 

\medskip
In definition~\ref{dfn-WlocTsheaf}, 
the involution 
on $V$ leads to an orthogonal vector bundle splitting $V=V^{\rho}\oplus
V^{-\rho}$, where $V^{\rho}$ consists of the vectors fixed by $\rho$
and $V^{-\rho}$ consists of the vectors fixed by $-\rho$. 
We write $D(V^{\rho})$ and $S(V^{-\rho})$ for the 
disk and sphere bundles associated with $V^{\rho}$ and $V^{-\rho}$,
respectively. The vertical tangent bundle of the projection 
\[  
\begin{array}{ccc}
D(V^{\rho})\,\times_{T\times X}\, S(V^{-\rho}) & \lra & X
\end{array} 
\]
inherits a preferred $\Theta$--orientation from $V$,  
described in detail at the end of section~\ref{subsec-regularization}.

\medskip
\begin{dfn} 
\label{dfn-WTsheaf} {\rm For $T$ in $\sK$, a 
sheaf $\cW_T$ on $\sX$ is defined as follows. For
$X$ in $\sX$, an element of $\cW_T(X)$ consists of  
\begin{description}
\item[(i)] a smooth graphic bundle $q\co M\to X$ of  
closed $d$--manifolds, with a $\Theta$--orientation of its fiberwise 
tangent bundle;
\item[(ii)] an element $(V,\rho)$ of $\cW_{\loc,T}(X)$;
\item[(iii)] a smooth embedding over $X$ respecting the
fiberwise tangential $\Theta$--orientations,   
\[ 
\begin{array}{ccc}
e\co  D(V^{\rho})\,\times_{T\times X}\, S(V^{-\rho}) 
& \lra & \quad M \,. 
\end{array}
\]
\end{description} 
} 
\end{dfn}

\bigskip
The sheaves in definitions~\ref{dfn-WlocTsheaf} and~\ref{dfn-WTsheaf}
depend contravariantly on the variable $T$ in $\sK$. This is 
clear in the case of definition~\ref{dfn-WlocTsheaf}: A 
morphism $(k,\ep)\co S\to T$ 
in $\sK$ induces a map from $\cW_{\loc,T}(X)$ to $\cW_{\loc,S}(X)$ 
given by pullback of vector bundles along the map
$k\times\id$ from $S\times X$ to  $T\times X$. (More precisely, 
for $(V,\rho)\in \cW_{\loc,T}(X)$ we define
$(k,\ep)^*(V,\rho)=(V',\rho')$ where $V'$ is the
\emph{restriction} of $V$ to $k(S)\times X$ viewed as a vector 
bundle on $S\times X$, and $\rho'$ is the restriction of $\rho$. 
This ensures that the projection $V'\to X$ is still 
a graphic map.) \newline
The case of definition~\ref{dfn-WTsheaf} is much more interesting. Let 
$(k,\ep)\co S\to T$ be a morphism in $\sK$. 
If $k$ is bijective, there is an obvious identification 
$\cW_T\cong \cW_S$ and this is the induced map. Therefore we may 
assume that $k$ is an inclusion $S\hookrightarrow T$. Then we can 
reduce to the case where $T\smin S$ has exactly one element, $a$. 
This case has two subcases: $\ep(a)=+1$ and $\ep(a)=-1$.

\medskip
\begin{dfn}
\label{dfn-WTsheafinduced1}
{\rm Let $(k,\ep)\co S\to T$ be a morphism in $\sK$ 
where $k$ is an inclusion and  $T\smin S=\{a\}$ 
with $\ep(a)=+1$. We describe the induced map 
\[ \cW_T(X)  \lra \cW_S(X). \] 
Let $(q,V,\rho,e)$ be an element 
of $\cW_T(X)$, with $q\co M\to X$. Map this 
to an element of 
$\cW_S(X)$ by keeping $q\co M\to X$, restricting $V$ to $S\times X$
and restricting $\rho$ and $e$ accordingly. 
}
\end{dfn}

\medskip
\begin{dfn}
\label{dfn-WTsheafinduced2}
{\rm Let $(k,\ep)\co S\to T$ be a morphism in $\sK$ 
where $k$ is an inclusion and  $T\smin S=\{a\}$ 
with $\ep(a)=-1$. For $X$ in $\sX$, the induced map 
\[ \cW_T(X)  \lra \cW_S(X) \] 
is defined as follows. Let $(q,V,\rho,e)$ be an element 
of $\cW_T(X)$, with $q\co M\to X$. Map this 
to the element $(q^{\flat},V^{\flat},\rho^{\flat},e^{\flat})$ of 
$\cW_S(X)$ where 
\begin{enumerate}
\item[(i)] $q^{\flat}\co M^{\flat}\to X$ is the bundle obtained from 
$q\co M\to X$ by fiberwise surgery on the embedded bundle 
of thickened spheres
$e\big(D(V^{\rho}|X_a)\times_{X_a} S(V^{\,-\rho}|X_a)\big)$,
where $X_a$ means $a\times X$;
\item[(ii)] $(V^{\flat},\rho^{\flat})$ is the restriction of $(V,\rho)$ to 
$S\times X$; 
\item[(iii)] $e^{\flat}$ is obtained from $e$ by restriction. 
\end{enumerate}
} 
\end{dfn}
 
\begin{rmk} {\rm For now, the main point is that the fiberwise 
surgery in (i) amounts to removing the interior of the embedded 
thickened sphere bundle and gluing in a copy of 
$D(V^{\,-\rho}|X_a)\times_{X_a} S(V^{\rho}|X_a)$
instead. More details will be given later.  
Note that when $V^{\,-\rho}=0$, the embedded 
thickened sphere bundle whose interior we have to remove is empty.
In this case the fiberwise surgery consist in adding a (disjoint) 
copy of the sphere bundle $S(V)|X_a$ to $M$. 
}
\end{rmk}

\medskip 
There is a forgetful map 
of sheaves $\cW_T\to \cW_{\loc,T}$. It has 
the concordance lifting property, so 
that by corollary~\ref{cor-morerelativesheafhomotopy}, 
the representing spaces of its fibers are the homotopy fibers
of the induced map or representing spaces
\[ |\cW_T|\to |\cW_{\loc,T}|. \]
It is easy to see that the representing space of any fiber 
of $\cW_T\to \cW_{\loc,T}$ is a
classifying space for certain bundles of compact 
$\Theta$--oriented $d$--manifolds with a prescribed boundary.

\subsection{Morse singularities, Hessians and surgeries}
\label{subsec-Morse}

We begin by recalling some well known facts about 
elementary and multi-elementary Morse functions. 
The reader is referred to \cite[ch.I]{Milnor63} and \cite{Milnor65}
for more details in the non-parame\-tri\-zed situation.
By an \emph{elementary} 
Morse function we shall mean a proper smooth map $E\to\RR$ which 
is regular on $\partial E$ and has exactly 
one critical point in $E\smin\partial E$ which is nondegenerate. By a 
\emph{multi-elementary} Morse function we mean a proper smooth map 
$E\to\RR$ which is regular on $\partial E$ and has finitely many 
critical points in $E\smin\partial E$, all nondegenerate and all 
with the same critical value. 

Fix a finite dimensional real vector space $V$ 
with an inner product (i.e., a positive definite bilinear form)
and a linear isometric involution $\rho\co V\to V$. It is convenient 
to have a name for such a triple: we call it a \emph{Morse vector 
space}. The function $f_V\co V\to\RR$ given by 
\begin{equation*}
\label{eqn-standardMorse}
f_V(v)= \lag v,\rho v\rag  
\ttag
\end{equation*}
is a Morse function on $V$ with exactly one critical point.
If we write $V=V^{\rho}\oplus V^{-\rho}$,
then the fomula for $f_V$ becomes 
\[ f_V(v)= \|v_+\|^2-\|v_-\|^2 \]
where $v_+$ and $v_-$ are the components of $v$ in $V^{\rho}$ and
$V^{-\rho}$, respectively. The gradient of $f_V$ on $V$ 
is everywhere perpendicular to the gradient of 
$v\mapsto \|v_+\|^2\|v_-\|^2$, so that the latter function 
is constant on the trajectories of the gradient flow of $f_V$. 
This motivates the following definition. 

\begin{dfn}
\label{dfn-saddle}
{\rm $\saddle(V,\rho)=\{v\in V \,\big|\, 
\|v_+\|^2\|v_-\|^2 \le 1\}$.}
\end{dfn} 

\medskip
If $V^{\rho}=0$ or $V^{-\rho}=0$, then $\saddle(V,\rho)=V$. 
In the remaining cases, the formula 
\begin{equation*}
\label{eqn-longtraceformula}
v\mapsto (\|v_-\|v_+,\,\|v_-\|^{-1}v_-,\,f_V(v)) \ttag
\end{equation*}
defines a smooth embedding of $\saddle(V,\rho)\smin V^{\rho}$ in  
$D(V^{\rho})\times S(V^{-\rho})\times\RR$,
with complement $0\times S(V^{-\rho})\times[0,\infty[$. 
It respects boundaries and it is a map 
over $\RR$, where we use 
the restriction of $f_V$ on 
the source and the function $(x,y,t)\mapsto t$ on 
the target. \newline
Dually, the formula 
\begin{equation*}
\label{eqn-duallongtraceformula}
v\mapsto (\|v_+\|v_-,\,\|v_+\|^{-1}v_+,\,f_V(v)) \ttag
\end{equation*}
defines a smooth embedding of $\saddle(V,\rho)\smin V^{-\rho}$ in  
$D(V^{-\rho})\times S(V^{\rho})\times\RR$,
with complement $0\times S(V^{\rho})\times\,]-\infty,0]$. 
It respects boundaries and it is a map over $\RR$.

\medskip
The map $f_V$ 
in~(\ref{eqn-standardMorse})
restricted to $\saddle(V,\rho)$ is a good local model for elementary Morse 
functions. 
Let $M$ be any smooth compact manifold and let 
\begin{equation*}
\label{eqn-surgerydata}
e\co D(V^{\rho})\times S(V^{-\rho})\to M\smin\partial M  \ttag
\end{equation*}
be a codimension 
zero embedding (``surgery data''), assuming $\dim(V)=\dim(M)+1$.
Then in $M\times\RR$ we have an embedded 
copy of $D(V^{\rho})\times S(V^{-\rho})\times\RR$. We can 
remove its interior and glue in $\saddle(V,\rho)$ instead, 
using formula~(\ref{eqn-longtraceformula}) to identify the 
boundary of $\saddle(V,\rho)$ with the boundary 
of $D(V^{\rho})\times S(V^{-\rho})\times\RR$. The result is a 
smooth manifold $\trace(e)$
of dimension $\dim(M)+1$. More precisely:

\begin{dfn}
\label{dfn-longtrace} {\rm The \emph{long trace} of $e$, denoted 
$\trace(e)$, is the pushout of the two smooth 
codimension zero embeddings 
\begin{equation*} 
\label{eqn-pushoutlongtrace} 
\xymatrix@C+5pt@M+5pt@R-35pt{
\saddle(V,\rho)\smin V^{\rho} 
\ar[rr]^-{(e\times\id)\circ(\ref{eqn-longtraceformula})} & & 
(M\times\RR) \smin \, e(0\times S(V^{-\rho}))\times[0,\infty[\,, \\
\saddle(V,\rho)\smin V^{\rho} \incto[rr] & &  {\saddle(V,\rho).
\rule{39mm}{0mm}}
}
\ttag
\end{equation*}
}
\end{dfn}

\medskip 
For example, if $V^{-\rho}=0$, then $\saddle(V,\rho)=V$ 
and $\saddle(V,\rho)\smin V^{\rho}$ is empty, so that 
$\trace(e)$ becomes the disjoint union of $M\times \RR$ and $V=V^{\rho}$.
Note that $M$ can be empty in this case. 
If $V^{\rho}=0$, then $M$ contains a codimension zero copy 
of $S(V)$. The long trace is obtained by removing  
$S(V)\times[0,\infty[$ from the copy of $S(V)\times\RR$
in $M\times\RR$ and adding a single point instead, so that 
$\trace(e)$ becomes the
disjoint union of $(M\smin\im(e))\times\RR$ and $V=V^{-\rho}$.  

\medskip
The description~\ref{dfn-longtrace}
determines
a structure of smooth manifold 
on $\trace(e)$ and shows that $\trace(e)$ comes with a 
(smooth) elementary Morse 
function, the \emph{height function}, which is the projection to
$\RR$ on the complement of $V^{\rho}$ and equal to $v\mapsto\lag
v,\rho v\rag$ on
the glued-in copy of $\saddle(V,\rho)$. 
The unique critical point is the origin of 
$V^{\rho}\subset \trace(e)$. The corresponding critical value is $0$.
 
\medskip
Conversely, suppose that $N$ is any smooth manifold 
with boundary and $g\co N\to \RR$ is an elementary Morse 
function, with critical value $0$ and unique critical point 
$z\in N\smin\partial N$. Choose a Morse vector space $V$, a
codimension zero embedding $h\co (V,0)\to (N,z)$ and $\delta>0$ such that 
$gh(v)=f_V(v)$ for all $v\in V$ with $\lag v,v\rag<\delta$. 
This is possible by the Morse-Palais lemma; see 
for example \cite{Lang72}. At the price of replacing $g$ by
$3\delta^{-1}g$, we can assume $\delta=3$, so that $gh$ agrees with 
$f_V$ on a neighborhood of 
$D(V^{\rho})\times D(V^{-\rho})$. Now choose a smooth vector 
field $\xi$ on $N$ which extends 
$h_*(\textup{grad}(gh))$ on $h(D(V^{\rho})\times D(V^{-\rho}))$, 
is tangential to $\partial N$ and satisfies 
$\lag dg,\xi\rag >0$ on $N\smin z$. We then 
have a unique smooth embedding $\iota\co \saddle(V,\rho)\to N$ which 
extends $h$ on $ D(V^{\rho})\times D(V^{-\rho})$, maps gradient flow 
trajectories of $f_V$ to flow trajectories of $\xi$ and satisfies 
$g\circ\iota=f_V$ on $\saddle(V,\rho)$. This identifies $N$ with a long 
trace. 

\medskip
The long trace construction has some obvious generalizations.
For example, we can allow simultaneous 
surgeries on a finite number of pairwise disjoint 
thickened spheres. In this case the surgery data consist of  
a finite set $T$, a riemannian vector bundle $V$ on $T$ 
with an isometric involution $\rho$, where $\dim(V)=\dim(M)+1$, 
and a smooth embedding 
\[  e\co D(V^{\rho})\times_T S(V^{\rho}) \lra M\smin\partial M\,. \]
Then $\trace(e)$ is defined as the manifold obtained from $M\times\RR$ 
by deleting the embedded copy of 
\[ D(V^{\rho}_t)\times S(V^{-\rho}_t)\times \RR  \]
for each $t\in T$, and substituting $\saddle(V_t,\rho)$ for it
using formula~(\ref{eqn-longtraceformula}) to do 
the gluing. There is a
canonical height function on $\trace(e)$. It is a Morse 
function with one critical point for each $t\in T$. The only 
critical value is $0$ (if $T\ne \emptyset$).  
\newline
We shall use a parametrized version of the previous construction. 
Let $q\co M\to X$ be a bundle of smooth 
compact $n$-manifolds, let $V\to T\times X$ be a riemannian vector 
bundle of fiber dimension $n+1$ with isometric involution $\rho$, and let 
\[ e\co D(V^{\rho})\times_{T\times X} S(V^{-\rho}) \lra M\smin\partial M \]
be a smooth embeding over $X$. We can regard $e$ as a family 
of embeddings $e_x$ for $x\in X$, each from a disjoint union 
of finitely many thickened spheres 
to a fiber $M_x$ of $q$.   
The manifolds $\trace(e_x)$ for $x\in X$ are the 
fibers of a smooth bundle 
\begin{equation*}
\label{eqn-parametrizedlongtrace}
E\,=\,\trace(e) \lra X\,. \ttag  
\end{equation*} 
It comes equipped 
with a smooth height function $f\co \trace(e) \lra \RR$ which 
is fiberwise Morse; if $T\ne \emptyset$, then the unique critical 
value is 0.

\bigskip 
So far we have looked at ways to create nondegenerate critical 
points, starting with a regular function such as a projection 
$M\times\RR\to \RR$. For us the opposite process, that of removing 
or ``regularizing'' nondegenerate critical points of a Morse function 
$N\to \RR$, will be more important. One approach to this is to go 
through the long trace construction in reverse, but  
this method is unfortunately not very practical when the tangent 
bundle $TN$ carries a $\Theta$--orientation. 
We will therefore use another method. \newline 
Assuming that $f\co N\to \RR$ 
is an elementary Morse function with unique critical value $0$, 
we look for a Morse vector 
space $V$ and a 
codimension zero embedding $\lambda\co \saddle(V,\rho)\to N\smin\partial N$ 
with the property $f\lambda=f_V$. In addition to that, 
we choose a \emph{proper} regular smooth function 
\[ f_V^{+}\co \saddle(V,\rho)\smin V^{\rho} \lra \RR \]
which agrees with $f_V$ on some open subset  
of $\saddle(V,\rho)$ containing the entire boundary and 
the subset $\{w\in\saddle(V,\rho)\mid f_V(w)\le -1\}$. It will be 
shown in a moment that such a function exists and that it is 
essentially unique. Then we let $N^{\reg}=N\smin\lambda(V^{\rho})$
and define 
\[ f^{\reg}\co N^{\reg}\to \RR \]
by $f^{\reg}(x)=f(x)$ for $x\notin \im(\lambda)$ and  
$f^{\reg}(\lambda(w))=f_V^{+}(w)$
for $w\in \saddle(V,\rho)\smin V^{\rho}$. The function $f^{\reg}$ 
is smooth, proper and regular. Any $\Theta$--orientation 
on $TN$ can obviously be restricted to $TN^{\reg}$. 

\medskip
A fairly explicit construction of an $f_V^{+}$ is as follows. 
Choose a diffeomorphism $\psi$ from $\RR$ to $]\!-\infty,0\,[$ such that 
$\psi(t)=t$ for $t<-1/2$. Choose a smooth non-decreasing function 
$\varphi\co [0,1]\to [0,1]$ such that $\varphi(x)=x$ for $x$ 
close to $0$ and $\varphi(x)=1$ for $x$ close to $1$. Let 
\[ \psi_x(t)= \varphi(x)t+(1-\varphi(x))\psi(t) \]
for $x\in[0,1]$. Then $\psi_0=\psi$ embeds $\RR$ in $\RR$ 
with image $]\!-\infty,0\,[$, whereas each $\psi_x$ for $x>0$ 
is a diffeomorphism $\RR\to \RR$. Hence we can define $f_V^{+}$
by 
\begin{equation*}
\label{eqn-smile}
f_V^{+}(v)= \psi_x^{-1}(t) \ttag
\end{equation*}
where $t=f_V(v)$ and $x=\|v_-\|^2\|v_+\|^2$. This choice of $x$ and $t$ 
is suggested by the identification~(\ref{eqn-longtraceformula}).
Similarly, the formula 
\begin{equation*}
\label{eqn-frown}
f_V^{+}(v)= -\psi_x^{-1}(-t) \ttag
\end{equation*}
where $t=f_V(v)$ and $x=\|v_-\|^2\|v_+\|^2$ defines a proper, smooth 
regular function from $\saddle(V,\rho)\smin V^{-\rho}$ to $\RR$. This 
agrees with $f_V$ on a neighborhood of the boundary and on a 
neighborhood of $\{v\in\saddle(V,\rho)\mid f_V(v)\ge +1\}$. 

\medskip
We point out that there are diffeomorphisms 
\begin{equation*}
\label{eqn-newlongtraceformulae}
\begin{array}{ccc}
\sigma_V^+\co \saddle(V,\rho)\smin V^{\rho} & 
\lra & D(V^{\rho})\times S(V^{-\rho})\times\RR  \\
\rule{0mm}{4.5mm}\sigma_V^-\co \saddle(V,\rho)\smin V^{-\rho} & 
\lra & D(V^{-\rho})\times S(V^{\rho})\times\RR
\end{array} \ttag
\end{equation*}
analogous to~(\ref{eqn-longtraceformula}) 
and~(\ref{eqn-duallongtraceformula}). The formulae are 
\[ \sigma_V^+(v)=
(\|v_-\|v_+,\,\|v_-\|^{-1}v_-,\,f^+_V(v))\,,\qquad 
\sigma_V^-(v) = (\|v_+\|v_-,\,\|v_+\|^{-1}v_+,\,f^-_V(v)).
\]

\medskip
\begin{rmk}
\label{rmk-regulartool} {\rm 
Let $V=(V,\rho)$ be a Morse vector space.  
The space of smooth, proper, regular functions on 
$\saddle(V,\rho)\smin V^{\rho}$ which agree with $f_V$ on and near 
$\partial(\saddle(V,\rho))$ and $\{w\in\saddle(V,\rho)\mid
f_V(w)\le -1\}$ is contractible. (This space can be defined as 
the representing space of a sheaf on $\sX$.) We will not use this 
result anywhere and leave the proof to the reader. 
} 
\end{rmk}

\medskip
We finish this section with a useful naturality property 
of $\saddle(V,\rho)$.

\begin{prp}
\label{prp-longtraceintertwine} Suppose given a smooth map 
$e\co\RR\to \RR$ and $a,b\in\RR$ such that 
$e(a)=b$. Assume $0<e'(x)\le 1$ for all $x\in\RR$.   
Then there is a smooth embedding  
$
\tau\co \saddle(V,\rho)
\to \saddle(V,\rho)
$ with $\tau(0)=0$ and $\tau'(0)=\sqrt{e'(a)}\cdot \id_V$
such that 
\[ (f_V+b)\circ\tau=e\circ(f_V+a).
\]
\end{prp}

\proof Without loss of 
generality, $a=b=0$; otherwise replace $e$ by $e_1$ where 
$e_1(x)=e(x+a)-b$, note that $e_1(0)=0$ and that 
$f_V\circ\tau=e_1\circ f_V$ implies 
$(f_V+b)\circ\tau=e\circ(f_V+a)$. Assuming 
$e(0)=0$ therefore, we have to define $\tau$ in such a way that 
$f_V\circ\tau=e\circ f_V$. We remark that $e$ is an orientation 
preserving embedding since $e'(x)>0$ for all $x$. \newline
First define $u\co\RR\to\RR$ by $u(x)=e(x)/x$ for $x\ne 0$ and 
$u(0)=e'(0)$. Then $u$ is smooth, as can be seen from 
\[ e(x)=\int_0^xe'(t)\,\,dt = x\int_0^1 e'(xs)\,\,ds\,. \] 
We have $0<u(x)\le 1$ for $x\in \RR$
and $e(x)=u(x)\cdot x$. Let 
\[ \tau(w) = (u(f_V(w)))^{1/2}w \]
for $w\in \saddle(V,\rho)$.  
Then $f_V(\tau(w))=u(f_V(w))\cdot f_V(w)=e(f_V(w))$, so 
$f_V\circ\tau=e\circ f_V$. It remains to show that $\tau$ is 
an embedding. 
Write $q(w)=(u(f_V(w)))^{1/2}$ so that $\tau(w)=q(w)\cdot w$. 
The product rule gives 
\[ \tau'(w)(h) = (q'(w)(h))\cdot w +q(w)\cdot h  \]
for $h$ in the tangent space $T_wV$. For $w=0$ and $h\ne 0$
the right--hand side is clearly nonzero. For $w\ne 0$ the 
right--hand side can only vanish if $h$ is a scalar multiple of 
$w$. It is therefore enough to try $h=w$. This gives $\tau'(w)(w)$ in the 
left--hand side, which is the derivative of $t\mapsto \tau(tw)$
at $1\in\RR$. If this vanishes, then the derivative of 
\[ t \mapsto f_V(\tau(tw)) \]
at $1\in \RR$ also vanishes. But $f_V(\tau(tw))=e(f_V(tw))$, 
and since $e'$ is everywhere nonzero, 
it follows that $f'_V(w)(w)=0$ by the chain 
rule. Since $f_V$ is a quadratic form, this forces 
$f_V(w)=f_V(tw)=0$. But then $\tau(tw)=(u(0))^{1/2}tw$ which, 
as a function of $t$, certainly has a nonzero derivative at $1\in\RR$, 
contradiction. Hence $\tau'(w)$ is invertible for every $w$. 
Since $\tau$ also maps each line segment through $0\in V$ to itself, 
it follows immediately that $\tau$ is an embedding. \qed

\subsection{Right--hand column}
In this and the next sections we shall use vector bundles $V\to Y$
equipped with a fiberwise inner product $\lag\,,\,\rag$ and a
fiberwise linear isometric involution $\rho\co V\to V$. We call a 
vector bundle with this additional structure  
a \emph{Morse vector bundle}.

Our most important examples of Morse vector bundles are as follows. 
Let $(\pi,f)$ be an element of $\cW_{\loc}(X)$, with $\pi\co E\to X$. 
The restriction of the vertical tangent bundle $T^{\pi}E$ to the 
fiberwise singularity set 
$\Sigma=\Sigma(\pi,f)$ comes with an everywhere nondegenerate 
symmetric bilinear form $\frac{1}{2}H$, where $H$ is the vertical Hessian of
$f$, that is, the second derivative in the fiber direction. 
See \cite[I,\S2]{Milnor63}. We can choose
an orthogonal direct sum decomposition of $T^{\pi}E|\Sigma$ into 
a positive definite subbundle and a negative definite subbundle. 
(The choice is usually not unique, but the space of all such choices 
is contractible.) By changing 
the sign of $\frac{1}{2}H$ on the negative definite subbundle, 
we make $T^{\pi}E|\Sigma$ into a Morse vector 
bundle, with an isometric involution which is $-\id$ on the preferred 
negative definite summand and $+\id$ on the positive definite summand.
Note in addition that $\pi|\Sigma$ is an \'{e}tale map $\Sigma\to X$ 
and that the restriction of $(\pi,f)$ to $\Sigma$ is a proper map 
from $\Sigma$ to $X\times\RR$.

\begin{dfn}
\label{dfn-LWlocsheaf} 
 {\rm Let $\cL_{\loc}$ be the following sheaf on $\sX$. For $X$ in
$\sX$, an element of $\cL_{\loc}(X)$ is a triple $(p,g,V)$ where
\begin{description}
\item[(i)] $p$ is a graphic and \'{e}tale map 
from some smooth $Y$ to $X$~;
\item[(ii)] $g$ is a smooth function $Y\to\RR$~;
\item[(iii)] $V\stackrel{\omega}{\lra}Y$ is a $(d+1)$-dimensional  
Morse vector bundle with a $\Theta$--orientation.
\end{description}
\emph{Conditions}: The map $(p,g)\co Y\to X\times\RR$ is proper and
$p\omega\co V\to X$ is a graphic map.  
}
\end{dfn}

\begin{dfn} 
\label{dfn-LWmulocsheaf}
{\rm An element of $\cW^{\mu}_{\loc}(X)$ consists of an element $(\pi,f)$ of
$\cW_{\loc}(X)$ with $\pi\co E\to X$, an element $(p,g,V)$ 
of $\cL_{\loc}(X)$ with $p\co Y\to X$, and an isomorphism over $X$
of the vector bundle $V\to Y$ with the vector bundle 
$T^{\pi}E|\Sigma(\pi,f)\to \Sigma(\pi,f)$. \emph{Condition:} 
the vector bundle isomorphism preserves the $\Theta$--orientations 
and carries the function $f_V$ on $V$ to $w\mapsto \frac{1}{2}H(w,w)$
on $T^{\pi}E|\Sigma(\pi,f)$. 
}
\end{dfn}

\begin{lem} The forgetful map $\cW^{\mu}_{\loc}\to \cW_{\loc}$ is a 
weak equivalence. 
\end{lem} 

\proof This is a straightforward application of 
proposition~\ref{prp-relsurjectivity}. \qed

\medskip 
There is also a forgetful map $\cW^{\mu}_{\loc}\to \cL_{\loc}$. 
We now describe a homotopy inverse for this.  
Fix $X$ in $\sX$ and let 
$(p,g,V)$ be an element of $\cL_{\loc}(X)$, with $p\co Y\to X$.
Let $E=V$ and let $\pi\co E\to \RR$ 
agree with the composition $V\to Y\to X$.
Let $f\co E\to \RR$ be given by 
\begin{equation*}
\begin{array}{ccc}
& f(v) = g(y)+ f_V(v), & \qquad\qquad f_V(v)=\lag v,\rho v\rag  
\end{array} \ttag
\end{equation*}
for $y\in Y$ and $v$ in the fiber of $V$ over $y$. Then $(\pi,f)$ 
is an element of $\cW^{\mu}_{\loc}(X)$ with 
$\Sigma(\pi,f)\cong Y$. There are obvious identifications 
of $\Sigma(\pi,f)$ with $Y$ and of $T^{\pi}E|\Sigma(\pi,f)$ with 
$V$. The rule $(p,g,V)\mapsto (\pi,f,p,g,V,...)$
is therefore a map $\cL_{\loc}(X)\to \cW^{\mu}_{\loc}(X)$. 

\begin{prp} 
\label{prp-easylocal} The map 
$\cL_{\loc} \to \cW^{\mu}_{\loc}$ so defined is a weak 
equivalence; consequently the forgetful map $\cW^{\mu}_{\loc}\to \cL_{\loc}$
of is also a weak equivalence.  
\end{prp}

\proof We are going to use the relative surjectivity criterion
of proposition~\ref{prp-relsurjectivity}. To deal with the 
absolute case first, we assume given $X$ in $\sX$ and 
$(\pi,f)\in \cW^{\mu}_{\loc}(X)$, with $\pi\co E\to \RR$ and $f\co E\to\RR$. 
Let $\Sigma=\Sigma(\pi,f)$ be the fiberwise singularity set of 
$f$. Choose a vertical tubular neighborhood $V$ of $\Sigma$ in $E$ 
(see definition~\ref{dfn-verticaltube}
and lemma~\ref{lem-singprojetale}). Note that 
$V$ is identified with the normal bundle of $\Sigma$ in $E$, which 
is identified with $T^{\pi}E|\Sigma$; so $V$ has the structure 
of a Morse vector bundle. In turn, by proposition~\ref{lem-singloc}, 
the element $(\pi,f)$ in $\cW^{\mu}_{\loc}(X)$ is concordant to 
$(\pi^{(1)},f^{(1)})$
where $\pi^{(1)}$ and $f^{(1)}$ are the restrictions of $\pi$ and 
$f$ to $V$, respectively.  
The next step is to improve $f^{(1)}$. \newline
Let $\psi\co \RR\to [0,1]$ be a smooth non-increasing function such that 
$\psi(t)=1$ for $t<1+\ep$ and $\psi(t)=0$ for $t>2-\ep$, for 
some small $\ep>0$. For $t\in \RR$ let $f^{(t)}$ be
given by 
\[ v\mapsto \left\{\begin{array}{cc} 
fp(v)+ \psi(t)^{-2}(f(\psi(t)v)-fp(v)) & \quad\textup{for $\psi(t)>0$
and $v\in V$}  \\ 
fp(v)+\frac{1}{2}H(pv)(v,v) & \quad\textup{for $\psi(t)=0$ and
$v\in V$}
\end{array} \right. 
\]
where $H(pv)$ denotes the vertical Hessian of $f$ at $p(v)$. 
Let $\pi^{(t)}=\pi^{(1)}$ for convenience. 
Then $t\mapsto (\pi^{(t)},f^{(t)})$ defines a concordance
from $(\pi^{(1)},f^{(1)})$ 
to $(\pi^{(2)},f^{(2)})$. Now $(\pi^{(2)},f^{(2)})$ has a canonical 
lift  
to an element of $\cL_{\loc}(X)$, since $V$ is a Morse vector
bundle. We have now 
established the absolute case of the relative surjectivity 
condition of~\ref{prp-relsurjectivity} for our map 
$\cL_{\loc}\to \cW_{\loc}$. The relative case is not much more 
difficult and we leave it to the reader. \qed  

\medskip
We next come to the homotopy colimit decompositions of the 
right--hand column of~(\ref{eqn-stratificdiagram}), based on 
the following key observation. 

\begin{lem} 
\label{lem-localbounds}
Let $(p,g,V)\in \cL_{\loc}(X)$, with $p\co Y\to X$.
For every $x\in X$ and every $b>0$ 
there exist a neighborhood $U$ of $x$ in $X$ 
such that, on every component of $p^{-1}(U)$, the function 
$g$ is either bounded below by $-b$ 
or bounded above by $b$. 
\end{lem}

\proof Chose a descending sequence of open 
balls $U_i$ for $i=0,1,2,3,\dots$ 
forming a neighborhood basis for $x$ in $X$. 
If the statement is false, then there exists $b>0$ and  
connected subsets $K_i\subset Y$ 
for $i=0,1,2,3,\dots$ such that $p(K_i)\subset U_i$
and $g(K_i)\supset[-b,b]$ for all $i$. 
Choose $z_i\in K_i$ such that $g(z_i)=0$. 
The sequence $z_0,z_1,z_2,\dots$ in $Y$ 
must have a convergent (infinite) subsequence, 
because $(p,g)\co Y\to X\times\RR$ is proper and 
the two image sequences in $X$ and $\RR$
converge. Let $z_{\infty}\in Y$ be the point which the 
subsequence converges to. Then $p(z_{\infty})=x$ and
$g(z_{\infty})=0$. Now $p\co Y\to X$
is \'{e}tale. Hence, for sufficiently large $i$, 
there are unique neighborhoods $U_i'$ of $z_{\infty}$ in $Y$ such that 
$p$ maps $U_i'$ diffeomorphically to $U_i$. It follows
that $z_i\in U'_i$ for infinitely many $i$ and hence $K_i\subset U'_i$ 
for infinitely many $i$. But it is also clear that 
the diameter of $g(U'_i)$ tends to zero as $i$ tends to 
infinity; hence the lim inf of the diameters of the intervals 
$g(K_i)$ is zero, which contradicts our 
assumption. \qed

\begin{dfn} 
\label{dfn-LWlocbetterstrata}
{\rm Fix $S$ in $\sK$. We define a sheaf 
$\cL_{\loc,S}$ on $\sX$. For $X$ in $\sX$, an element of 
$\cL_{\loc,S}(X)$ is an element 
$(p,g,V)$ of $\cL_{\loc}(X)$, where $p$ has source $Y$, together with 
a continuous function $\delta\co Y\lra \{-1,0,+1\}$,
and a diffeomorphism 
\[ h\co S\times X\lra \delta^{-1}(0)\subset Y \]
over $\{0,1,\dots,d+1\}\times X$.   
\emph{Condition}: Every $x\in X$ has a neighborhood $U$ in $X$ 
such that $g$ admits a lower bound on $p^{-1}(U)\cap\delta^{-1}(+1)$
and an upper bound on $p^{-1}(U)\cap \delta^{-1}(-1)$. 
}  
\end{dfn}

\medskip
In definition~\ref{dfn-LWlocbetterstrata}, the function 
$\delta$ clearly has to be constant on each component of $Y$.
Note that the Morse vector bundle structure on $V\to Y$ 
determines a map $Y\to \{0,1,\dots,d+1\}$ given by the Morse 
index: $y\mapsto \dim(V_y^{-\rho})$. 
This is what we mean when referring to $Y$ as a space
over $\{0,,1,\dots,d+1\}\times X$. The local 
boundedness condition on $g$ has an alternative ``global'' 
formulation as follows: a continuous function $b\co X \to \RR$
must exist such that $-bp\le g$ on $\delta^{-1}(+1)$ 
and $bp\ge g$ on $\delta^{-1}(-1)$. (We do not ask for a 
constant bound $b$ because we need to ensure that $\cL_{\loc,S}$
is a sheaf.)

\medskip
A morphism $(k,\ep)\co R\to S$ in $\sK$ induces a map
$\cL_{\loc,S}\to 
\cL_{\loc,R}$ 
taking an element $(p,g,V,\delta,h)$ of $\cL_{\loc,S}(X)$ to
$(p,g,V',\delta',h')$
where $V'$ is obtained from $V$ by pulling back, 
$h'(r,x)=h(k(r),x)$ for $(r,x)\in R\times X$ and 
\begin{equation*}
\label{eqn-inducedmapLloc}
\delta'(y)= \left\{\begin{array}{cl}
\ep(s)  & \textup{if $y=h(s,x)$ where }s\in S\smin k(R),\,x\in X \\
\delta(y) & \textup{otherwise}.\\ 
\end{array} 
\right. 
\ttag
\end{equation*}
This makes the rule $T\mapsto \cL_{\loc,T}$ into a contravariant functor 
from $\sK$ to the category of sheaves on $\sX$. Moreover, for each 
$T$ in $\sK$ there is a forgetful map $\cL_{\loc,T}\to \cL_{\loc}$,
and the maps 
$\cL_{\loc,T}\to \cL_{\loc,S}$ induced by morphisms $S\to T$ in $\sK$ 
are over $\cL_{\loc}$. This leads to a canonical map of sheaves 
\begin{equation*}
\label{eqn-middlerightcolumn}
v\co \hocolimsub{T\textup{ in }\sK} \cL_{\loc,T}
\quad \lra \quad\cL_{\loc}\,.  \ttag 
\end{equation*}

\medskip
\begin{prp}
\label{prp-tamelifting} 
The map $v$ in~(\ref{eqn-middlerightcolumn}) is a weak equivalence. 
\end{prp}

\proof Let $\cL^{\delta}_{\loc}$ be the following sheaf on $\sX$
with category structure.
An object of 
$\cL^{\delta}_{\loc}(X)$ is an element $(p,g,V)$ of $\cL_{\loc}(X)$,
with $p\co Y\to X$, together with a continuous function 
$\delta\co Y\to\{-1,0,+1\}$ subject to the following condition: 
\begin{itemize}
\item[{}]
Every $x\in X$ 
has a neighborhood $U$ in $X$ such that $g$ admits a lower 
bound on $p^{-1}(U)\cap \delta^{-1}(+1)$, an upper bound on 
$p^{-1}(U)\cap \delta^{-1}(-1)$, and both an upper and a lower 
bound on $p^{-1}(U)\cap \delta^{-1}(0)$. 
\end{itemize}
Given two such objects, $(p,g,V,\delta_a)$ and $(p,g,V,\delta_b)$
with the same underlying $(p,g,V)$, we write 
$(p,g,V,\delta_a)\le (p,g,V,\delta_b)$ if $\delta_a^{-1}(+1)
\subset \delta_b^{-1}(+1)$ and $\delta_a^{-1}(-1)
\subset \delta_b^{-1}(-1)$. In this situation there is 
a unique morphism from $(p,g,V,\delta_a)$ to $(p,g,V,\delta_b)$, 
otherwise there is none. Thus the category 
$\cL^{\delta}_{\loc}(X)$ is a poset. \newline 
The map $v$ in~(\ref{eqn-middlerightcolumn}) can now be factorized 
as follows: 
\begin{equation*}
\label{eqn-factorofv}
\xymatrix{
{\rule{0mm}{6mm}\hocolimsub{T\textup{ in }\sK} \cL_{\loc,T}} \ar[r]^-{v_1} &
{\beta\cL^{\delta}_{\loc}} \ar[r]^{v_2} & {\cL_{\loc}}
} \ttag
\end{equation*}
Here $v_2$ is induced by the forgetful map 
$\cL^{\delta}_{\loc}\to {\cL_{\loc}}$. (Compare 
proposition~\ref{prp-localarmlets}.) To describe $v_1$ 
we recall that $\hocolim_T\,\cL_{\loc,T}$ was defined as 
\[\begin{array}{c}
\beta(\sK\op\ssmallint\cL_{\loc,\bullet}).
\end{array}\]
An object in $(\sK\op\ssmallint\cL_{\loc,\bullet})(X)$
consists of an object $T$ in $\sK$ and an element $a$ in
$\cL_{\loc,T}(X)$. A morphism from $(T,a)$ to $(S,b)$ is a 
morphism $S\to T$ in $\sK$ taking $a$ to $b$. An object 
$(T,a)$ in $(\sK\op\ssmallint\cL_{\loc,\bullet})(X)$ with $a=(p,g,V,\delta,h)$
determines an object $(p,g,V,\delta)$
in $\cL^{\delta}_{\loc}(X)$.  
This canonical association is a functor, for each $X$, 
and as such induces $v_1$. The next two lemmas complete the
proof. \qed

\begin{lem} The map $v_1$ of (\ref{eqn-factorofv}) is a weak
equivalence. \end{lem} 

\proof We show that the functor 
$(\sK\op\ssmallint\cL_{\loc,\bullet})(X) \lra \cL^{\delta}_{\loc}(X)$
is an equivalence of categories when $X$ is simply connected. 
Indeed for an object 
$(p,g,V,\delta)$ of $\cL^{\delta}_{\loc}(X)$, 
the subset $Y_0=\delta^{-1}(0)$ of $Y$ is closed 
and $g\co Y_0\to \RR$ is locally bounded. Thus $p\co Y_0\to X$ 
is a proper \'{e}tale map, hence a covering. Since $X$ is simply 
connected there is a  
diffeomorphism $h\co S\times X\to Y_0$, giving an object of 
$(\sK\op\ssmallint\cL_{\loc,\bullet})(X)$. \newline   
In particular, we have an equivalence of categories for the 
extended simplices, $X=\Delta^k_e$ where $k\ge 0$. It follows that 
$|\sK\op\ssmallint\cL_{\loc,\bullet})| \lra |\cL^{\delta}_{\loc}|$
is a weak homotopy equivalence, cf. 
section~\ref{subsec-sheaveswithcats}. \qed

\begin{lem} The map $v_2$ of (\ref{eqn-factorofv}) is a weak
equivalence. \end{lem}

\proof The proof is completely analogous to the proof of 
proposition~\ref{prp-localarmlets}. We note that  
given objects $(p,g,V,\delta_1)$ and
$(p,g,V,\delta_2)$ in $\cL^{\delta}_{\loc}(X)$ with the same
underlying $(p,g,V)\in \cL_{\loc}(X)$, there always 
exists an object $(p,g,V,\delta_3)$ in $\cL^{\delta}_{\loc}(X)$ such that 
\[ 
\begin{array}{ccc}
(p,g,V,\delta_3) & \le & (p,g,V,\delta_1) \\
(p,g,V,\delta_3) & \le & (p,g,V,\delta_2). 
\end{array}
\]
Namely, let $\delta_3(z)=+1$ if and only if 
$\delta_1(z)=+1=\delta_2(z)$; let $\delta_3(z)=-1$ if and only if 
$\delta_1(z)=-1=\delta_2(z)$, and let $\delta_3(z)=0$ in the 
remaining cases. \newline
Now we apply proposition~\ref{prp-relsurjectivity} to $v_2$. 
Given $(p,g,V)\in \cL_{\loc}(X)$, we can by lemma~\ref{lem-localbounds}
find a locally 
finite covering of $X$ by open subsets $U_j$, where $j\in J$, 
such that $(p,g,V)\,| U$ has a lift $\varphi_{jj}$ to 
$\ob(\cL^{\delta}_{\loc})(U_j)$
for all $j$. With the observation just above, it is 
easy to extend the collection of the $\varphi_{jj}$ to a collection of
objects $\varphi_{RR}\in \ob(\cL^{\delta}_{\loc})(U_R)$, in such a way 
that $\varphi_{RR} \le \varphi_{QQ}|U_R$ 
whenever $Q\subset R$. The collection of these $\varphi_{RR}$ is then an
element of $\beta\cL^{\delta}_{\loc}(X)$. This establishes the 
absolute case of the hypothesis in ~\ref{prp-relsurjectivity}, 
and the verification 
is much the same in the relative case. \qed

\begin{dfn}
\label{dfn-LlocTtoWlocT}  
 {\rm Fix $T$ in $\sK$. We define a map from 
$\cL_{\loc,T}$ to $\cW_{\loc,T}$ by 
\[ \cL_{\loc,T}(X) \ni(p,g,V,\delta,h)\quad\mapsto\quad h^*(V)\in 
\cW_{\loc,T}(X) . \]
There is an equally simple map in the other direction,
$\cW_{\loc,T}\to \cL_{\loc,T}$. Indeed,
we can identify $\cW_{\loc,T}(X)$ with 
the subset of $\cL_{\loc,T}(X)$ consisting of the elements 
$(p,g,V,\delta,h)\in \cL_{\loc,T}(X)$ which have $h=\id_{T\times X}$
and $\delta\equiv 0$, $g\equiv 0$. 
}
\end{dfn}

\begin{lem}
\label{lem-stratabeauty}
The inclusion $\cW_{\loc,T}\to \cL_{\loc,T}$ is a weak equivalence. 
\end{lem}

\proof We use proposition~\ref{prp-relsurjectivity}.  
Given $(p,g,V,\delta,h)\in \cL_{\loc,T}(X)$ 
with $p\co Y\to X$, choose 
a smooth $\psi\co [0,1/2[\, \to [0,\infty[\,$ such that $\psi(s)=0$ 
for $s$ close to $0$ and $\psi(s)$ tends to $+\infty$ for $s\to 1/2$.
Choose another smooth $\varphi\co [0,1]\to[0,1]$ such that 
$\varphi(s)=1$ for $s$ close to 0 and $\varphi(s)=0$ for $s$ close to 1.   
Then define a concordance
\[(\bar p,\bar g,\bar V,\bar\delta,\bar h)\in
\cL_{\loc,T}(X\times\RR) \]
in the following way. The source of $\bar p$
is the union of $Y\times \,]-\infty,1/2[\,$ and 
$h(T\times X)\times\,]0,\infty[\,$. The
formula for $\bar p$ is $\bar p(y,s)=(p(y),s)$. 
(To ensure that $\bar p$ 
is graphic, we should define the 
source of $\bar p$ and $\bar g$ 
as a subset of the pullback of $p\co Y\to X$ along the 
projection $X\times\,]0,1\,[\,\lra X$. 
See definition~\ref{dfn-graphicmaps}.) The formula for 
$\bar g$ is $\bar g(y,s):=g(y)\cdot\varphi(s)$ if $y$ is in $h(T\times X)$ and 
$\bar g(y,s):= g(y)+\delta(y)\psi(s)$ otherwise.
The vector bundle $\bar V$ is the pullback of $V$ under the projection.
The formula for $\bar h$ is $\bar h(t,x,s):=(h(t,x),s)$ and the
formula for $\bar\delta$ is $\bar\delta(y,s)=\delta(y)$. By
inspection, $(\bar p,\bar g,\bar V,\bar\delta,\bar h)$ is a concordance
from $(p,g,V,\delta,h)\in \cL_{\loc,T}(X)$ to an element 
$(p^{\flat},g^{\flat},V^{\flat},\delta^{\flat},h^{\flat})\in
\cL_{\loc,T}(X)$ where $h^{\flat}$ is a homeomorphism and $g^{\flat}\equiv 0$. 
With some renaming 
we can arrange $h^{\flat}$ to be an identity map, so that   
$(p^{\flat},g^{\flat},V^{\flat},\delta^{\flat},h^{\flat})\in 
\cW_{\loc,T}(X)$. 
If a closed subset $C$ of $X$ is given, and the restriction of 
$(p,g,V,\delta,h)$ to some open neighborhood $U$ of $C$ is already 
in $\cW_{\loc,T}(U)$, then the concordance just constructed is constant 
on $U$, giving the relative surjectivity 
condition in proposition~\ref{prp-relsurjectivity}. \qed 

\begin{cor} 
\label{cor-stratabeauty}
The map $\cL_{\loc,T}\to \cW_{\loc,T}$
of definition~\ref{dfn-LlocTtoWlocT} is a weak 
equivalence. 
\end{cor} 

\proof The composite map, from $\cW_{\loc,T}$ to $\cL_{\loc,T}$ 
and back to $\cW_{\loc,T}$, is clearly a weak equivalence. 
\qed

Summarizing, we have established the weak equivalences of the right
hand column of diagram~(\ref{eqn-stratificdiagram}), and conclude:

\begin{thm} There is a homotopy equivalence
$\displaystyle |\cW_{\loc}| \,\simeq\, \hocolimsub{T\textup{ in } \sK}
\,\, |\cW_{\loc,T}|$. 
\end{thm}

\subsection{Upper left hand column: Couplings} 

\begin{dfn} 
\label{dfn-Wmusheaf}
{\rm An element of $\cW^{\mu}(X)$ is an element $(\pi,f,p,g,V,\dots)$
of $\cW^{\mu}_{\loc}(X)$ such that $(\pi,f)\in \cW(X)$. 
}
\end{dfn}

\begin{dfn}
\label{dfn-coupling}
{\rm A \emph{coupling} between an element $(\pi,f)$
of $\cW(X)$ with $\pi\co E\to X$ and an element 
$(p,g,V)$ of 
$\cL_{\loc}(X)$ with $\omega\co V\to Y$ is a smooth embedding 
$\lambda\co \saddle(V,\rho)\to E$ which satisfies 
$f\lambda(v)=f_V(v)+g(\omega(v))$ for $v\in \saddle(V,\rho)$, 
has $\im(\lambda)\supset \Sigma(\pi,f)$ and respects 
$\Theta$--orientations of the vertical tangent bundles along 
fiberwise singularity sets. 
} 
\end{dfn}

\emph{Remarks, explanations and reminders.}  
The condition $f\lambda(v)=f_V(v)+g(\omega(v))$ implies that
the embedding $\lambda$ takes the zero section of $V$ to the 
fiberwise singularity set $\Sigma(\pi,f)$.
The condition $\im(\lambda)\supset \Sigma(\pi,f)$ forces 
an identification of the vector bundle $\omega\co V\to Y$ 
with $T^{\pi}E|\Sigma(\pi,f)\lra
\Sigma(\pi,f)$. These are the 
vertical tangent bundles along fiberwise singularity sets referred to 
in definition~\ref{dfn-coupling}. Both are $\Theta$--oriented 
vector bundles. 

\medskip
\begin{rmk}
\label{rmk-disturbing}
{\rm The embedding $\lambda\co \saddle(V)\to E$ 
need not have a closed image, because the \'{e}tale map 
$Y\to X$ need not be a closed map. But $\im(\lambda)$ is locally 
compact, therefore locally closed in $E$.
}
\end{rmk}

\begin{dfn}
\label{dfn-Lsheaf} 
{\rm 
For $X$ in $\sX$, an element of $\cL(X)$ is a triple consisting of 
an element in $\cW(X)$, an element in $\cL_{\loc}(X)$ 
and a coupling $\lambda$ between the two. 
}
\end{dfn}

\medskip
\begin{prp}
\label{prp-LWforget}
The forgetful map $\cL\to \cW^{\mu}$ is a weak equivalence. 
\end{prp} 

\proof 
Again we use the relative surjectivity criterion of 
proposition~\ref{prp-relsurjectivity} and again we begin 
with the absolute case. Fix $X$ in $\sX$ and 
$(\pi,f)\in \cW^{\mu}(X)$, with $\pi\co E\to X$.
We want to lift the concordance class of $(\pi,f)$ to a class
in $\cL[X]$. As in the proof of
proposition~\ref{prp-easylocal},
we begin by choosing a vertical tubular neighborhood $V\subset E$ 
of $\Sigma=\Sigma(\pi,f)$, with vector bundle projection
\[\omega\co V\to \Sigma \]
over $X$. Then $V$ is canonically identified with 
$T^{\pi}E|\Sigma$, and so is a Morse vector bundle with 
quadratic function $f_V\co V\to \RR$, corresponding to half 
the Hessian on $T^{\pi}E|\Sigma$.  
By the Morse-Palais lemma~\cite{Lang72}, 
we can set up the vector bundle structure on $V$ in such a way that 
$f(v) = f_V(v)+f\omega(v)$ 
in a neighborhood $U$ of the zero section of $V$.  Without loss of 
generality, the neighborhood $U$  
contains all $v\in \saddle(V,\rho)$ for which 
$|f\omega(v)|\le 1$  and $|f_V(v)|\le 2$. 
(If not, replace $f$ by $\psi f$ where 
$\psi\co E\to[1,\infty[\,$ is a suitable smooth function which factors through 
$\pi\co E\to X$. Multiply the inner product on $V$ by $\psi$, too.
The pairs $(\pi,f)$ and $(\pi,\psi f)$ are 
clearly concordant.)  \newline 
Now choose a smooth embedding $e\co\RR\to\RR$ with 
$\im(e)=\,]-1,1[\,$ and $0<e'\le 1$ throughout. Then
$(\pi,f)$ is concordant to $(\pi^{\sharp},f^{\sharp})$, where 
$\pi^{\sharp}$ is the restriction of $\pi$ to 
$E^{\sharp}= f^{-1}(\im(e))$
and $f^{\sharp}$ is $e^{-1}f$ on $E^{\sharp}$. 
Let $\Sigma^{\sharp}=\Sigma\cap E^{\sharp}$
and $V^{\sharp}=V|\Sigma^{\sharp}$. Let
\[ K= \left\{v\in \saddle(V^{\sharp},\rho)\,\,\big|\,\,\,
|f_V(v)+f\omega(v)|<1\,\right\}. \]
For $v\in K$ we have $|f\omega(v)|<1$ and $|f_V(v)|<2$, so 
$K\subset U$ by our assumptions and consequently 
$f|K=f_V|K + f\omega|K$. It follows that $K\subset E^{\sharp}$.
Using proposition~\ref{prp-longtraceintertwine}, but writing $\lambda$
for $\tau$, we can construct an embedding
\[ \lambda\co \saddle(V^{\sharp},\rho) \lra K \]
relative to and over $\Sigma^{\sharp}$, such that 
$(f_V+f\omega)\circ\lambda=e\circ(f_V+e^{-1}f\omega)=
e\circ(f_V+f^{\sharp}\omega)$. 
This can also be viewed as an 
embedding of $\saddle(V^{\sharp},\rho)$ in $E^{\sharp}$.
We have
\[f^{\sharp}\lambda=e^{-1}f\lambda=e^{-1}(f_V+f\omega)\lambda=
e^{-1}e(f_V+f^{\sharp}\omega)=f_V+f^{\sharp}\omega \]
on $\saddle(V^{\sharp},\rho)$.  
That is, $\lambda$ is a coupling, in the sense of 
definition~\ref{dfn-coupling}, of $(\pi^{\sharp},f^{\sharp})\in \cW(X)$ with 
$(p,g,V^{\sharp})\in \cL_{\loc}(X)$ where $p=\pi^{\sharp}|\Sigma^{\sharp}$
and $g=f^{\sharp}|\Sigma^{\sharp}$. Note that $\lambda$ 
identifies $V^{\sharp}$ with $T^{\pi}E|\Sigma^{\sharp}$, as explained 
in the remarks following definition~\ref{dfn-coupling}, so that 
$V^{\sharp}$ inherits a Morse vector bundle structure and a 
$\Theta$--orientation from $T^{\pi}E|\Sigma^{\sharp}$. 
(The new Morse structure on $V^{\sharp}$ does not quite agree with 
the restriction of the Morse structure on $V$ wich we used 
earlier in this proof. In fact the two structures agree up to a scalar 
factor given by a strictly positive function $\Sigma^{\sharp}\to \RR$.) 
The coupling $\lambda$ promotes the pair 
$(\pi^{\sharp},f^{\sharp})$ to an element of $\cL(X)$ and thereby 
establishes the absolute case of the relative surjectivity 
condition. \newline
The relative case is only slightly more difficult. We sketch it. 
Again fix $X$ in $\sX$ and 
$(\pi,f)$ in $\cW(X)$, with $\pi\co E\to X$. Let $C\subset X$ 
be closed. We want to find an element in $\cL(X)$ 
whose image in $\cW(X)$ is concordant rel $C$ to $(\pi,f)$.
This can be constructed essentially as in the absolute
case, except for one change which consists in replacing the 
embedding $e\co\RR\to\RR$ above by a smooth family of smooth embeddings 
$e_x\co\RR\to\RR$, depending on $x\in X$. Then we have the option to  
choose $e_x=\id_{\RR}$ for $x$ in a small neighborhood of $C$,
while having $\im(e_x)=\,]-1,1[\,$ for $x$ outside a slightly 
larger neighborhood of $C$. \qed 

\medskip The forgetful map $\cL\to \cL_{\loc}$ is not surjective in general, 
nor does it have the concordance lifting property. However
certain ``easy'' concordances in $\cL_{\loc}$ 
can be lifted across the forgetful map $\cL\to \cL_{\loc}$, and this
fact will be needed later. 

\begin{lem} 
\label{lem-easyconcordances} Let $(\pi,f,p,g,V,\lambda)$ be an element of
$\cL(X)$. Let $(\bar p,\bar g,\bar V)\in \cL_{\loc}(X\times\RR)$ 
be a concordance whose initial position is $(p,g,V)\in
\cL_{\loc}(X)$. Suppose that there is a 
pullback diagram 
\begin{equation*}
\xymatrix@C=45pt{
\bar Y \ar[r] \ar[d]^-{\bar p} & Y \ar[d]^-p \\
X\times\RR \ar[r]^-{\textup{proj.}}& X\,. 
} 
\end{equation*}
Then $(\bar p,\bar g,\bar V)$ lifts to a concordance 
$(\bar\pi,\bar f,\dots)\in \cL(X\times\RR)$ whose initial position
is $(\pi,f,p,g,V,\lambda)\in \cL_{\loc}(X)$. If 
the concordance $(\bar p,\bar g,\bar V)$ is relative to a closed 
subset $A$ of $X$, then $(\bar\pi,\bar f,\dots)$ can also be taken 
relative to $A$.
\end{lem}

\proof The statement involves $\RR$ in two ways: as a target 
for functions such as $f$ and $g$, and as a time--like axis 
which parametrizes concordances.   
To reduce confusion, we will write $\RR_{\tau}$ instead of $\RR$ if we mean
the time axis. \newline  
The restriction of $(\pi,f)$ to $\partial(\im(\lambda))$ 
is a submersion $\partial(\im(\lambda))\to X\times\RR$.
This follows from the equation $f\lambda=f_V+g\omega$
(where $\omega\co V\to Y$ is the projection)
and either~(\ref{eqn-longtraceformula}) 
or~(\ref{eqn-duallongtraceformula}). \newline
It is therefore possible to find an outward collar 
for $\im(\lambda)$ in $E$ (the source of $\pi$) 
which is ``over'' $X\times\RR$.
We mean by that a smooth codimension 
zero embedding $u$ of $\partial(\im(\lambda))\times[0,1]$ in 
$E\smin\intr(\im(\lambda))$ which extends the inclusion 
of $\partial(\im(\lambda))\cong \partial(\im(\lambda))\times\{1\}$,
and which satisfies $\pi(u(z,t))=\pi(u(z,1))$ as well as 
$f(u(z,t))=f(u(z,1))$ for all $z\in\partial(\im(\lambda))$ 
and $t\in[0,1]$. Note that $u(\partial(\im(\lambda))\times\{0\})$
is the far end of the collar. \newline
We now construct our 
concordance $(\bar\pi,\bar f,\dots)$ as follows.  
Let $\bar E=E\times\RR_{\tau}$ and let $\bar\pi=\pi\times\id\co \bar E\to 
X\times\RR_{\tau}$. Elements of $\bar E$ should be relabelled to 
ensure that $\bar\pi$ is graphic, but we will not pay much 
attention to that now. Since $\bar Y$ is identified with
$Y\times\RR_{\tau}$, we may also identify $\bar V$ with $V\times\RR_{\tau}$, 
so that $\bar\omega\co\bar V\to \bar Y$ is identified with 
$\omega\times\id\co V\times\RR_{\tau}\to Y\times\RR_{\tau}$. Now we can
define $\bar\lambda$ by $\bar\lambda(v,t)=(\lambda(v),t)\in
E\times\RR_{\tau}=\bar E$. \newline 
It remains to define $\bar f$ on $\bar E$. 
For $z\in E$ outside $\im(\lambda)\cup\im(u)$ and any $t\in
\RR_{\tau}$ we let 
$\bar f(z,t)=f(z)$. For $z=\lambda(v)\in \im(\lambda)$ we must 
define 
\[ \bar f(z,t)=f_V(v)+\bar g(\omega(v),t)
=f(z)+\bar g(\omega(v),t)-g(\omega(v)) \,.\]
This leaves the case $z\in \im(u)$, say $z=u(\lambda(v),s)$ with
$v\in\partial(\saddle(V,\rho))$ and $s\in [0,1]$. In that case we say 
$\bar f(z,t)= f(z)+\bar g(\omega(v),st)-g(\omega(v))$. \qed 

\medskip
\begin{dfn}
\label{dfn-L'Tsheaf}
{\rm For $T$ in $\sK$, we define a sheaf $\cL'_T$
as the pullback of 
\[ 
\xymatrix@=50pt{
\cL \ar[r]^-{\textup{ forget }} &  \cL_{\loc} & 
\ar[l]_-{\textup{ forget }} \cL_{\loc,T}\,. 
}
\]
}
\end{dfn} 

\medskip 
The forgetful maps $\cL'_T\to \cL$ for $T$ in $\sK$ 
determine a canonical map $u$ from the sheaf 
$\hocolim_T\,\cL'_T$ to $\cL$. 

\begin{prp}
\label{prp-firsthocodeco}
The map 
$\displaystyle
u\co \hocolimsub{T\textup{ in }\sK}\,\cL'_T \lra \cL$
is a weak equivalence. 
\end{prp} 

\proof The proof is completely analogous to that of 
proposition~\ref{prp-tamelifting}. There is a factorization of $u$
having the form
\begin{equation*}
\label{eqn-factorofu}
\xymatrix{
{\rule{0mm}{6mm}\hocolimsub{T\textup{ in }\sK} \cL'_{T}} \ar[r]^-{u_1} &
{\beta\cL^{\delta}} \ar[r]^{u_2} & {\cL}
} \ttag
\end{equation*}
where $\cL^{\delta}$ is defined as the pullback of
$\cL\lra \cL_{\loc}\longleftarrow \cL_{\loc}^{\delta}$. One shows 
that $u_1$ and $u_2$ are weak equivalences.  \qed

\subsection{Lower left hand column: Regularization}
\label{subsec-regularization}
In order to make this section more accessible, we 
assume to begin with that $\Theta$ is a discrete space and 
discuss the general case afterwards.

\medskip 
Let $(\pi,f,p,g,V,\delta,h,\lambda)$ be an element of $\cL'_T(X)$ with 
\begin{equation*}
\label{eqn-L'Tdatalist}
\begin{array}{ccl}
(\pi,f)\co E\to X\times\RR\,, & (p,g)\co Y\to X\times\RR\,, 
& V\stackrel{\omega}{\lra} Y\,, \\
\delta\co Y\to \{-1,0,+1\}\,, & 
h\co T\times X\cong \delta^{-1}(0)\,, &
\lambda\co \saddle(V,\rho)\lra E\,. 
\end{array} \ttag
\end{equation*}
We adopt the notation $Y_+=\delta^{-1}(+1)$, $Y_-=\delta^{-1}(-1)$, 
$Y_0=\delta^{-1}(0)$ and let $V_+$, $V_-$, $V_0$ be the restrictions 
of the Morse bundle $V$ to these three (open and closed) 
subspaces of $Y$.

\begin{dfn}
\label{dfn-easyLTsheaf} 
{\rm $\cL_T$ is the subsheaf of $\cL'_T$ 
consisting of the elements $(\pi,f,p,g,V,\delta,h,\lambda)$ as above with 
$g|Y_0\equiv 0$. 
}
\end{dfn}

\begin{prp}
\label{prp-easyprimesnoprimes} 
The inclusion $\cL_T\to \cL'_T$ is a weak equivalence. 
\end{prp}

\proof This 
is a direct application of 
proposition~\ref{prp-relsurjectivity} in 
conjunction with lemma~\ref{lem-easyconcordances}. \qed

\medskip 
For an element $(\pi,f,\dots)$ of $\cL_T(X)$ as above we define 
the regularization $(\pi^{\reg},f^{\reg})$ 
with $\pi^{\reg}\co E^{\reg}\to X$ and $f^{\reg}\co E^{\reg}\to \RR$ by
\begin{equation*}
\label{eqn-regulardata}
\begin{array}{ccl}
E^{\reg}& =& E\smin\lambda(V_+^{\rho}\cup V_0^{\rho}\cup V_-^{-\rho})\,, \\ 
\pi^{\reg} & = & \pi|E^{\reg}\,, \\
f^{\reg}(z)& = &\rule{0mm}{10mm}\left\{ \begin{array}{cc} 
f(z) & \textrm{ if }v\notin \im(\lambda) \\
f^+_V(v) & \textrm{ if }z=\lambda(v) \textrm{ and }
v\in V_+\cup V_0 \\
f^-_{V}(v) & \textrm{ if }z=\lambda(v)\textrm{ and }
v\in V_-\,\,\,. 
\end{array} \right. 
\end{array} \ttag
\end{equation*}  
The maps $f_V^+$ and $f_V^-$ were defined 
in~(\ref{eqn-smile}) and~(\ref{eqn-frown}), respectively. 
We note that $E^{\reg}$ is an open subset of $E$ 
despite remark~\ref{rmk-disturbing}. Moreover, 
$\pi^{\reg}\co E^{\reg}\to \RR$ 
is a smooth submersion and $f^{\reg}$
is regular on each fiber of $\pi^{\reg}$. In fact
\[ (\pi^{\reg},f^{\reg})\co E^{\reg}\lra X\times\RR \]
is a smooth proper submersion. Hence by Ehresmann's fibration 
lemma, we have

\begin{prp}
\label{prp-goodreg}
The map $(\pi^{\reg},f^{\reg})\co E^{\reg}\lra X\times\RR$
is a smooth bundle of closed $d$--manifolds. 
\end{prp} 

It follows that the inverse image of $-1$ under $f^{\reg}$ is a
bundle $q\co M\to X$ of closed $d$--manifolds. Since the restriction 
of $f^{\reg}\circ\lambda$ to $\saddle(V_0,\rho)$ is $f^+_V$ and 
since $f^+_V$ agrees with $f_V$ on levels less than or equal to $-1$, 
the restriction of $\lambda$ gives an embedding
\begin{equation*}
\label{eqn-embeddingfound} 
\begin{array}{ccc}
e\co \{v\in\saddle(V_0,\rho)\mid f_V(v)=-1\} & \lra & M 
\end{array} \ttag
\end{equation*}
The source of $e$ is identified with $D(V_0^{\rho})
\times_{T\times X}S(V_0^{-\rho})$ by
formula~(\ref{eqn-longtraceformula}). The diffeomorphism 
$h\co T\times X\to Y_0$ gives the required element $h^*(V_0)
\in \cW_{\loc,T}(X)$. Since we are (still) assuming that $\Theta$ 
is discrete, the $\Theta$--orientation on $T^{\pi}E^{\reg}$
determines a $\Theta$--orientation on $T^qM$ in a straightforward 
manner. Also, the embedding $e$ respects the
$\Theta$--orientations of the fiberwise tangent bundles. \newline 
Starting from an element in $\cL_T(X)$, we have now produced an
element of $\cW_T(X)$ consisting of $q\co M\to X$ and the 
embedding $e$. 

\bigskip 
It is convenient 
to introduce two subsheaves $\cL^!_T$ and $\cL^{!!}_T$ of $\cL_T$. 
For $\cL^!_T$ we add to the data in~\ref{dfn-easyLTsheaf} 
the condition that $g\ge 1$ on $Y_+$
and $g\le -1$ on  $Y_-$. For the sheaf $\cL^{!!}_T$ 
we add the stronger condition that $\delta\equiv 0$, so that $Y_+$ and
$Y_-$ are empty (and $g\equiv 0$). 

\begin{lem}
\label{lem-movelevels}
The inclusions $\cL^!_T\to \cL_T$ and 
$\cL^{!!}_T\to \cL^!_T$ are weak equivalences. 
\end{lem}

\proof A direct application of~\ref{lem-easyconcordances} shows that 
the inclusion $\cL^!_T\to \cL_T$ is a weak equivalence. For the
inclusion $\cL^{!!}_T\to \cL^!_T$ we use lemma~\ref{lem-basicshrink}.
Given an element $(\pi,f,\dots)$ of $\cL^!_T(X)$ as 
in~(\ref{eqn-L'Tdatalist}), 
we choose a suitable smooth $e\co X\times\RR\to \RR$ such that 
each $e_x\co \RR\to \RR$ defined by $e_x(t)=e(x,t)$ is a smooth 
orientation preserving embedding for all $t\in\RR$, with $e_x(0)=0$. 
In addition we require $0<e'_x\le 1$ for all $x\in X$, with a view to 
proposition~\ref{prp-longtraceintertwine}, and that the image 
of $e_x$ does not contain any nonzero critical values of the 
Morse function $f|E_x$. (For example, if $-1<e_x<1$, then 
$\im(e_x)$ does not contain any nonzero critical values of 
$f|E_x$.) Define $E^{(1)}\subset E$, $\pi^{(1)}$ and 
$f^{(1)}$ exactly as in lemma~\ref{lem-basicshrink}. Let
$V^{(1)}=V_0$. Define $\lambda^{(1)}\co 
\saddle(V^{(1)},\rho)\to E^{(1)}$ by composing 
$\lambda\co \saddle(V,\rho)\to E$ with an embedding 
$\tau\co \saddle(V_0,\rho)\to \saddle(V_0,\rho)$ over $X$ as in 
proposition~\ref{prp-longtraceintertwine}, so that 
$f_V(\tau(v))=e_x(f_V(v))$ for $x\in X$, $y\in Y_0$ with 
$p(y)=x$ and $v\in V_y$. These new 
data, and others obtained by restriction from the old data,  
make up an element $(\pi^{(1)},f^{(1)},\dots)$ of $\cL^{!!}_T(X)$ 
which is concordant in $\cL^!_T$ to $(\pi,f,\dots)$. In fact the proof 
of proposition~\ref{prp-longtraceintertwine} carries over to this 
situation without much change. Therefore
$\cL^{!!}_T[X]\to \cL^!_T[X]$ is surjective. The same argument 
gives surjectivity in the relative case, 
$\cL^{!!}_T[X,A;s]\to \cL^!_T[X,A;s]$, assuming $A\subset X$ 
is closed and $s\in \colim_U\cL^!_T(U)$ where $U$ runs over 
the open neighborhoods of $A$ in $X$. The only detail to watch here 
is that we need $e_x=\id_{\RR}$ for $x$ in a sufficiently 
small neighborhood of $A$. We complete the proof by applying 
proposition~\ref{prp-relsurjectivity}. \qed

\medskip 
\begin{prp}
\label{prp-easyjustreg}
The map $\cL_T\to \cW_T$ defined above (in the case where 
$\Theta$ is discrete) is a weak 
equivalence. 
\end{prp}

\proof By lemma~\ref{lem-movelevels}, it is enough to verify 
that the composition $\cL^{!!}_T\to \cL_T\to \cW_T$ is a weak 
equivalence. But this is almost obvious from
section~\ref{subsec-Morse}. Namely, the long trace construction 
gives us a map of concordance sets $\cW_T[X]\to \cL^{!!}_T[X]$
which is inverse to $\cL^{!!}_T[X]\to \cW_T[X]$. This works 
equally well in a relative setting, so that 
proposition~\ref{prp-relsurjectivity} can be used. 
The only thing to watch here is the $\Theta$--orientation issue. For this, 
fix an element $(q,V,e)$ of $\cW_T(X)$ with $q\co M\to X$ 
and write 
\[ 
\begin{array}{ccc}
e\co \{v\in\saddle(V,\rho)\mid f_V(v) = -1\} & \lra &  M\,. 
\end{array}
\]
Let $E$ be the long trace of $e$, with projection $\pi\co E\to X$. 
Then $E$ contains a copy of $C=M\sqcup_{\im(e)}\saddle(V,\rho)$. A
$\Theta$--orientation of $T^{\pi}E|C$ is already specified. The 
inclusion $C\to E$ is a homotopy equivalence, so that there is no 
obstruction to extending the $\Theta$--orientation of $T^{\pi}E|C$ 
to a $\Theta$--orientation of $T^{\pi}E$. \qed

\begin{rmk} {\rm Here we resolve a set-theoretical issue
related to the question of how exactly $T\mapsto \cW_T$ should 
be regarded as a functor and why the lower square in 
diagram~\ref{eqn-stratificdiagram} is commutative. 
The following item should be added to 
definition~\ref{dfn-WTsheaf}:
{\it a choice of pushout (in the category of sets) 
\[ M\sqcup_{\im(e)}\saddle(V,\rho) \]
with graphic projection to $X$.} \newline
Next we update definition~\ref{dfn-WTsheafinduced2} by 
explaining that $M^{\flat}$ should be regarded as the subspace 
of the pushout 
\[ M\sqcup_{\im(e_a)}\saddle(V_a,\rho_a) \]
consisting of all elements of $M\smin\im(e)$ and those 
$v\in\saddle(V_a,\rho_a)$ with $f_V^{-}(v)=-1$. By the above,
this gives us an explicit choice of underlying set for $M^{\flat}$. 
Note also that 
$M^{\flat}\sqcup_{\im(e^{\flat})}\saddle(V^{\flat},\rho^{\flat})$ is 
identified with a subset of $M\sqcup_{\im(e)}\saddle(V,\rho)$, 
so that we have a preferred choice of underlying set for it.
This choice is to be added to the data $q^{\flat}\co M^{\flat}\to X$ and 
$e^{\flat}$ to give an element of $\cW_S(X)$.
\newline
With these adjustments, $T\mapsto \cW_T$ is a functor. Now
the map $\cL_T\to \cW_T$ of proposition~\ref{prp-easyjustreg} has to 
be adjusted as well, but this is straightforward. Indeed, 
returning to the notation used in~(\ref{eqn-regulardata})  
and~(\ref{eqn-embeddingfound}), we have that $M$ and 
$\lambda(\saddle(V_0,\rho))$ are subsets 
of $E$ with intersection $\im(e)$. Hence their union in $E$ 
is an explicit pushout $M\sqcup_{\im(e)}\saddle(V_0,\rho)$ in 
the category of sets, with a graphic projection map to $X$. By inspection, 
the (updated) maps $\cL_T\to \cW_T$ now constitute a natural transformation 
between two contravariant functors, $T\mapsto \cL_T$ and $T\mapsto
\cW_T$, from $\sK$ to the category of sheaves on $\sX$.
}
\end{rmk}

\bigskip
We end the section with the promised discussion of
general $\Theta$--orientations. We start again with the data 
list~(\ref{eqn-L'Tdatalist}) for an element of $\cL'_T(X)$, 
but do not assume that $\Theta$ is discrete. The coupling $\lambda$ 
identifies $T^{\pi}E|\im(\lambda)$ with $\omega^*V|\saddle(V,\rho)$. 
The differential 
\[  df\co T^{\pi}E \lra f^*(T\RR) \]
is surjective over $E\smin \Sigma\subset E$, where
$\Sigma=\Sigma(\pi,f)$. Over $\im(\lambda)\smin\Sigma$ it has a preferred 
splitting, since $T^{\pi}E|\im(\lambda)$ is a Riemannian vector
bundle. We add the following two items
to definition~\ref{dfn-easyLTsheaf}:
\begin{description}
\item[(A)] A vector bundle splitting of 
$df\co T^{\pi}E|E\smin\Sigma \to f^*(T\RR)|E\smin\Sigma$ which extends 
the preferred splitting over $\im(\lambda)\smin \Sigma$; 
\item[(B)] the condition that $\lambda$ preserve the given 
$\Theta$--orientations of the vertical tangent bundles (and not just 
their restrictions to the fiberwise singularity sets, as in 
definition~\ref{dfn-coupling}).
\end{description}

\begin{prp} 
\label{prp-primesnoprimes} The forgetful map $\cL_T\to \cL'_T$ 
is a weak equivalence. 
\end{prp}

\proof Like the proof of proposition~\ref{prp-easyprimesnoprimes},
this one is a direct application of 
proposition~\ref{prp-relsurjectivity} in 
conjunction with lemma~\ref{lem-easyconcordances}. \qed

\medskip
Given an element of $\cL_T(X)$, consisting of data as 
in~(\ref{eqn-L'Tdatalist}) and items (A) and (B) just above, 
we produce a $d$--manifold bundle $q\co M\to X$ and an embedding 
\[ e\co D(V_0^{\rho})\times_{T\times X}S(V_0^{-\rho})\lra M \]
as before. We note that 
$T^{\pi}E|E^{reg}$ has a canonical splitting, 
\[ 
\begin{array}{ccc}
T^{\pi}E|E^{\reg} & \cong & \ker(df^{\reg})\times \RR \,.
\end{array}
\]
Indeed, over points $z\in E^{\reg}$ not in $\im(\lambda)$ we can use 
the data of item (A) and 
over points  $z\in \im(\lambda)\cap E^{\reg}$ we can use the Riemannian 
vector bundle structure on the fiberwise tangent bundle of 
$\im(\lambda)\to X$. The matching condition in (A) ensures 
that this gives a continuous splitting. 
Since $M\subset E^{\reg}$, we deduce a canonical vector bundle splitting 
\[ 
\begin{array}{ccc} 
T^{\pi}E|M \cong & \cong & T^qM\times \RR\,.
\end{array}
\]
The $\Theta$--orientation on $T^{\pi}E$ therefore induces a 
$\Theta$--orientation on $T^qM\times\RR$, which amounts to 
a $\Theta$--orientation on $T^qM$ itself. \newline
In the same way, the codimension 1 inclusion of   
$\{v\in \saddle(V_0,\rho) \mid f_V(v)=-1\}$ in $\saddle(V,\rho)$
with preferred normal line bundle leads to a $\Theta$--orientation on
the vertical tangent bundle of
\[ \{v\in \saddle(V_0,\rho) \mid f_V(v)=-1\}
\,\,\cong\,\,  D(V_0^{\rho})\times_{T\times X}S(V_0^{-\rho})\,. \]
This is our standard choice of a $\Theta$--orientation on 
the vertical tangent bundle of the source of $e$. With this 
choice $e$ clearly respects the $\Theta$--orientations. Hence 
$(q,V_0,e)$ is a triple satisfying the requirements for an element 
of $\cW_T(X)$ in definition~\ref{dfn-WTsheaf}.

\medskip 
\begin{prp}
\label{prp-justreg}
The map $\cL_T\to \cW_T$ defined above is a weak 
equivalence. 
\end{prp}

\proof The proof of proposition~\ref{prp-easyjustreg} goes through 
without essential changes. \qed

\bigskip 
This completes the construction of
diagram~(\ref{eqn-stratificdiagram}) in the general case.

\subsection{Using the concordance lifting property}
\begin{lem} 
\label{lem-WWlocholift} 
For fixed $T$ in $\sK$, 
the forgetful map $\cW_T\to \cW_{\loc,T}$ has the 
concordance lifting property. 
\end{lem}

\proof We first consider the slighty easier case where $\Theta$ is discrete. 
Fix $X$ in $\sX$.  Any Morse vector bundle $V$ on $T\times X\times\RR$ 
is isomorphic to the pullback of a Morse vector bundle on 
$T\times X$ along the projection $T\times X\times[0,1]\to T\times X$.
It follows easily that any concordance starting at an element  
$z$ of $\cW_{\loc,T}(X)$ is trivial up to an isomorphism of Morse 
vector bundles with $\Theta$--orientation. 
More precisely, suppose that the concordance is given 
by a Morse vector bundle $V$ on $T\times X\times\RR$. For
$s\in\RR$ let $V(s)$ be the pullback (as in~\ref{dfn-graphicmaps}) 
of $V$ along the map $(t,x)\mapsto(t,x,s)$ from $T\times X$ 
to $T\times X\times\RR$. Then we have $V(s)_{(x,t)}=V(1)_{(x,t)}$ for 
$(x,t,s)$ in a neighborhood of $X\times T\times [1,\infty[\,$, and 
similarly $V(s)_{(x,t)}=V(0)_{(x,t)}$ for 
$(x,t,s)$ in a neighborhood of $X\times T\times\,]-\infty,0]$. 
There exists a Morse vector bundle $W$ on $X\times T$ 
and isomorphisms $j_s\co V(s)\to W$ depending smoothly on $s$, 
such that $j_s|V(s)_{(x,t)}=j_1|V(s)_{(x,t)}$ for 
$(x,t,s)$ in a neighborhood of $X\times T\times [1,\infty[\,$, 
and  $j_s|V(s)_{(x,t)}=j_0|V(s)_{(x,t)}$ for 
$(x,t,s)$ in a neighborhood of $X\times T\times\,]-\infty,0]$.
A choice of such a family $(j_s)_{s\in\RR}$ determines, 
for each $y\in\cW_T(X)$ which lifts 
$z$, a lifted concordance starting at $y$. \newline 
In the general case, when $\Theta$ is arbitrary, we begin with the 
construction of a lifted concordance as above, first without worrying
about tangential $\Theta$--orientations. We then have to  
make a choice of $\Theta$--orientation on the fiberwise tangent 
bundle of a manifold bundle of the form 
\[ q\times\RR\co M\times\RR \lra X\times\RR\,. \] 
This is prescribed over the union of $U$ and $\im(e)\times\RR$, 
where $U$ is a neighborhood (germ) of $M\times\,]-\infty,0]$ and 
$e$ is an embedding as in definition~\ref{dfn-WTsheaf}. Since 
the inclusion of 
\[ M\times\,]-\infty,0] \,\,\cup \,\, \im(e) \]
in $M\times\RR$ is a homotopy equivalence, the problem 
can be solved. \qed

\medskip
Now fix an element $(V,\rho)$ in $\cW_{\loc,T}(\pt)$. That is, $V$ is
a $(d+1)$-dimensional $\Theta$--oriented Morse vector bundle on 
$T$. For each $t\in T$, the dimension of ${V_t}^{-\rho}$ is equal to the 
label of $t$ in $\{0,1,\dots,d+1\}$. The following is true by
definition and trivial reformulations.

\begin{lem}
\label{lem-stratfibers}
The fiber of the forgetful map $\cW_T\to \cW_{\loc,T}$
over $(V,\rho)\in \cW_{\loc,T}(\pt)$ is weakly equivalent 
to the sheaf which takes 
an $X$ in $\sX$ 
to the set of all pairs $(q,e)$ where
\begin{description}
\item[(i)] $q\co M\to X$ is a smooth 
graphic bundle of closed $d$--manifolds with a
$\Theta$--orientation of the vertical tangent bundle $T^qM$~;
\item[(ii)] $e\co  D(V^{\rho})\times_T S(V^{-\rho})\times X \,\lra\,
M$ is a smooth embedding over $X$ which is 
fiberwise $\Theta$--orientation preserving. \qed
\end{description}
\end{lem} 

\begin{cor}
\label{cor-stratfibers} 
The fiber of the forgetful map $\cW_T\to \cW_{\loc,T}$
over $(V,\rho)\in \cW_{\loc,T}(\pt )$ is weakly equivalent to 
the sheaf which takes an $X$ in $\sX$ 
to the set of all smooth graphic bundles $q\co M\to X$ of 
tangentially $\Theta$--oriented compact $d$--manifolds with collared 
boundary, where the boundary bundle $\partial M\to X$
is identified with 
\[ 
\begin{array}{ccc}
-(S(V^{\rho})\times_T S(V^{-\rho}))\times X & \lra & X. 
\end{array}
\]
\end{cor}

\proof To get from data $(q,e)$ as in lemma~\ref{lem-stratfibers} to the
kind of bundle described in 
corollary~\ref{cor-stratfibers}, delete the interior 
of $\im(e)$ from the total space of the bundle 
$q$. To get from a bundle $M\to X$ as in 
corollary~\ref{cor-stratfibers} to the data described in 
lemma~\ref{lem-stratfibers}, form the union $M'$ of 
$M$ and $(D(V^{\rho})\times_T S(V^{-\rho}))\times X$
along $(S(V^{\rho})\times_T S(V^{-\rho}))
\times X$. \newline 
The minus sign in front of $S(V^{\rho})\times_T S(V^{-\rho})$ 
indicates a $\Theta$--orientation which is ``opposite'' to the 
one inherited from the tangent bundle of 
$D(V^{\rho})\times_T S(V^{-\rho})$. Equivalently, 
the preferred $\Theta$--orientations on the vertical tangent 
bundles of $M$ and $(D(V^{\rho})\times_T S(V^{-\rho}))\times X$ 
must match to produce a $\Theta$--orientation on the vertical 
tangent bundle of $M'$. 
\qed 

\begin{rmk}
\label{rmk-bifurc} 
{\rm The description of the (homotopy) fiber 
in corollary~\ref{cor-stratfibers} uses only the part of $T$ 
lying over $\{1,2,\dots,d\}\subset \{0,1,2,\dots,d,d+1\}$, 
since spheres of dimension $-1$ 
are empty.} 
\end{rmk}

\subsection{Introducing boundaries} 
\label{subsec-boundaries}
Here we are concerned with a slight generalization of 
diagram~(\ref{eqn-stratificdiagram}). It is obtained by replacing
all families of $(d+1)$--manifolds in sight by families of 
$(d+1)$--manifolds with a prescribed boundary. For the purposes 
of this section we indicate the change by a superscript ``$\partial$'' as 
in $\cdW$; later, in sections~\ref{sec-surgery} and~\ref{sec-stab}, 
the superscript will be dropped.

\medskip
We assume $d>0$ and fix a closed nonempty smooth $(d-1)$--manifold 
$C$ with a $\Theta$--orientation of the tangent bundle $TC$.
The ``prescribed boundary'' which we have in mind will be
$C\times\RR$. We assume also that $C$ is nullbordant in the 
following sense: there exists a compact smooth $d$--manifold $K$
with collared boundary $\partial K=C$ and a $\Theta$--orientation 
of $TK$ which extends the specified one on 
$TC\times\RR\cong TK|C$. (Here we use the outward normal field 
along $C$ to identify $TC\times\RR$ with $TK|C$.)

\begin{dfn}
\label{dfn-deltaWsheaf}
{\rm An element of $\cdW(X)$ is a pair $(\pi,f)$ 
as in~\ref{dfn-Wsheaf}, with $\pi\co E\to X$ and $f\co E\to \RR$, 
except for the following: $\partial E$ need not be empty. We 
ask for a diffeomorphism germ over $X\times\RR$ which identifies
a neighborhood of $\partial E$ in $E$ with a 
neighborhood of $X\times C\times \RR\times\{0\}$ 
in $X\times C\times \RR\times[0,\infty[\,$,
respecting the $\Theta$--orientations.
}
\end{dfn}

\medskip
The same change made in the definitions of the sheaves $\cW_{\loc}$, 
$\cW^{\mu}$, $\cW^{\mu}_{\loc}$,
$\cL$ and $\cL_T$ produces $\cdW_{\loc}$, 
$\cdW^{\mu}$, $\cdW^{\mu}_{\loc}$
$\cdL$ and $\cdL_T$, respectively. (There is also a small change 
in definition~\ref{dfn-coupling}: 
we require $\im(\lambda)\cap \partial E=\emptyset$.) 
No changes are needed in the definitions of 
$\cL_{\loc}$ and $\cL_{\loc,T}$, that is, we put 
$\cdL_{\loc}=\cL_{\loc}$ and $\cdL_{\loc,T}=\cL_{\loc,T}$. 
The case of $\cdW_T$ is slightly different: 

\begin{dfn}
\label{dfn-deltaWTsheaf} 
{\rm An element of $\cdW_T(X)$ is a triple $(q,V,e)$ 
with $q\co M\to X$, as in~\ref{dfn-WTsheaf} except for 
the following. We ask for a diffeomorphism germ 
over $X$ which identifies a neighborhood of $\partial M$ in $M$ with a 
neighborhood of $X\times C\times\{0\}$ in $X\times
C\times[0,\infty[\,$, respecting the $\Theta$--orientations.
We require $\im(e)\cap \partial M=\emptyset$.
}
\end{dfn}

The $\partial$--variant of
diagram~(\ref{eqn-stratificdiagram}) 
is 

\begin{equation*}
\label{eqn-deltastratificdiagram}
\xymatrix{ 
\cdW  \ar[r]  & {\cdW_{\loc}} \\
\cdW^{\mu}  \ar[u]_{\simeq} \ar[r]  & 
{\cdW^{\mu}_{\loc}} \ar[u]_{\simeq} \ar[d]^{\simeq} \\
\cdL \ar[r] \ar[u]_{\simeq} & {\cL_{\loc}}    \\
{\hocolimsub{T\textup{ in }\sK}\,\cdL_T \ar[u]_{\simeq}}
\ar[d]^{\simeq} \ar@<.8ex>[r] &
{\hocolimsub{T\textup{ in }\sK}\,\cL_{\loc,T}} \ar[u]_{\simeq} 
\ar[d]^{\simeq} \\
{\hocolimsub{T\textup{ in }\sK}\,\cdW_T} 
\ar@<.8ex>[r]  & 
{\hocolimsub{T\textup{ in }\sK}\,\cW_{\loc,T}\,.} 
} \ttag
\end{equation*}

To prove that all the maps labelled ``$\simeq$'' are indeed 
weak equivalences, one could proceed roughly as in 
the no--boundary situation. Another method is to show 
that diagrams~(\ref{eqn-deltastratificdiagram}) 
and~(\ref{eqn-stratificdiagram}) can be related by a chain 
of natural transformations which are termwise weak equivalences. 
We now explain how this works for the top left--hand terms, 
and give some indications for the remaining terms. \newline
The first thing we need to know is that $\cdW(\pt)$ is nonempty.
Indeed, an element in $\cdW(\pt)$ is given by $K\times\RR$, 
where $K$ is the nullbordism for $C$ mentioned earlier, 
with the projection map $K\times\RR\to \RR$ and a 
$\Theta$--orientation on $TK\times T\RR$ which can 
be described as the opposite of the one specified 
earlier. (It has to extend the preferred
$\Theta$--orientation on $TC\times\RR\times T\RR$ 
under the identification $TC\times\RR\cong TK|C$ which 
is determined by the \emph{inward} normal field along $C$.) \newline
Modulo set theoretic adjustments as in the proof
of theorem~\ref{thm-spaceWtobordism} and especially the remark
following that proof, the formula 
\[((\pi,f),(\psi,g))\mapsto
(\pi\sqcup\psi,f\sqcup g) \]
defines maps $\cW\times \cW\to \cW$ and $\cW\times \cdW\to
\cdW$. We can view these as an addition law on $\cW$ and an action 
of $\cW$ on $\cdW$, respectively (up to weak equivalences). From 
section~\ref{sec-Vassiliev} we know that $|\cW|\simeq|h\cW|$ and 
from section~\ref{sec-bottomrow} we know that the monoid $\pi_0|h\cW|$ 
is a group. It follows that $\pi_0\cW=\cW[\pt]$ 
is a group under the above addition law. \newline
We now claim that for 
fixed $z\in \cdW(*)$, the 
restriction of the action map to $\cW\times z$ is a weak
equivalence 
\[u\co \cW\times z \to \cdW\,. \]  
Conversely we have a map, essentially from 
$\cdW$ to $\cW$, given by gluing in the ``nullbordism'' $K\times\RR$. 
More precisely, we define 
\[ v\co \cdW\to \cW \]
by taking $(\pi,f)\in \cdW^{\lambda}(X)$ with $\pi\co E\to \RR$ to 
$(\pi\sqcup p_X,f\sqcup p_{\RR})\in \cW(X)$, where $p_X$ and 
$p_{\RR}$ are the projections from $X\times K\times\RR$ to 
$X$ and $\RR$, respectively. The common source of 
$\pi\sqcup p_X$ and $f\sqcup p_{\RR}$ is the pushout 
$E\sqcup_{X\times C\times\RR}(X\times K\times\RR)$. 
(In addition, there are set theoretic precautions to be taken; the
details are left to the reader.)  \newline 
Taking representing spaces, we obtain $|u|$ from 
$|\cW\times z|\simeq|\cW|$ to $|\cdW|$ and $|v|$ from 
$|\cdW|$ to $|\cW|$, both well defined up to homotopy.
It is easy to verify that 
the homotopy classes of $|v||u|$ and $|u||v|$ are both given by
translation with the concordance class $[v(z)]\in \cW[\pt]$. Since 
$\cW[\pt]$ is a group, this shows that $|u|$ and $|v|$ are homotopy
equivalences, i.e. $u$ and $v$ are weak equivalences. 

\medskip
Rows 2 and 3 of diagrams~(\ref{eqn-stratificdiagram}) 
and~(\ref{eqn-deltastratificdiagram}) can be compared in the same 
fashion. For rows 4 and 5 some extra ideas are required. For example, 
to compare $\hocolim_T \cW_T$ and $\hocolim_T\,\cdW_T$ we use 
the sheaf 
\[ \hocolimsub{(S,T)} \cW_S\times\cW_T \]
where $S$ and $T$ are objects of $\sK$ with $S\cap T=\emptyset$. 
The disjoint union maps (with the usual adjustments) and 
substitutions $U=S\sqcup T$ induce 
\[ \hocolimsub{(S,T)} \cW_S\times\cW_T \lra \hocolimsub{U} \cW_U\,. \]
But the map given by specialization to the coordinates, 
\[  \hocolimsub{(S,T)} \cW_S\times\cW_T
\lra  (\hocolimsub{S} \cW_S)\,\times\,(\hocolimsub{T} \cW_T) \]
is a weak equivalence, so that we end up with an addition law 
on $|\hocolim_S\cW_S|$. In the same way, 
we can make an action (up to homotopy) of $|\hocolim_S \cW_S|$
on $|\hocolim_S \,\cdW_S|$. Then a choice of an element 
$z\in \cdW_{\emptyset}(\pt)$ leads, via the action, to a map
\[  \hocolimsub{S} \cW_S\times z \lra \hocolimsub{S} \cdW_S \]
which, by the same reasoning as before, turns out to be a weak 
equivalence. \qed

\medskip 
Corollary~\ref{cor-deltastratfibers} has a variant ``with boundaries''
which looks as follows. 
   
\begin{cor}
\label{cor-deltastratfibers} 
The fiber of the forgetful map $\cdW_T\to \cW_{\loc,T}$
over $(V,\rho)\in \cW_{\loc,T}(\pt )$ is weakly equivalent to 
the sheaf which takes an $X$ in $\sX$ 
to the set of all smooth graphic bundles $q\co M\to X$ of 
tangentially $\Theta$--oriented compact $d$--manifolds with collared 
boundary, where the boundary bundle $\partial M\to X$
is identified with 
\[ 
\begin{array}{ccc}
-(C\,\sqcup\,S(V^{\rho})\times_T S(V^{-\rho}))\times X & \lra & X. 
\end{array}
\]
\end{cor}

\section{The connectivity problem} 
\label{sec-surgery}
\subsection{Overview and definitions} 
Throughout this section we work with the sheaves $\cdW$, $\cdW_S$ 
introduced in section~\ref{subsec-regularization} (which depend 
on the choice of a $(d-1)$--manifold $C$, as specified there). 
But we drop the superscripts and simply write $\cW$, $\cW_S$. We need the 
following extra condition on $\Theta$. (This is satisfied by the 
examples listed in~\ref{exm-Theta}, except for the case
$\Theta=\pi_0\GL\times Y$ when $Y$ is not path--connected.)

\begin{asp} 
\label{asp-transitive action}
{\rm The action of $\pi_0\GL$ on $\pi_0\Theta$ is transitive.}
\end{asp}

\bigskip
The previous section gave us decompositions  of $\cW$ and $\cW_{\loc}$  
into pieces $\cW_S$ and $\cW_{\loc,S}$,
respectively, and a description of the homotopy fibers of 
the forgetful maps 
\[ \cW_S \lra \cW_{\loc,S} \]
as certain bundle theories,
cf.~corollary~\ref{cor-deltastratfibers}. For a given $S$ in $\sK$, 
the $d$--manifolds involved are typically not connected.
In this section we remedy this by showing that upon 
taking the homotopy colimit over $S$, we can in fact assume 
that the relevant $d$--manifolds are connected.

\begin{dfn}
\label{dfn-moreTillmann}
 {\rm For $X$ in $\sX$ let 
$\cW_{c,S}(X)\subset \cW_S(X)$ consist of the triples $(q,V,e)$ 
as in definition~\ref{dfn-WTsheaf}/\ref{dfn-deltaWTsheaf}, 
with $q\co M\to X$ etc., such that the bundle   
$M\smin \im(e) \lra X$ has connected fibers.}
\end{dfn}

Then $\cW_{c,S}$ is a subsheaf of $\cW_S$ and
$|\cW_{c,S}|$ is a union of connected components of 
$|\cW_S|$. The forgetful map from $\cW_{c,S}$ to $\cW_{\loc,S}$ still has 
the concordance lifting property. By analogy with 
corollary~\ref{cor-deltastratfibers}, we have 
the following analysis of its fibers.  

\begin{cor}
\label{cor-connstratfibers}
The fiber of the forgetful map $\cW_{c,S}\to \cW_{\loc,S}$
over $V\in \cW_{\loc,S}(\pt )$ is weakly equivalent to 
the sheaf which takes an $X$ in $\sX$ 
to the set of all smooth graphic bundles $q\co M\to X$ of 
tangentially $\Theta$--oriented smooth 
compact connected $d$--manifolds, where the 
boundary of each fiber $M_x$ is identified with  
\[ -(C\,\sqcup\,S(V^{\rho})\times_S S(V^{-\rho})).  \]
\end{cor}

It would therefore be nice to have a statement saying that 
the inclusion of $\hocolim_S \, \cW_{c,S}$ in $\hocolim_S \, \cW_S$
is a weak equivalence. Unfortunately such a statement is 
nonsensical if we insist on letting $S$ run through 
the entire category $\sK$. We have a contravariant functor
$S\mapsto \cW_S$
from $\sK$ to the category of sheaves on $\sX$, 
but we do \emph{not} have a subfunctor $S\mapsto \cW_{c,S}$. 
It is not the case that the map 
\[ (k,\ep)^*\co \cW_T\to \cW_S \]
induced by a morphism $(k,\ep)\co S\to T$ in $\sK$ will always 
map the subsheaf 
$\cW_{c,T}$ to the subsheaf $\cW_{c,S}$. Let us take a more 
careful look at this phenomenon. \newline 
We may assume that $k$ is an inclusion and that $T\smin S$ 
has exactly one element $t$,  
with label $\lambda(t)\in\{0,1,\dots,d,d+1\}$ and sign 
$\ep(t)\in\{\pm 1\}$. 
Fix $(q,V,e)$ in $\cW_T(X)$, with $q\co M\to X$ and let 
$(q',V',e')$ be the image of 
$(q,V,e)$ in $\cW_S(X)$, 
with $q'\co M'\to X$. For each $x\in X$ there 
is a canonical embedding
\[ M_x\smin \im(e_x) \,\lra \, M'_x\smin \im(e'_x).  \]
The complement of its image is identified with 
\[
\begin{array}{cccc}
 D(V_{(t,x)}^{\rho})\times S(V_{(t,x)}^{-\rho}) & \textup{ if }
&  \ep(t)=+1, & \textup{ and }\\
 S(V_{(t,x)}^{\rho})\times D(V_{(t,x)}^{-\rho}) & \textup{ if }
&  \ep(t)=-1, & 
\end{array}
\]
where $V_{(t,x)}$ is the fiber of $V$ over $(t,x)\in T\times X$.   
We have a problem when the complement is nonempty but 
has empty boundary, because then it will contribute an 
additional connected component. This 
happens precisely when $(\lambda(t),\ep(t))=(d+1,+1)$ and 
when $(\lambda(t),\ep(t))=(0,-1)$. In all other cases, there 
is no problem. 

\medskip
Now our indexing category  
$\sK$ is equivalent to a product $\sK'\times\sK''$. The categories 
$\sK'$ and $\sK''$ can be described as full 
subcategories of $\sK$: namely, $\sK'$ is spanned by 
the objects $S$ whose reference map $S\to \{0,1,2,\dots,d+1\}$ has 
image contained in $\{0,d+1\}$ and $\sK''$ is spanned by 
the objects $S$ whose reference map $S\to \{0,1,2,\dots,d+1\}$ has image 
contained in $\{1,2,\dots,d\}$. \newline
For homotopy colimits of functors from a product 
category to spaces (or to sheaves on $\sX$) 
there is a Fubini principle. In our case it states that 
\begin{equation*}
\hocolimsub{T\textup{ in }\sK} \, \cW_T 
\,\,\simeq\,\, \hocolimsub{Q\textup{ in }\sK'}
\,\hocolimsub{S\textup{ in }\sK''}\, \cW_{Q\amalg S}. \ttag
\end{equation*}

\begin{lem}
\label{lem-slices} 
For any morphism $(k,\ep)\co P\to Q$ in $\sK'$, 
the commutative square 
\[ 
\xymatrix@+5pt{
{\hocolimsub{S\textup{ in }\sK''}\, \cW_{Q\amalg S}}
    \ar@<.8ex>[r]^-{(k,\ep)^*}  \ar[d]  &
{\hocolimsub{S\textup{ in }\sK''}\, \cW_{P\amalg S}} \ar[d]  \\
{\hocolimsub{S\textup{ in }\sK''}\, \cW_{\loc, Q\amalg S}}
 \ar@<.8ex>[r]^-{(k,\ep)^*} &  
{\hocolimsub{S\textup{ in }\sK''}\, \cW_{\loc, P\amalg S}}
}
\]
is homotopy cartesian (after passage to representing spaces). 
\end{lem}

\begin{thm}
\label{thm-surgery} The inclusion 
\[
\hocolimsub{S\textup{ in }\sK''}\,\cW_{c,\,S}
\,\lra\,
\hocolimsub{S\textup{ in }\sK''}\,\cW_S    
\]
is a weak equivalence. 
\end{thm}

\bigskip
Theorem~\ref{thm-surgery} is the main result of 
the section.  We develop a surgery method to prove it. The idea 
is to make nonconnected $d$--manifolds connected by means 
of multiple surgeries on embedded (thickened) $0$-spheres. 
Then we need to know that such multiple 0-surgeries on 
a $d$--manifolds are essentially unique. In order to state the 
uniqueness, we view them as the objects of a category.

\subsection{Categories of multiple surgeries }
\label{subsec-multicat}
In this section we fix a compact, smooth, nonempty 
$d$--manifold $M$ with a $\Theta$--orientation of $TM$. 
Unless otherwise stated, $\RR^{d+1}$ will be regarded as a Morse 
vector space with the standard inner product and involution 
$\rho(x_1,\dots,x_d,x_{d+1})=(x_1,\dots,x_d,-x_{d+1})$. We shorten 
$D((\RR^{d+1})^{\rho})\times S((\RR^{d+1})^{-\rho})$
to $D^d\times S^0$. This is normally identified with the space 
$\{v\in \saddle(\RR^{d+1},\rho) \,|\,\lag v,\rho(v)\rag=-1\}$. Hence any 
$\Theta$--orientation on the tangent bundle of $\RR^{d+1}$ will 
induce one on the tangent bundle of $D^d\times S^0$.  Cf. the 
discussion leading up to proposition~\ref{prp-justreg}. 

We make a slight change in the definitions  
of $\cW_S$ and $\cW_{\loc,S}$. Namely, where 
definitions~\ref{dfn-WlocTsheaf} and~\ref{dfn-WTsheaf} ask for a Morse 
vector bundle $\omega\co V\to Y$ with a $\Theta$--orientation 
on $V$ itself, we will be satisfied with a $\Theta$--orientation on 
the Morse vector bundle $T^{\omega}V\cong \omega^*V$. This change 
does not affect the homotopy types of $|\cW_S|$ and 
$|\cW_{\loc,S}|$. 

\begin{dfn} 
\label{dfn-surgecat} {\rm 
Let $\sC_M$ be the topological category defined as follows. An object 
consists of a finite set $T$, a $\Theta$--orientation on the tangent 
bundle of $\RR^{d+1}\times T$ and a smooth embedding 
$e\pd_T$ of $D^d\times S^0\times T$ in $M\smin\partial M$ which
respects the $\Theta$--orientations and satisfies the following
condition: Surgery on $e\pd_T$ results in a 
\emph{connected} $d$--manifold. 
A morphism from $(S,e\pd_S)$ to $(T,e\pd_T)$ is an injective map 
$k\co S\to T$ such that $k^*e\pd_T=e\pd_S$. \newline 
The set of objects of $\sC_M$ is topologized as a subspace 
of the disjoint union, over all $T$, of the spaces 
\[
\begin{array}{cl}
& (\textup{space of smooth embeddings }
D^d\times S^0\times T \lra M\smin\partial M) \\
\times & (\textup{space of $\Theta$--orientations on the tangent 
bundle of $\RR^d\times T$}).
\end{array}
\]
The total morphism set $\mor(\sC_M)$ is topologized as 
a closed subset of $\ob(\sC_M)\times\ob(\sC_M)$ via the map 
(source,target).
}
\end{dfn}

\begin{prp} 
\label{prp-surgery} The space $B\sC_M$ is contractible. 
\end{prp}

The proof requires a lemma.

\begin{lem}
\label{lem-etalereduction}
Let $\sigma\co N\to X$ be a  
submersion of smooth manifolds without boundary, $\dim(N)>\dim(X)$. 
Suppose that for each $x\in X$ there 
exists a contractible open neighborhood $W$ of $x$ in $X$, a finite set 
$Q$ and a map $Q\times W\to N$ over $X$ inducing a surjection 
from $Q\cong\pi_0(Q\times W)$ to $\pi_0(N_y)$ for every $y\in W$. 
Then there exists a locally finite covering of $X$ by 
contractible open sets $W_j$,  where $j\in J$, and
finite sets $Q_j$, and a smooth 
embedding 
\[ a \co \coprod_j Q_j\times W_j \lra N \]
over $X$, such that the restriction of $a$ to 
$Q_j\times W_j$ induces 
surjections $Q_j\to \pi_0(N_x)$, for each $j\in J$ 
and $x\in W_j$. 
\end{lem} 

\begin{exm} {\rm The submersion $\RR^2\smin(0,0) \lra \RR~;~(x,y)\mapsto x$
satisfies the hypothesis of lemma~\ref{lem-etalereduction}. 
The submersion $\RR\smin 0\to \RR~;~x\mapsto x$ 
does not. Surjectivity is not directly related to the issue; 
the projection from $(\RR\times\{0,1\})\smin(0,0)$ to $\RR$ 
is a surjective submersion which also fails to satisfy the 
hypothesis of lemma~\ref{lem-etalereduction}. 
}
\end{exm}

\proof[Proof of lemma~\ref{lem-etalereduction}.] Note first that 
the statement is not completely trivial. Using 
the hypothesis, we could start with a locally finite covering  
of $X$ by contractible open sets $W_j$, and choose finite sets 
$Q_j$ and maps $a_j\co Q_j\times W_j\to N$ over $X$ inducing surjections 
$Q_j\to \pi_0(N_y)$ for every $y\in W_j$. This would give us a map 
\[ a\co \coprod_j Q_j\times W_j \lra N \]
which is an immersion. Unfortunately there is no guarantee that 
it is an embedding. To solve this problem we will partition a 
``large'', dense open subset $U$ of $N$ into ``levels'' indexed by 
the real numbers, and arrange that $a$ maps distinct 
connected components of $\coprod Q_j\times W_j$ to distinct levels
of $U$. Then $a$ is an embedding. \newline
The jet transversality theorem, 
applied to sections of the vertical tangent bundle of $N$,
implies that we can find a $k\gg 0$
and a smooth $f\co N\to \RR$ such that the 
fiberwise $k$-jet prolongation $j^k_{\sigma}f\co N\to J^k_{\sigma}(N,\RR)$ 
is nowhere 0. 
Let $U\subset N$ consist of all $z\in N$ such that $f|N_{\sigma(z)}$ 
is regular at $z$. Then $U$ is open in $N$ 
and $U_x:= U\cap N_x$ is dense in $N_x$, for each $x\in X$.
Hence the inclusions $U_x\to N_x$ induce surjections 
$\pi_0(U_x)\to \pi_0(N_x)$. 
The hypotheses on $\sigma$ now give us a  
a covering of $X$ by contractible 
open subsets $W_j$, and for each $W_j$ 
a finite set $Q_j$ and a map $a_j\co Q_j\times W_j\to U$ over 
$X$ such that the induced composite 
map $Q_j\to \pi_0(U_x)\to \pi_0(N_x)$ is onto 
for every $x\in W_j$. We can assume that the 
$W_j$ are the open stars of the vertices in a sufficiently fine 
triangulation of $X$, in which case the covering is locally
finite. But in addition we can easily arrange 
that $fa_j$ is constant on $q\times W_j$ for each $q\in Q_j$, 
and that the resulting map $\coprod_j\,Q_j\to \RR$ is injective. 
Then the map $a$ which equals $a_j$ on $Q_j\times W_j$ satisfies 
all our requirements. \qed

\bigskip
In the proof of theorem~\ref{thm-surgery}, we will use a sheaf 
version $\cC_M$ of $\sC_M$. For connected 
$X$ in $\sX$ let $\cC_M(X)$ be the (discrete) category
whose objects are the pairs $(T,e\pd_T)$ 
where $T$ is a finite set together with a $\Theta$--orientation on the 
tangent bundle of $\RR^{d+1}\times T$, and 
\[ 
\begin{array}{ccc}
e\pd_T\co D^d\times S^0\times T \times X & \lra & 
(M\smin\partial M)\times X 
\end{array}
\] 
is a smooth embedding over $X$, respecting the tangential 
$\Theta$--orientations and subject to the condition 
that fiberwise surgery on $e\pd_T$ results in a bundle of 
connected manifolds. A morphism from $(S,e_S)$ 
to $(T,e_T)$ is an injective map $k\co S\to T$ such that
$k^*e_T=e_S$. \newline
Since $\ob(\cC_M(\Delta^k_e))$ is the set of smooth maps from 
$\Delta^k_e$ to the space $\ob(\sC_M)$, one gets a 
functor of topological categories $|\cC_M|\to \sC_M$ 
which induces a degreewise homotopy equivalence of 
the nerves and therefore a homotopy equivalence 
$B|\cC_M\op|\cong B|\cC_M|\to B\sC_M$. (Here it is best to 
define $B\sC_M$ as the fat realization \cite{Segal74} of the nerve 
of $\sC_M$, ignoring the degeneracy operators.) 

\proof[Proof of proposition~\ref{prp-surgery}.] We show that 
$\beta\cC_M\op$ is weakly equivalent to the terminal sheaf taking every
$X$ in $\sX$ to a singleton. (Without loss of generality, we can 
assume that the fixed indexing set $J$ which is used in the 
$\beta$--construction is uncountable.) By
proposition~\ref{prp-relsurjectivity}, this reduces to the 
following
\begin{itemize}
\item[]
\emph{Claim.} Let $X$ in $\sX$ be given with a closed subset $A$
and a germ $s\in \colim_U\,\beta\cC_M\op(U)$, where $U$ ranges over the 
neighborhoods of $A$ in $X$. Then $s$ extends to an element of 
$\beta\cC_M\op(X)$.
\end{itemize}
To verify this, choose an open neighborhood $U$ of $A$ in $X$ such
that the germ $s$ can be represented by some $s_0\in \beta\cC_M\op(U)$. 
The information contained in $s_0$ includes a locally finite 
covering of $U$ by open subsets 
$U_j$ for $j\in J$. (Making $U$ smaller if necessary, we can assume
that this is locally finite in the strong sense that every $x\in X$ 
has a neighborhood which meets only finitely many $U_j$.)
It also includes a choice of 
object $\psi_{RR}\in \ob(\cC_M(U_R))$ for each finite nonempty 
subset $R$ of $J$. (There are also 
morphisms $\psi_{RS}\in \mor(\cC_M(U_S))$, but they are of course determined 
by their sources $\psi_{RR}|U_S$ and targets $\psi_{SS}$.) 
Next, choose an open $X_0\subset X$ 
such that $U\cup X_0=X$ and the closure of $X_0$ in $X$ avoids $A$. \newline
Let $N$ be the open subset of $(M\smin\partial M)\times X_0$ obtained 
by removing from $(M\smin\partial M)\times X_0$ the closures of the 
embedded disk bundles determined by the various 
$\varphi_{RR}|U_R\cap X_0$. By making $U$ and $X_0$ and the $U_j$ 
smaller if necessary, but taking care that the $U_j$ remain the 
same near $A$, we can arrange that the projection $N\to X_0$ 
satisfies the hypothesis of lemma~\ref{lem-etalereduction}.
\newline 
By the lemma, there exists a locally finite covering of $X_0$ by 
contractible open sets $U'_j$, and finite sets $Q_j$ and an embedding
$a$ of $\coprod_j Q_j\times U'_j$ in $N$, over $X_0$ , such that 
$a$ induces surjections $Q_j\to \pi_0(N_x)$ 
for each $j$ and $x\in U'_j$. (Again, making $X_0$ smaller if
necessary, we can assume
that this is locally finite in the strong sense that every $x\in X$ 
has a neighborhood which meets only finitely many $U'_j$.)
We can also choose a smooth embedding $b$ of  
$\coprod_j Q_j\times U'_j$ in $N$,
over $X_0$, inducing constant maps $Q_j\to \pi_0(N_x)$ for each $j$ 
and all $x\in U'_j$, and such that $\im(a)\cap \im(b)=\emptyset$. 
(For example, the distinct sheets of $b$ restricted to $Q_j\times U'_j$  
can be chosen very close to a selected sheet of $a$.) 
Since the $U'_j$ are contractible, the normal bundles of $a$ and 
$b$ can be trivialized (as $d$-dimensional vector bundles),
and so the ``union'' of $a$ and $b$ extends to a smooth and 
fiberwise $\Theta$--orientation preserving embedding 
\[ 
\begin{array}{ccc}
c\co D^d\times S^0\times \coprod_j (Q_j\times U'_j)
& \lra  & N
\end{array}
\]
over $X_0$, with suitably chosen $\Theta$--orientations on 
the vertical tangent bundles of the projections 
$\RR^{d+1}\times Q_j\times U'_j\lra Q_j\times U'_j$.
(Such trivializations and such $\Theta$--orientations can be 
chosen thanks to our assumption that $\pi_0\GL$ acts transitively 
on $\pi_0\Theta$.) For each $j$ with nonempty $U'_j$, the restriction 
of $c$ to the summand 
\[ D^d\times S^0\times Q_j\times U'_j \]
is an object $\varphi_{jj}$ of $\cC_M(U'_j)$. 
Because $J$ is uncountable, we can arrange that $U'_j$ is empty whenever $U_j$ 
is nonempty.
\newline 
We are now ready to define an explicit element in $\beta\cC_M\op(X)$
which extends the germ $s$. Let $Y_j=U_j$ if $U_j$ is nonempty, 
$Y_j=U'_j$ if $U'_j$ is nonempty, 
and $Y_j=\emptyset$ for all other $j\in J$. Then the $Y_j$ form a 
locally finite  
open covering of $X$. For finite $R\subset J$ 
with nonempty $Y_R$, we can write $Y_R=U_S\cap U'_T$ for 
disjoint subsets $S,T$ of $R$ with $S\cup T=R$. Let 
$\varphi_{RR}\in \ob(\cC_M(Y_R))$ be the coproduct (which exists 
by construction) of $\psi_{SS}|Y_R$ and the $\varphi_{jj}|Y_R$ for 
$j\in T$. The covering $j\mapsto Y_j$ together with the data 
$\varphi_{RR}$ for finite nonempty $R\subset J$ is an element in 
$\beta\cC_M\op(X)$ which extends the germ $s$. \qed 

\medskip
The category $\sC_M$ in definition~\ref{dfn-surgecat}
is not quite ideal for our purposes. Suppose that $\partial M$ 
is identified with $-C$. Then 
$M$ determines a vertex $v$ in $|\cW_{\emptyset}|$.
Each object $(T,e_T)$ in $\sC_M$ determines a path in 
\[ \hocolimsub{S\textup{ in }\sK''}\,|\cW_S| \]
starting at $v$ and ending somewhere in the subspace $|\cW_{c,\emptyset}|$.
(The path is composed of two edges. The first edge leads from $v$ 
to the vertex in $|\cW_T|$ represented by $M$ and the embedding $e_T$.
The second edge leads from there to the vertex in $\cW_{c,\emptyset}$ 
represented by $M(e_T)$, the manifold obtained from $M$ by surgery 
on $e_T$. Note that we view $T$ as a set over $\{1,2,\dots,d\}$ by mapping 
each $t\in T$ to $1\in\{1,2,\dots,d\}$. We are making use of two morphisms 
$\emptyset \to T$ in $\sK''$, one having $\ep(t)=+1$ for all $t\in T$
and the other having $\ep(t)=-1$ for all $t\in T$.)
One might hope that morphisms in $\sC_M$ determine in the same manner
homotopies between paths in 
\[ \hocolimsub{S\textup{ in }\sK''}\,|\cW_S| \]
starting at $v$ and ending somewhere in the 
subspace $|\cW_{c,\emptyset}|$. That is not the case, but the 
problem is easy to fix. We use Segal's edgewise 
subdivision of of $\sC_M$.

\medskip
\begin{dfn}
\label{dfn-edgewise}
 {\rm For any category $\sD$, the edgewise subdivision 
$\es(\sD)$ of $\sD$ is another category defined as follows. 
An object of $\es(\sD)$ is a morphism $f\co c_0\to c_1$ in $\sD$. 
A morphism in $\es(\sD)$ from an object $f\co c_0\to d_0$
to an object $g\co d_0\to d_1$ is a commutative square in $\sD$ of the form
\[
\xymatrix{ 
c_0  \ar[r]^f & c_1 \ar[d] \\
d_0 \ar[u] \ar[r]^g & d_1\,. 
}
\]
}
\end{dfn}

\medskip
It is well known that $B(\es(\sD))$ is homeomorphic 
to $B\sD$, if $\sD$ is a discrete category. More precisely, 
by \cite[Lm.2.4]{GoodwillieKleinWeiss02} 
the nerve of $\es(\sD)$ is isomorphic as a simplicial set 
to the edgewise subdivision of the nerve of $\sD$, 
and this implies by \cite{Segal73} that the realizations are 
homeomorphic. In the case of a simplicial category $\sD$ 
one can argue degreewise. The general case of a topological 
category can in most cases be reduced to the case of a 
simplicial category. In particular:

\begin{cor}
\label{cor-surgery}
The classifying space of $\es(\sC_M)$ is 
contractible. 
\end{cor} 

With the notation introduced just before definition~\ref{dfn-edgewise},
each object $(T,U,e_T)$ of $\es(\sC_M)$ determines a path in 
$\hocolim_S\,|\cW_S|$ starting at $v\in |\cW_{\emptyset}|$ and ending 
somewhere in the subspace 
$\hocolim_S\,|\cW_{c,S}|$, where $S$ runs through $\sK''$.
(The path is again composed of two edges. The first edge leads from $v$ 
to the vertex in $|\cW_T|$ represented by $M$ and the embedding $e_T$.
The second edge leads from there to the vertex in $\cW_{c,T\smin U}$ 
represented by $M(e_U)$, the manifold obtained from $M$ by surgery 
on $e_U$, together with the embedding $e_{T\smin U}$ obtained 
by restricting $e_T\,$.) 
It turns out that morphisms in $\es(\sC_M)$ do indeed 
determine homotopies between such paths starting at $v$ 
and ending somewhere in the subspace 
$\hocolim_S\,|\cW_{c,S}|$. For now we leave the verification 
to the reader; in a moment it will be made explicit 
in a more general setting.

\medskip

\subsection{Parametrized multiple surgeries} 
We reformulate corollary~\ref{cor-surgery}  
in a parametrized setting and deduce
theorem~\ref{thm-surgery} from the reformulation. 

\begin{dfn} {\rm Fix an object $S$ in $\sK''$. 
Let $(T,U)$ be a pair of finite sets with $U\subset T$ and 
$T\cap S=\emptyset$. We view $T$ as an object of $\sK''$
by mapping each $t\in T$ to the element $1$ of $\{0,1,2,\dots,d+1\}$.
We introduce a sheaf $\cW_{S;T,U}$
on $\sX$ with a forgetful map to $\cW_S$. \newline
For $X$ in $\sX$, an element in $\cW_{S;T,U}(X)$ is an element 
$(q,V,e)$ of $\cW_{S\sqcup T}(X)$ with $q\co M\to X$ etc.,  
where the restriction of $V$ to $T\times X$ is identified 
with a trivial (Morse) vector bundle $\RR^{d+1}\times T\times X$. 
\emph{Condition:} Fiberwise surgery on 
\[ e\pd_U\co D^d\times S^0\times U\times X \lra 
M\smin \im(e_S) \]
results in a bundle of connected $d$--manifolds; here $e\pd_U$ and 
$e\pd_S$ denote the restrictions of $e$ to the portions 
of the source lying over $U\times X$ and $S\times X$, respectively. 
}
\end{dfn} 

Let $\sP$ be the category whose objects are pairs of finite 
sets $(T,U)$ with $U\subset T$, where a morphism from $(Q,R)$ 
to $(T,U)$ is an injective map $h\co Q\to T$ with $U\subset h(R)$.
(This is equivalent to the edgewise subdivision of the category of
finite sets and injective maps.)
A morphism $(Q,R)\to (T,U)$ in $\sP$ induces a map of sheaves
$\cW_{S;T,U}  \lra \cW_{S;Q,R}$, so that 
there is a contravariant functor from $\sP$ to sheaves on $\sX$ given by 
\[  (T,U)\,\mapsto \,\cW_{S;T,U} . \]

\medskip
\begin{cor}
\label{cor-parametricsurgery} The forgetful maps $\cW_{S;T,U}\to \cW_S$
induce a homotopy equivalence 
\[ \hocolimsub{(T,U)}\,|\cW_{S;T,U}|\quad \simeq\quad 
|\cW_S|\,. \]
\end{cor}

\proof Fix an element $(q,V,e)$ in $\cW_S(\pt)$, with $q\co M\to X$.
It is enough to show that the homotopy fiber of 
\[ \hocolimsub{(T,U)}\,|\cW_{S;T,U}| \lra |\cW_S| \]
over the point corresponding to $(q,V,e)$ is contractible. In this situation 
the processes of forming homotopy colimits and taking homotopy fibers 
commute. Moreover each of the 
forgetful maps $\cW_{S;T,U}\to \cW_S$ has the concordance lifting
property, so by proposition~\ref{prp-sheavespullback},
the homotopy fiber which we are
interested in is weakly equivalent to 
\begin{equation*}
\label{eqn-hocolim1}
\hocolimsub{(T,U)} |\,\fiber_M(\cW_{S;T,U}\to \cW_S)|. \ttag
\end{equation*}
Let $M_S$ be the compact surface obtained from $M$ 
by deleting $\intr(\im(e))$. It is clear that each 
expression $|\,\fiber_M(\cW_{S;T,U}\to \cW_S)|$ in~(\ref{eqn-hocolim1}) 
can be replaced by the naturally homotopy equivalent 
\[ \ob_{(T,U)}\es(\sC_{M_S}), \]
the space of objects in $\es(\sC_{M_S})$ of definitions~\ref{dfn-surgecat}
and~\ref{dfn-edgewise} whose underlying injection of finite sets 
is the inclusion $U\to T$. The 
homotopy colimit now becomes the classifying space of the 
transport category 
\[ \sP\op\ssmallint \ob_{\bullet}\es(\sC_{M_S})\,, \]
cf.\ section~\ref{subsec-hocolim}, 
where the bullet stands for objects $(T,U)$ of $\sP$. 
This transport category is clearly equivalent to 
$\es(\sC_{M_S})$, and its classifying space is therefore 
contractible by corollary~\ref{cor-surgery}. \qed

\proof[Proof of theorem~\ref{thm-surgery}.]
Using the homotopy invariance 
property of homotopy direct limits, we obtain from 
corollary~\ref{cor-parametricsurgery} a 
homotopy equivalence of spaces
\[ 
\xymatrix@C+12pt@M+4pt{
{\eta_+ :\hocolimsub{S\textup{ in }\sK''}
 \hocolimsub{(T,U)\textup{ in }\sP}\, |\cW_{S;T,U}|}
\ar@<.8ex>[r] & {\hocolimsub{S\textup{ in }\sK''}\, |\cW_S|\,.} 
}
\]
We compare this with the map 
\begin{equation*}
\label{eqn-dothesurgery1}
\xymatrix@C+12pt@M+4pt{
{\eta_-\co \hocolimsub{S\textup{ in }\sK''}
 \hocolimsub{(T,U)\textup{ in }\sP}\, |\cW_{S;T,U}|}
\ar@<.8ex>[r] &  {\hocolimsub{R\textup{ in }\sK''}\, |\cW_{c,R}|} 
} \ttag
\end{equation*}
induced by the composite maps
\begin{equation*}
\label{eqn-dothesurgery2}
\xymatrix@C+8pt@M+10pt{
\cW_{S;T,U} \ar[r] &  \cW_{S\cup T} \ar[r]^{(-)^*} &
\rule{0mm}{3.6mm}\cW_{S\cup(T\smin U)}
} \ttag
\end{equation*}
and renaming $S\sqcup(T\smin U)$ as $R$. 
Here the first arrow in~(\ref{eqn-dothesurgery2})
is self-explanatory. The second 
is induced by the inclusion $S\sqcup(T\smin U)\to S\sqcup T$, with the 
sign function on $U$ which is $\equiv-1$.
Thus the first arrow amounts to adding the surgery data 
corresponding to labels in $T$ (but not performing any 
surgeries), while the second amounts  
to performing the surgeries corresponding to labels in 
$U\subset T$. It follows that the composite map
in~(\ref{eqn-dothesurgery2}) 
lands in the subsheaf $\cW_{c,S\sqcup(T\smin U)}$, as required
in~(\ref{eqn-dothesurgery1}). 
The map $\eta_-$ in~(\ref{eqn-dothesurgery1})
is clearly a retraction, with a 
canonical section $\zeta$ which identifies each $\cW_{c,R}$ 
with $\cW_{R;\emptyset,\emptyset}$. 
The target of $\eta_-$ is contained in the target of $\eta_+$~, 
so we may ask whether $\eta_-$ and $\eta_+$ are homotopic 
as maps to $\hocolim_S\, |\cW_S|$.  
This is indeed the case, by remark~\ref{rmk-naturalhocolim2} and
the fact that each $\cW_{S;T,U}$ fits into a natural commutative diagram 
\[
\xymatrix@+10pt{
& \cW_{S;T,U} \ar[dl]_{\textup{forget}} \ar[d] 
\ar[dr]^{(\ref{eqn-dothesurgery2})} & \\
\cW_S & \ar[l]_{(+)^*}  \cW_{S\sqcup T} \ar[r]^{(-)^*} &
\cW_{S\sqcup(T\smin U)}\,.
}
\]
The homotopy restricts to a 
constant homotopy from $\eta_+\zeta$ to $\eta_-\zeta$. Consequently, 
it is a deformation retraction of $\hocolim_S\, |\cW_S|$  
to $\hocolim_S\, |\cW_{c,S}|$. \qed

\subsection{Annihiliation of $d$-spheres} 
The goal is to prove lemma~\ref{lem-slices}. Most of the proof 
is based on some elementary product decompositions. 

\begin{lem}
\label{lem-locsplitoff} Let $T=T_1\cup T_2$ be a disjoint union, where
$T_1$ is an object of $\sK'$ and $T_2$ is an object of $\sK$. 
There are weak equivalences, natural in $T_2$ for fixed $T_1$~, 
\[
\begin{array}{cccccc}
\cW_T  & \lra & \cW_{\loc, T_1}\times \cW_{T_2}\,, \qquad &
\cW_{\loc,T} & \lra & \cW_{\loc, T_1}\times \cW_{\loc,T_2}. 
\end{array}
\]
\end{lem}

\proof The second map is induced by the inclusions $T_1\to T$ 
and $T_2\to T$. It should be clear that it is a weak 
equivalence. Note that sign functions on $T_2$ and $T_1$ 
are not needed. \newline 
The first coordinate of the 
first map is again induced by the inclusion $T_1\to T$. 
The second coordinate of the first map, 
\[ \cW_T \lra \cW_{T_2}~, \]
is defined as follows. Let $(q,V,e)$ be an element of $\cW_T(X)$ 
as in definition~\ref{dfn-WTsheaf}, with $q\co M\to X$. For $a\in T_1$~,
the bundle   
\[ D(V^{\rho}_a)\times_{X_a} S(V^{-\rho}_a) \]
(where $X_a= a\times X$ and $V_a=V|X_a$) is either empty or 
a bundle of $d$-spheres. In any case it has empty boundary 
and its image under $e$ is a union of connected components of $M$. 
Let $M'$ be obtained from $M$ 
by deleting these components, for all $a\in T_1$. 
Let $V'$ be the restriction of $V$ to 
$T_2\times X$ and let $e'$ be the restriction of $e$ to
\[ \coprod_{b\in T_2} D(V^{\rho}_b)\times_{X_b} S(V^{-\rho}_b)\,. \]
Then $(q',V',e')\in
\cW_{T_2}(X)$. This determines the map $\cW_T \lra \cW_{T_2}$. 
Again it should be clear that the resulting map 
\[ \cW_T \lra \cW_{\loc,T_1}\times \cW_{T_2} \]
is a weak equivalence: it is easy to write down an inverse 
for the induced map on homotopy groups. \qed

\medskip
\proof[Proof of lemma~\ref{lem-slices}.] Applying 
lemma~\ref{lem-locsplitoff} with and noting that 
homotopy colimits commute with the functors $\cW_{\loc,Q}\times\,$
and $\cW_{\loc,P}\times\,$, and passing to representing spaces
we can rewrite 
the commutative diagram in lemma~\ref{lem-slices} in the form
\[ 
\xymatrix@+5pt{
{|\cW_{\loc,Q}|\,\times|\!\hocolimsub{S\textup{ in }\sK''}\, 
\cW_S\,|}
    \ar@<.8ex>[r]^-{(k,\ep)^*}  \ar[d]^{\id\times\ell}  &
{|\cW_{\loc,P}|\,\times|\!\hocolimsub{S\textup{ in }\sK''}\, 
\cW_S\,|} \ar[d]^{\id\times \ell}  \\
{|\cW_{\loc,Q}|\,\times|\!\hocolimsub{S\textup{ in }\sK''}\, 
\cW_{\loc,S}\,|}
 \ar@<.8ex>[r]^-{(k,\ep)^*} &  
{|\cW_{\loc,P}|\,\times
|\!\hocolimsub{S\textup{ in }\sK''}\, \cW_{\loc,S}\,|}
}
\]
where $\ell\co \hocolim_S\,\cW_S\lra \hocolim_S\,\cW_{\loc,S}$ is the 
forgetful map. For $y\in \cW_{\loc,Q}(\pt)$ 
and $z\in \hocolim_S\,\cW_S(\pt)$, the homotopy fiber of 
the left--hand vertical arrow over $(y,z)$ is therefore identified 
with $\hofiber_z(\ell)$ 
and the homotopy fiber of the right--hand vertical arrow 
over the image point $((k,\ep)^*y,z)$ is also identified 
with $\hofiber_z(\ell)$. However, with these identifications 
the map
\[ u\co \hofiber_z(\ell) \lra \hofiber_z(\ell) \]
induced by the horizontal arrows in the diagram is not 
obviously the identity. To understand what it is, 
we can assume that $Q\smin P$ has exactly one element $a$. 
Associated with this we have a label $n_a\in\{0,d+1\}$ and a 
value $\ep(a)\in \{-1,+1\}$. By inspection, if 
$(n_a,\ep(a))$ is $(0,+1)$ or $(d+1,-1)$, then the map $u$ 
is the identity. To describe what happens in the remaining cases, 
we note that by choosing $y$ we have also selected an element 
$(p,W,g)\in \cW_{\loc,\{a\}}(\pt)$ where $W$ is a vector bundle 
with inner product over a point, i.e., a vector space with 
inner product. We may identify $W$ with $\RR^{d+1}$. 
Now the map $u$ is given by disjoint union with 
$S(W)=S^d$, assuming that $(n_a,\ep(a))$ is $(d+1,+1)$ or
$(0,-1)$. More precisely, for each $S$ in 
$\sK''$ and $X$ in $\sX$, we have a map 
\[ \bar u\co \cW_S(X) \to \cW_S(X) \]
given by $(q,V,e)\mapsto (q^{\sharp},V,e)$ where 
$q\co M\to X$ is a bundle of $d$--manifolds etc., and $q^{\sharp}$
is obtained from $q$ by disjoint union with a trivial sphere 
bundle $S^d\times X\to X$. This is natural 
in the variables $X$ and $S$. It covers the identity 
map of $\cW_{\loc,S}(X)$ and so induces $u$ above. Hence 
it only remains to show that $\bar u$ is a weak equivalence.

\begin{lem} 
\label{lem-spherekilling} The map 
\[ \bar u\co \hocolimsub{S\textup{ in }\sK''}\, \cW_S
\, \lra\,\hocolimsub{S\textup{ in }\sK''}\, \cW_S\]
given by disjoint union of all $d$--manifolds in
sight with $S^d$ is a weak equivalence. 
\end{lem} 

\proof We reason as in section~\ref{subsec-boundaries}. 
This will require two variants of $\cW_S$ as in
definition~\ref{dfn-WTsheaf}, one where we use 
$-C$ as the prescribed boundary and another where we use 
$C\sqcup -C$, in other words, the boundary of $C\times[0,1]$.
To distinguish these, we write $\cdW_S$ for the first and $\cddW_S$ for the
second. \newline
Concatenation defines maps $\cddW_S\bar\times\, \cddW_T\to 
\cddW_{S\sqcup T}$ and $\cddW_S\bar\times \,\cdW_T\to 
\cdW_{S\sqcup T}$ where $\bar\times$ denotes a mildly enhanced version 
of the product. Hence 
\[ \hocolimsub{S\textup{ in }\sK''} |\,\cddW_S| \]
becomes a homotopy associative $H$--space. It has a homotopy 
unit given by the element $C\times [0,1]$ in
$\cddW_{\emptyset}(\pt)$. This $H$--space acts 
(in a homotopy associative manner) on 
\[ \hocolimsub{S\textup{ in }\sK''} |\,\cdW_S|. \]
Up to homotopy, the map $\bar u$ is given by translation with a 
single element $z$ of
\[   |\,\cddW_{\emptyset}| \subset \hocolimsub{S\textup{ in }\sK''}
  |\,\cddW_S|\,.\]
Namely, $z$ is the element defined by $S^d\sqcup (C\times[0,1])$. It is
therefore enough to show that $z$ is in the connected component of 
the homotopy unit, defined by $C\times[0,1]$. This amounts to saying 
that $S^d\sqcup (C\times[0,1])$ can be transformed into $C\times[0,1]$
by elementary surgeries of index $1,2,\dots,d$ only. In fact 
a single surgery of index $1$, that is, a surgery on a thickened 
$0$--sphere in $S^d\sqcup (C\times[0,1])$, is enough. (Let 
one component of the thickened $0$--sphere be in $S^d$ 
and the other in $C\times[0,1]$. Here at last we are using the 
assumption that $C\ne\emptyset$.) \qed 

\medskip
\section{Stabilization}
\label{sec-stab}

\subsection{Stabilizing the decomposition}

\emph{Conventions.} Throughout this section we assume $d=2$ and 
$\Theta=\pi_0\GL$ with the translation action of $\GL$. (This means 
that $\Theta$--orientations are ordinary orientations.) We 
continue to write $\cW$ and $\cW_S$ for $\cdW$ and $\cdW_S$,
respectively. The fixed boundary $C$ is 
$S^1\sqcup -S^1$.

\medskip
In section~\ref{sec-surgery}, we modified
the homotopy colimit decomposition of $|\cW|$  
obtained in section~\ref{sec-strat}. The goal was
to banish non-connected $d$--manifolds from the 
picture, as far as possible. Here 
we make a second and quite drastic modification to 
our homotopy colimit decomposition (assuming $d=2$ etc.); the goal 
is roughly to ensure that all surfaces in sight are of large genus, 
in addition to being connected. We achieve this by 
repeatedly concatenating with a standard surface of 
genus 1, with boundary $S^1\sqcup -S^1$. This standard 
surface can be viewed as an element $z\in \cW_{c,\emptyset}(\pt)$. 

\medskip 
For every $X$ in
$\sX$, the unique map $X\to \pt$ induces $\cW_{c,\emptyset}(\pt)\to 
\cW_{c,\emptyset}(X)$ 
and so allows us to think of $z$ as an element of 
$\cW_{c,\emptyset}(X)$. For $S$ in $\sK$ 
define $z^{-1}\cW_S$ and $z^{-1}\cW_{c,S}$
as the colimits, in the category of sheaves on $\sX$, of the diagrams
\[ \begin{array}{c}
\cW_S \stackrel{z\cdot}{\lra}\cW_S
\stackrel{z\cdot}{\lra}\cW_S
\stackrel{z\cdot}{\lra}\cW_S \stackrel{z\cdot}{\lra}\cdots, \\
\cW_{c,S} \stackrel{z\cdot}{\lra}\cW_{c,S}
\stackrel{z\cdot}{\lra}\cW_{c,S}
\stackrel{z\cdot}{\lra}\cW_{c,S} \stackrel{z\cdot}{\lra}\cdots, 
\end{array}
\]
respectively. The arrows labelled $z\cdot$ are given by concatenation with 
$z$. These colimits are obtained by sheafifying the 
naive colimits, which are presheaves. The sheafification process is 
very mild in this case. In particular, it does not alter the values 
on compact objects of $\sX$, such as spheres. Hence the representing
spaces of these colimits are homotopy equivalent to the colimits 
of the individual representing spaces: 
\[  |z^{-1}\cW_S| \simeq z^{-1}|\cW_S|\,, 
\qquad\quad |z^{-1}\cW_{c,S}| \simeq z^{-1}|\cW_{c,S}|\,. \]

\medskip
For an object $T$ in $\sK''$, corollary~\ref{cor-connstratfibers} 
implies that the homotopy fiber of the localization map 
$ |\cW_{c,T}| \lra |\cW_{\loc,T}| $
over any base point is homotopy equivalent to 
$\coprod_g\, B\Gamma_{g,\,2+2|T|}$. 

\medskip
\begin{lem} For $T$ in $\sK''$, any homotopy fiber of 
$|z^{-1}\cW_{c,T}| \lra |\cW_{\loc,T}|$
is homotopy equivalent to $\ZZ\times B\Gamma_{\infty,\,2+2|T|}$\,. \qed
\end{lem} 

\medskip 
Finally we have the stabilized versions of 
lemma~\ref{lem-slices} and theorem~\ref{thm-surgery}:

\medskip 
\begin{cor}
\label{cor-stabslices} 
For any morphism $(k,\ep)\co P\to Q$ in $\sK'$, 
the commutative square 
\[ 
\xymatrix@C=50pt@M+4pt@C+4pt{ 
{\hocolimsub{S\textup{ in }\sK''}\, |z^{-1}\cW_{Q\amalg S}| }
\ar@<.8ex>[r]^-{(k,\ep)^*}  \ar[d] & 
{\hocolimsub{S\textup{ in }\sK''}\, |z^{-1}\cW_{P\amalg S}| } \ar[d] \\
{\hocolimsub{S\textup{ in }\sK''}\, |\cW_{\loc, Q\amalg S}| }
\ar@<.8ex>[r]^-{(k,\ep)^*} & 
{\hocolimsub{S\textup{ in }\sK''}\, |\cW_{\loc, P\amalg S}| }
}
\]
is homotopy cartesian. 
\end{cor}

\begin{cor}
\label{cor-stabsurgery} The inclusion 
\[ \hocolimsub{T\textup{ in }\sK''}\, |z^{-1}\cW_{c,\,T}|
\,\lra\,
\hocolimsub{T\textup{ in }\sK''}\, |z^{-1}\cW_T| \] 
is a homotopy equivalence.  
\end{cor}

\medskip
Corollaries~\ref{cor-stabslices} and~\ref{cor-stabsurgery} 
are about a new homotopy colimit decomposition of $|\cW|$: 

\begin{lem} 
\label{lem-stabledecomposition}
$\displaystyle |\cW|\,\,\simeq\,|\,z^{-1}\cW| 
\,\,\simeq\,\,\hocolimsub{T\textup{ in }\sK}\, |z^{-1}\cW_T|$\,.
\end{lem}  

\proof Since $|\cW|$ is group complete, the inclusion 
$|\cW|\to z^{-1}|\cW|\simeq |z^{-1}\cW|$ is a homotopy 
equivalence. The second homotopy equivalence in the chain 
follows from $|z^{-1}\cW_T|\simeq z^{-1}|\cW_T|$ and 
\[
\hocolimsub{T\textup{ in }\sK}\, z^{-1}|\cW_T|
\,\,\simeq \,\, z^{-1}\big(\!\hocolimsub{T\textup{ in
  }\sK}\,|\cW_T|\,\big). \qquad\qed
\]

\subsection{Using the Harer-Ivanov stability theorem} 

\begin{lem}
\label{lem-Hareruse} The canonical map from 
$\ZZ\times B\Gamma_{\infty,2}$ to the homotopy fiber 
(over the base point) of the forgetful map 
\[ \hocolimsub{S\textup{ in }\sK''}\, |z^{-1}\cW_{c,S}|
\quad\lra\quad \hocolimsub{S\textup{ in }\sK''}\, |\cW_{\loc,S}|
\]
induces an isomorphism in homology with integer coefficients. 
\end{lem}

\proof For the object $S=\emptyset$ of $\sK''$, 
we have $|z^{-1}\cW_{c,S}|\simeq 
\ZZ\times B\Gamma_{\infty,2}$ and $|\cW_{\loc,S}|= \pt$, 
so that there is indeed a canonical map from 
$\ZZ\times B\Gamma_{\infty,2}$ to the homotopy fiber of 
\[ \hocolimsub{S\textup{ in }\sK''}\, |z^{-1}\cW_{c,S}|
\quad\lra\quad \hocolimsub{S\textup{ in }\sK''}\, |\cW_{\loc,S}|. \]
We now check that the hypothesis of 
corollary~\ref{cor-hocolimhomologyfibration} is satisfied.  
Let $(k,\ep)\co S\to T$ be a morphism in $\sK''$. We have to verify 
that, in the commutative square of spaces 
\[ 
\xymatrix@M+4pt@W+4pt{
|{z^{-1}\cW_{c,T}|} \ar[r] \ar[d]^-{(k,\ep)^*}   & |\cW_{\loc,T}| 
\ar[d]^-{(k,\ep)^*} \\
|{z^{-1}\cW_{c,S}|} \ar[r] & |\cW_{\loc,S}|,
}
\]
the induced map from any of the homotopy fibers in the upper row 
to the corresponding homotopy fiber in the lower row induces an 
isomorphism in homology. The homotopy fibers in question are related 
by a map 
\[  \ZZ\times B\Gamma_{\infty,\,2+2|T|}
\lra \ZZ\times B\Gamma_{\infty,\,2+2|S|}\, \]
given geometrically by attaching cylinders $D^1\times S^1$ or 
double disks $D^2\times S^0$ to those pairs of boundary 
circles which correspond to elements of $T\smin k(S)$. This 
map is an integral homology equivalence by the Harer-Ivanov stability 
theorem. Apply corollary~\ref{cor-hocolimhomologyfibration}. \qed

\begin{cor}
\label{cor-almostthere}
The canonical map from 
$\ZZ\times B\Gamma_{\infty,2}$ to the homotopy fiber 
(over the base point) of the forgetful map 
\[ \hocolimsub{S\textup{ in }\sK}\,|z^{-1}\cW_S|
\quad\lra\quad \hocolimsub{S\textup{ in }\sK}\, |\cW_{\loc,S}|
\]
induces an isomorphism in homology with integer coefficients. 
\end{cor}

\proof Combine lemma~\ref{lem-Hareruse} with 
corollaries~\ref{cor-stabsurgery} and~\ref{cor-stabslices}. \qed

\proof[Proof of theorem~\ref{thm-toprow}.]
By lemma~\ref{lem-stabledecomposition} and 
diagram~\ref{eqn-stratificdiagram}, we have 
\[  \hocolimsub{S\textup{ in }\sK}\,|z^{-1}\cW_S|\,\,\simeq\,\
|\cW|\,,
\qquad\quad \hocolimsub{S\textup{ in }\sK}\, |\cW_{\loc,S}|
\,\,\simeq\,\ |\cW_{\loc}|\,. \]
Therefore corollary~\ref{cor-almostthere} implies that 
the homotopy fiber of $|\cW|\to |\cW_{\loc}|$ receives a 
map from $\ZZ\times B\Gamma_{\infty,2}$ which induces an 
isomorphism in integer homology. But 
$|\cW|$ and $|\cW_{\loc}|$ are infinite loop spaces by 
theorems~\ref{thm-middlecolumn} and section~\ref{sec-bottomrow}, 
and the map $|\cW|\to |\cW_{\loc}|$ is an infinite loop map. 
Hence its homotopy fiber is an infinite loop space, and each 
of its components has an abelian fundamental group. Each of 
these fundamental groups is isomorphic to
$H_1(B\Gamma_{\infty,2\,};\ZZ)=0$. Summing up, all connected components
of the homotopy fiber in question are simply connected, and the 
homotopy fiber is therefore $\ZZ\times B\Gamma_{\infty,2}^+$. \qed

\setcounter{section}{0}
\renewcommand{\thesection}{\Alph{section}}

\section{More about sheaves}
\label{sec-spacesboxdetails}

\subsection{Concordance and the representing space}
\label{subsec-concordanceandrepresent}

Let $\cF$ be a sheaf on $\sX$.
We shall construct a natural transformation 
$\vartheta\co [X,|\cF|\,] \lra \cF[X]$,
and an inverse 
$\xi\co \cF[X]\to [X,|\cF|\,]$ for $\vartheta$. 

We start with the construction of $\xi$.
Fix $X$ in $\sX$ and an element $u\in \cF(X)$. Choose a smooth 
triangulation of $X$, with vertex set $T$. 
Suppose that $S\subset T$ is a distinguished subset 
(the vertex set of a simplex 
in the triangulation). Let 
\[ 
\begin{array}{rcl}
\Delta_e(S)  & = & \{ w\co S\to \RR\, \mid \,\Sigma_s w(s) = 1\}  \\
\Delta(S) & = & \{w\in \Delta_e(S)\, \mid\, w\ge 0\}.
\end{array}
\]
The triangulation gives us characteristic embeddings
$c_S\co \Delta(S)\to X$, one for each distinguished $S\subset T$. 
By induction on $S$, we can choose 
smooth embeddings 
\[c_{e,S}\co \Delta_e(S)\to X~, \] 
extending the $c_S$~, which are compatible: i.e., if $S$ is distinguished and
$R\subset S$ is nonempty, then $c_{e,S}$ agrees with $c_{e,R}$ 
on $\Delta_e(R)\subset \Delta_e(S)$. Let $u_S={c_{e,S}}^*(u)\in 
\cF(\Delta_e(S))$. 
\newline
Finally choose a total ordering of $T$. This leads to  
an identification of each $\Delta_e(S)$
with a standard extended simplex. Consequently it
promotes each $u_S$ to a simplex of the simplicial set 
$\uli n\mapsto \cF(\Delta_e^n)$. We then have a unique map 
$\xi(u)\co X\to |\cF|$ such that, for each $S$ as above with 
$|S|=n+1$, the  
diagram 
\[ 
\xymatrix{
\Delta(S) \ar[r]^{\cong}  \ar[d]^{c_S} &  
\Delta^n \ar[d]^{\textup{char. map for }u} \\
X \ar[r]^{\xi(u)} & |\cF| 
}
\]
commutes. It is straightforward to show that the resulting 
homotopy class of maps $X\to |\cF|$ depends only on the concordance 
class of $u\in \cF(X)$.

We remark that $\xi\co \cF[X]\to [X,|\cF|\,]$
so defined is a natural transformation. Indeed 
if $f\co X\to Y$ is a smooth embedding, then $f^*\xi=\xi f^*$
by inspection. An arbitrary morphism $g\co X\to Y$ in $\sX$ can be 
factored as $pf$, where $f\co X\to Y\times\RR^k$ is a smooth 
embedding for some $k$ and $p\co Y\times\RR^k\to Y$ is the 
projection. Let $s\co Y\to Y\times\RR^k$ be any smooth section of $p$. 
Then $s^*\xi=\xi s^*$, and consequently $p^*\xi=\xi p^*$ since 
$p$ is inverse to $s$ in the homotopy category of $\sX$. Therefore 
$g^*\xi=f^*p^*\xi=f^*\xi p^*=\xi f^*p^*=\xi g^*$. 

\medskip
The construction of an inverse $\vartheta$ for $\xi$ uses a 
simplicial approximation principle. To introduce notation 
for that, we suppose first that $L$ is a simplicial complex with a 
totally ordered vertex set $T$. For $n\ge 0$ let $L^s_n$ be the set 
of order preserving maps $f\co \{0,1,\dots,n\}\to T$ such that
$\im(f)$ is a simplex of $L$. Then $n\mapsto L^s_n$ is a simplicial 
set $L^s$ and the realization $|L^s|$ is homeomorphic to $L$. \newline 
Next, let $K$ be any simplicial complex and let $Q$ be a simplicial 
set. Our approximation principle states that, for any homotopy class of
maps from $K$ to $|Q|$, there exist a subdivision $L$ of $K$, with a 
total ordering of its vertex set, and a simplicial map $L^s\to Q$ 
such that the induced map from $|L^s|\cong L\cong K$ to $|Q|$ is in 
the prescribed homotopy class. (Probably the easiest proof proceeds 
by induction over the skeleta of $K$, with the original simplicial 
complex structure. We leave it to the reader.)

\bigskip
Next we construct $\vartheta\co [X,|\cF|\,] \lra \cF[X]$. 
We start with a choice of map $g\co X\to |\cF|$. By the above 
approximation principle, we may assume that $X$ comes with a 
smooth triangulation, with totally ordered vertex set $T$, 
and that $g$ is the realization of a simplicial map from $X^s$ 
to the simplicial set $n\mapsto \cF(\Delta^n_e)$.
In particular, each distinguished subset $S\subset T$ with $|S|-1=n$
determines a nondegenerate $n$--simplex $y_S$ of $X^s$ and then 
an element $g(y_S)\in \cF(\Delta^n_e)$. Now 
choose a smooth homotopy of smooth maps $h_t\co X\to X$, where 
$0\le t\le 1$, such that $h_0=\id$ and
\begin{enumerate}
\item for every $t$, the map
$h_t$ maps each simplex of the triangulation to itself;
\item each simplex 
has a neighbourhood in $X$ which is mapped to the simplex 
by $h_1$. 
\end{enumerate}   
Then for each $n\ge 0$ and each distinguished
subset $R\subset T$ with $|R|-1=n$ and a sufficiently small open 
neighborhood 
$V_R$ of $c_R(\Delta(R))$ in $X$, we obtain a smooth map 
$V_R\to \Delta_e(R)\cong\Delta_e^n$ by composing $h_1|V_R$ 
with the inclusion of $\Delta(R)$ in $\Delta_e(R)$. 
Using this map to pull back $g(y_R)\in \cF(\Delta_e^n)$, we obtain 
compatible elements $z_R\in \cF(V_R)$ which, by the sheaf 
property of $\cF$, determine a unique element $\vartheta(g)$ of 
$\cF(X)$. Again, it is straightforward to verify that the 
concordance class of $\vartheta(g)$ depends only on the 
homotopy class of $g$. 

\begin{prp}
\label{prp-xiandtheta}
The maps $\xi$ and $\vartheta$ are inverses of each other.
\end{prp}    

\proof Let $u\in \cF(X)$. We want to show that 
$\vartheta\xi(u)$ is concordant to $u$. With suitable choices in the 
constructions above, we have $V_R\supset \im(c_{e,R})$ for all 
distinguished $R$, and then 
$\vartheta\xi(u)$ equals ${h_1}^*(u)$, 
where $(h_t\co X\to X)_{0\le t\le 1}$ is the homotopy which appears in the 
definition of $\vartheta$. Since $h_1$ is smoothly homotopic to 
$h_0=\id_X$, this implies that $\vartheta\xi(u)$ is 
indeed concordant to $u$. Therefore 
\[ 
\vartheta\xi=\id\co \cF[X] \lra \cF[X] . \]
To show that $\xi\vartheta$ is the identity on $[X,|\cF|\,]$, we
return to the notation used and assumptions made in the 
construction of $\vartheta$. In particular we have $g\co X\to |\cF|$, 
induced by a simplicial map from $X^s$ to the simplicial set
$n\mapsto \cF(\Delta^n_e)$, and a homotopy of smooth maps 
$h_t$ from $X$ to $X$. In addition we have, from the construction of 
$\xi$, smooth embeddings $c_{e,S}\co \Delta_e(S)\to X$
extending the characteristic embeddings $c_S\co \Delta(S)\to X$
of the simplices of $X$. With some care, we can arrange that each 
$h_t\co X\to X$ takes the image of each $c_{e,S}$ to itself. We can 
also arrange that $h_t=h_0$ for $t$ close to $0$ and $h_t=h_1$ 
for $t$ close to $1$, and define $H\co X\times\RR\to X$ 
by 
\[
H(x,t)=\left\{\begin{array}{cl} 
h_t(x) & \qquad t\in [0,1] \\
h_1(x) & \qquad t\ge 1 \\
h_0(x) & \qquad t\le 0\,. 
\end{array}\right.
\]
We introduce the notation $\cF^{\RR}$ for the sheaf 
$Y\mapsto \cF(Y\times\RR)$ on $\sX$, and note that the embeddings 
$y\mapsto (y,0)$ and $y\mapsto (y,1)$ of $Y$ in $Y\times\RR$ 
determine maps of sheaves $\textrm{ev}_0~,\textrm{ev}_1\co \cF^{\RR}\to \cF$. 
The point of this is that our data so far determine a simplicial map 
$G$ from $X^s$ to $n\mapsto\cF^{\RR}(\Delta^n_e)$, and consequently 
a map $|G|\co X\to |\cF^{\RR}|$. Namely, for a nondegenerate
$n$--simplex $y_S$ of $X^s$ let $G(y_S)\in \cF^{\RR}(\Delta^n_e)$  
be the pullback of $g(y_S)\in \cF(\Delta^n_e)$ along      
\[ (c_{e,S})^{-1}\circ H\circ (c_{e,S}\times\id_{\RR})\co 
\Delta^n_e\times\RR\lra \Delta^n_e\,, \]
where we identify $\Delta^n_e$ with $\Delta_e(S)$.
Lemma~\ref{lem-evhomotopy} below implies that $g=|\textrm{ev}_0G|$ and
$\xi\vartheta(g)=|\textrm{ev}_1G|$ are homotopic. \qed

\begin{lem}
\label{lem-evhomotopy} The evaluation maps 
$|\textrm{ev}_0|,\,|\textrm{ev}_1|\co |\cF^{\RR}|\to|\cF|$ are homotopic. 
\end{lem}

\proof For an order preserving map $f\co\uli n\to \uli 1$ let 
$\bar f\co \Delta^n_e\to \Delta^n_e\times\RR$ be the unique affine embedding
which takes a vertex $v$ of $\Delta^n$ to $(v,f(v))$. The formula 
$(u,f)\mapsto \bar f^*(u)$ determines a simplicial homotopy, 
i.e., a simplicial map from 
$n\mapsto 
\cF(\Delta^n_e\times\RR)\times\mor_{\Delta}(\uli n,\uli 1)$
to $n\mapsto \cF(\Delta^n)$. The homotopy connects 
$\textrm{ev}_0$ with $\textrm{ev}_1$. \qed

\proof[Proof of proposition~\ref{prp-whatitrepresents}.]
The special case where the closed subset $A$ is empty 
is covered by proposition~\ref{prp-xiandtheta}.
The proof of the general case follows the same lines. 
To construct $\xi[u]$ for $u\in \cF(X,A;z)$, we choose a smooth 
triangulation of $X$ where each simplex which meets $A$ is contained 
in a fixed open neighborhood $Y$ of $A$ with $u|Y=z$. Conversely, 
for a relative homotopy class of maps $X\to |\cF|$ taking $A$ to $z$, 
we can find a smooth triangulation of $X$ with totally ordered 
vertex set and a simplicial map from $X^s$ to $n\mapsto \cF(\Delta^n_e)$
taking every nondegenerate simplex of $X^s$ which meets $A$ to $z$, 
and representing the relative homotopy class.   \qed

\subsection{Categorical properties}

\begin{prp}
\label{prp-productsheaf} The construction $\cF\mapsto |\cF|$ 
takes pullback squares of sheaves to 
pullback squares of compactly generated Hausdorff spaces. 
In particular it respects products.
\end{prp} 

\proof The functor $\cF\to |\cF|$ is a composition of two 
functors: one from sheaves to simplicial sets, and another 
from simplicial sets to compactly generated Hausdorff spaces. 
It is obvious that the first of these respects pullbacks. The 
second also respects pullbacks by \cite[\S3,
  Thm. 3.1]{GabrielZisman67}. 
\qed

\begin{dfn}
\label{dfn-coproductsheaf} 
{\rm The categorical coproduct $\cF_1\amalg \cF_2$ of 
two sheaves $\cF_1$ and $\cF_2$ on $\sX$ can be defined by 
$(\cF_1\amalg \cF_2)(X)=\prod_i \cF_1(X_i)\amalg
\cF_2(X_i)$
where $X_i$ denotes the path component of $X$ corresponding 
to an $i\in\pi_0(X)$.} 
\end{dfn}

\begin{prp} 
\label{prp-coproductsheaf}
$|\cF_1\amalg \cF_2|\cong |\cF_1|\amalg |\cF_2|$. 
\end{prp} 

\proof Note that $\Delta^n_e$ 
is path-connected for $n\ge 0$. \qed

\begin{prp}
\label{prp-sheavespullback} Suppose given sheaves $\cE, \cF, \cG$ 
on $\sX$ and morphisms (alias natural transformations) 
$u\co\cE\to\cG$, $v\co\cF\to \cG$. Let $\cE\times_{\cG}\cF$ 
be the fiber product (pullback) of $u$ and $v$. If $u$ 
has the concordance lifting property, definition~\ref{dfn-lifting},
then the projection 
$\cE\times_{\cG}\cF\to \cF$ has the concordance lifting property and the
following square is homotopy cartesian: 
\[
\xymatrix@W+1pt{
|\cE\times_{\cG}\cF| \ar[r] \ar[d] & |\cF| \ar[d]^v \\
|\cE|\ar[r]^u & |\cG|. 
}
\]
\end{prp}

We begin with a special case of 
proposition~\ref{prp-sheavespullback}. 
Given a natural transformation $u\co \cE\to\cG$
of sheaves on $\sX$ with the concordance lifting property, 
let $z$ be a point in $\cG(\pt )$ and let $\cE_z$ be the fiber 
of $u$ over $z$ (in the category of sheaves). Let 
$\hofiber_z\,|u|$ denote
the homotopy fiber of $|u|\co |\cE|\to |\cG|$ over 
the point $z$. 

\begin{lem}
\label{lem-relativesheafhomotopy}
 For any $y\in \cE_z(\pt )$, the homotopy set 
$\pi_n(\cE_z,y)$ is in canonical bijection with 
$\,\pi_n(\hofiber_z\,|u|, y)\,$.
\end{lem} 

\proof Because of the concordance lifting 
property, $\pi_n(\cE_z,y)$ can be identified with an 
appropriate relative homotopy group (or homotopy set) 
of the map of sheaves 
$u\co\cE\to\cG$. Representatives of the latter 
are elements $(s,h)\in \cE(S^n)\times \cG(S^n\times\RR)$, 
where $s$ has the value $y$ near the base point $S^n$ and
$h$ is a concordance (relative to a neighborhood of the base point)
from $u(s)$ to the constant $z$. We can 
identify this relative homotopy group 
(set) with a relative homotopy group (set) of 
the map of spaces $|u|\co |\cE|\to |\cG|$, which can then be 
identified with a homotopy group (set) of 
the homotopy fiber of $|u|$ over $z$. \qed

\begin{cor} 
\label{cor-morerelativesheafhomotopy}
In the situation of lemma~\ref{lem-relativesheafhomotopy},
the sequence  
\[ 
\xymatrix@M+2pt{
|\cE_z| \incto[r] &  |\cE| \ar[r]^{|u|} & |\cG| 
} 
\]
is a homotopy fiber sequence. 
\end{cor}

\proof The composite map from $|\cE_z|$ to $|\cG|$ is constant. 
This leads to a canonical map from $|\cE_z|$ to the homotopy 
fiber of $|u|\co|\cE|\to |\cG|$ over $z$. It is easy to verify 
directly that this induces a surjection on $\pi_0$. 
For each $y\in \cE_z(\pt )$, the induced map of 
homotopy sets  
\[ \pi_n(\cE_z,y) \lra \pi_n(\hofiber_z\,|u|, y) \]
is the one from lemma~\ref{lem-relativesheafhomotopy}. 
It is therefore always a bijection. \qed
  
\proof[Proof of proposition~\ref{prp-sheavespullback}.] We fix $z\in\cF(\pt )$ 
and obtain $v(z)\in \cG(\pt )$. The fiber of 
\[ \cE\times_{\cG}\cF \lra  \cF \]
over $z$ is identified with the fiber of $u\co\cE\to\cG$ 
over $v(z)$. Using corollary~\ref{cor-morerelativesheafhomotopy}
we can conclude that the homotopy fiber of $|\cE\times_{\cG}\cF| \lra |\cF|$ 
over $z$ maps to the homotopy fiber of $|u|\co|\cE|\to|\cG|$ 
over $v(z)$ by a homotopy equivalence. \qed

\subsection{Cocycle sheaves and classifying spaces}

This section contains the proof of~\ref{thm-cocyclesheaf}. 
To prepare for this we start with a variation on the standard 
nerve construction. 
Recall that $\sD\uli n$ is the poset of nonempty subsets of 
$\uli n= \{0,1,2,\dots,n\}$. There are 
functors $v_n\co\sD\uli n\to\uli n$ 
given by  $v_n(S)=\max(S)\in \uli n$. 

\medskip
\begin{lem}
\label{lem-nerveconversion} Let $\sC$ be a 
small category. Then the map of simplicial sets
\[(\,n\mapsto \hom(\uli n\op,\sC))
\quad\lra\quad (\,n\mapsto\hom(\sD\uli n\op,\sC)) \]
given by composition with $v_{\bullet}$ induces a 
homotopy equivalence of the geometric realizations. 
\end{lem} 

\proof The simplicial set $(n\mapsto\hom(\sD\uli n\op,\sC))$ 
is obtained by applying 
Kan's functor \emph{ex}, the right adjoint of the barycentric
subdivision, to $(n\mapsto\hom(\uli n\op,\sC))$. The statement 
is therefore a special case of \cite[3.7]{Kan57}. \qed

We note that the simplicial set $(n\mapsto \hom(\uli n\op,\sC))$
is precisely the nerve of $\sC$, denoted $N_{\bullet}\sC$ in 
section~\ref{sec-Vassiliev}. 

\medskip 
\begin{cor}
\label{cor-binerveconversion} 
Let $m\mapsto \sC_m$ be a simplicial 
category. The map of bisimplicial sets 
\[ 
\xymatrix@M+2pt{
(m,n) \ar@{|->}[r] & {\hom(\uli n\op,\sC_m)} \ar@<-11ex>@{=>}[d] \\
(m,n) \ar@{|->}[r] & {\hom(\sD\uli n\op,\sC_m)}
}
\]
given by composition with the functors $v_n\co\sD\uli n\to\uli n$
induces a homotopy equivalence of the geometric realizations.  \qed
\end{cor} 

\begin{lem}
\label{lem-simpinflation}
Let $S$ be an infinite set and let $K_{\bullet}$ 
be a simplicial set. For $n=1,2,\dots$ let $\emb(\uli n,S)$ be the 
set of injective maps from $\uli n$ to $S$. 
The geometric realization $|K_{\bullet}|$ is 
homotopy equivalent to the geometric realization of the 
incomplete simplicial set alias $\Delta$--set {\rm\cite{RourkeSanderson71}}
\[ n \mapsto K_n\times\emb(\uli n,S) . \]
\end{lem}

\proof 
There is a projection map $p$
from the realization of $n \mapsto K_n\times\emb(\uli n,S)$ to 
$|K_{\bullet}|$. We will show that it has contractible fibers. 
This gives the induction step in an inductive argument showing 
that the restrictions $p^{-1}(|K_{\bullet}|^n)\lra |K_{\bullet}|^n$ 
of $p$ are homotopy equivalences for all $n\ge 0$. \newline
Let $y$ be a point in an $m$-cell of $|K_{\bullet}|$, corresponding 
to some nondegenerate simplex in $K_m$. The fiber of $p$ over $y$ 
is homeomorphic to the classifying space of the poset $\sR$
whose elements are the nonempty finite subsets of $S$ equipped with a
total ordering and an order preserving
surjection to $\uli m$. For each finite subset 
$\sR'$ of $\sR$, there 
exists $T\in \sR$ which is disjoint from all $T'\in \sR'$, 
so that $T'\le T'\cup T\ge T$ in $\sR$ where $T'\cup T$ has the 
concatenated ordering. This implies that the inclusion of $|\sR'|$ 
in $|\sR|$ is homotopic to a constant map, with value equal 
to the vertex determined by $T$. Therefore $\sR$ is contractible, 
i.e., the fiber in question is contractible. \qed

\begin{cor} 
\label{cor-otherdiagonal} Let $J$ be the (fixed) infinite set 
from definition~\ref{dfn-sheafnerve}. The classifying space $B|\cF|$ is 
homotopy equivalent to the geometric realization of 
the incomplete simplicial set $K_{\bullet}$ given by 
\[ n \quad \mapsto \quad \hom(\sD\uli n\op,\cF(\Delta^n_e))
\times\emb(\uli n,J). \]
\end{cor} 

\medskip
We turn to the construction of a comparison map $\Psi$ from 
the incomplete simplicial set $K_{\bullet}$ 
in corollary~\ref{cor-otherdiagonal} to 
the simplicial set $\,\uli n\, \mapsto \beta\cF(\Delta^n_e)$.
An $n$-simplex in $K_{\bullet}$ consists of a functor 
\[ \varphi\co \sD\uli n\op\lra \cF(\Delta^n_e) \]
and an injective map 
$\lambda\co\uli n \to J$. The pair $(\varphi,\lambda)$ carries exactly the 
same information as an element in $\beta\cF(\Delta^n_e)$ 
whose underlying $J$-indexed open covering is given by 
$j\mapsto \Delta^n_e$ if $j=\lambda(t)$ for some $t\in\uli n$ 
and $j\mapsto\emptyset$ otherwise. To make these data 
functorial, i.e., compatible with face operators,
we need to replace the nonempty open sets in the open covering by smaller 
ones, according to the rule 
\begin{equation*}
\label{eqn-coverprecaution}
  j=\lambda(t) \quad \mapsto \quad \{\,(x_0,x_1,\dots,x_n)\in \Delta^n_e\mid
x_t>0\}. \ttag
\end{equation*}
The remaining data can be restricted and we now have an 
element $\Psi(\varphi,\lambda)\in \beta\cF(\Delta^n_e)$. 
The construction $\Psi$ respects the face operators. We 
restate theorem~\ref{thm-cocyclesheaf} as

\begin{thm} 
\label{thm-Thetabijective}
The map $\Psi$ induces a homotopy equivalence from
$B|\cF|\simeq|K_{\bullet}|$ to $|\beta\cF|$.
\end{thm} 

\proof We proceed 
by constructing an inverse at the homotopy set level, 
a natural map $\Lambda$ from $\beta\cF[X]\cong [X,|\beta\cF|\,]$
to $[X,|K_{\bullet}|\,]$. Let $(\sY,\varphi\bbul)$
be an element of $\beta\cF(X)$, so that $\sY$ is a locally finite open 
covering of the manifold $X$. Choose a smooth triangulation of $X$, 
with vertex set contained in $J$.
For each finite nonempty subset $S$ of $J$ spanning a 
simplex of the triangulation, choose a smooth map 
$c_{e,S}\co \Delta_e(S)\to X$ extending the characteristic 
inclusion $c_S\co \Delta(S)\to S^n$,
in such a way that $c_{e,S}$ agrees with 
$c_{e,R}$ on a face $\Delta_e(R)\subset\Delta_e(S)$.
All this is to be done in such a way that   
a map $\kappa\co J\to J$ can be found satisfying
\[ c_{e,S}(\Delta_e(S)) \subset Y_{\kappa(v)} \]
whenever $S\subset J$ spans a simplex and $v\in S$. 
Then for each pair of nonempty subsets $Q,R$ 
of $S$ with $Q\subset R$, the pullback under $c_{e,S}$ of the morphism 
$\varphi_{\kappa(Q)\kappa(R)}$ in $\cF(Y_{\kappa(R)})$ is a morphism 
in $\cF(\Delta_e(S))$. Together these morphisms define an element
\[ x_S\in \hom(\sD(S)\op,\cF(\Delta_e(S))) \, . \]
Finally we choose a total ordering of the vertex set of the 
triangulation. This promotes each $x_S$ to a simplex of $K_{\bullet}$, 
and as these simplices are compatibly constructed they determine a 
map from $X$ (viewed as a simplicial complex) to $|K_{\bullet}|$. 
It follows from lemma~\ref{lem-simpinflation} that the homotopy class 
of that map does not depend on the way in which the vertex set 
of the triangulation is embedded in $J$, and then it 
is altogether clear that the homotopy class depends only on the 
concordance class of $(\sY,\varphi\bbul)\in \beta\cF(X)$. Hence we have defined 
$\Lambda\co \beta\cF[X]\to [X,|K_{\bullet}|\,]$. \newline
We need to know that the composition $\Psi\Lambda\co \beta\cF[X]
\to \beta\cF[X]$ is the identity. 
Suppose that an element of $\beta\cF[X]$ is represented by a pair 
$(\sY,\varphi\bbul)$, where $\sY$ is a $J$--indexed open covering 
of $X$. Then, by construction and painful inspection,  
$\Psi\Lambda$ of that element is represented by a pair
$(\sY',\varphi\bbul')$ for which $\kappa\co J\to J$ can be found 
such that $Y'_j\subset Y_{\kappa(j)}$ for $j\in J$ and 
$\varphi'_{RS}$ is the restriction of $\varphi_{\kappa(R)\kappa(S)}$
to $Y'_{S}$, for finite nonempty $R,S\subset J$ with $R\subset S$. 
[What makes the inspection difficult is that the target of 
$\Psi$ is really $[X,|\beta\cF|\,]$ and we compose with the 
identification $\vartheta\co [X,|\beta\cF|\,]\to \beta\cF[X]$
of section~\ref{subsec-concordanceandrepresent}. At the level 
of representatives, we are therefore dealing with 
$\vartheta\Psi\Lambda$. Each nonempty $Y'_j$ is a subset of the 
open star $\st(j)$ of the vertex labelled $j$ in the triangulation of 
$X$ used to construct $\Lambda$ of $(\sY,\varphi\bbul)$, above.
In fact $Y'_j$ is the inverse image of $\st(j)$ under $h_1\co X\to X$, 
and $h_1$ is part of a homotopy $(h_t)$ as specified in 
the description of $\vartheta$.] 
Thus we have a situation where one element of $\beta\cF(X)$
``refines'' another. Lemma~\ref{lem-morphismstoconcordances} below 
then guarantees that the 
two elements are concordant. Hence $\Psi\Lambda=\id$. \newline 
Next we show that $\Lambda\co \beta\cF[X]\to [X,|K_{\bullet}|\,]$ is 
onto for any $X$ in $\sX$. 
As explained in section~\ref{subsec-concordanceandrepresent}, 
any element of $[X,|K_{\bullet}|\,]$ can be represented by
a simplicial map $f\co X^s\to K_{\bullet}$ where $X^s$ is the simplicial
set associated to some smooth triangulation of $X$ with 
totally ordered vertex set. We subject $f$ to a smoothing 
procedure, familiar from section~\ref{subsec-concordanceandrepresent}, 
which will result in a ``better'' simplicial map 
$f^{\sim}\co X^s\to K_{\bullet}$ representing the same homotopy class. 
The smoothing procedure begins with the choice of a homotopy 
of smooth maps $h_t\co X\to X$, where 
$0\le t\le 1$, such that $h_0=\id$ and 
\begin{enumerate} 
\item for every $t$, the map
$h_t$ maps each simplex of the triangulation to itself;
\item each simplex 
has a neighbourhood in $X$ which is mapped to the simplex 
by $h_1$. 
\end{enumerate} 
We also choose compatible smooth embeddings $c_{e,R}\co \Delta_e(R)\to X$
extending the characteristic embeddings $c_R\co \Delta(R)\to X$
of the simplices of $X$, and such that each 
$h_t\co X\to X$ takes the image of each $c_{e,R}$ to itself.
For an $n$--simplex $z_R$ of $X^s$, with vertex set $R\subset T$, define  
$f^{\sim}(z_R)$ by composing the functor
\[ (c_{e,R}^{-1}\circ h_1\circ c_{e,R})^*\co \cF(\Delta^n_e)\to
\cF(\Delta^n_e) \]
with $f(z_R)\co\sD\uli n\op\to \cF(\Delta^n_e)$~; here we have identified
$S$ with $\uli n$ and $\Delta_e(S)$ with $\Delta^n_e$ as usual. 
As in the proof of 
proposition~\ref{prp-xiandtheta}, lemma~\ref{lem-evhomotopy} 
shows that $f$ and $f^{\sim}$ represent the same homotopy class. 
But now it is easy to exhibit an element $(\sY,\varphi\bbul)$ of 
$\beta\cF(X)$ which is mapped to the class of $f^{\sim}$ by $\Lambda$. We 
may assume that the vertex set $T$ of the triangulation is contained in
$J$ and let 
\[  Y_j = \intr((h_1)^{-1}(\st(j)))  \]
for $j\in T$, where $\st(j)$ is the open star of the vertex $j$. 
All other $Y_j$ are empty. To obtain $\varphi_{RS}$, assuming 
$Y_S\ne \emptyset$, we note that 
\[ h_1(Y_S) \subset \bigcap_{j\in S}\st(j) \]
so that $S$ is the vertex set of an $n$--simplex $z_S$ in $X$ for 
some $n$. We therefore have 
\[ c_{e,S}^{-1}\circ h_1\co Y_S\to \Delta_e(S) \]
and we can use it to pull back the morphism in $\cF(\Delta_e(S))
\cong \cF(\Delta^n_e)$
which is the image of $R\subset S$ under the functor 
$f(z_S)$. The result is $\varphi_{RS}$, a morphism in $\cF(Y_S)$. 
Following the instructions above for finding a representative 
for $\Lambda$ of $(\sY,\varphi\bbul)$, we get precisely $f^{\sim}$. 
The conclusion is that $\Lambda$ is indeed surjective. \newline
The final step is to note that $\Psi$, as a map from 
$|K_{\bullet}|$ to $|\beta\cF|$, induces a surjection 
$\pi_1(|K_{\bullet}|,z) \to \pi_1(|\beta\cF|,\Psi(z))$
for any choice of base vertex $z\in |K_{\bullet}|$. We leave this 
verification to the reader: Given an element $u$ in 
$\pi_1(|\beta\cF|,\Psi(z))$, an element $v$ of 
$\pi_1(|K_{\bullet}|,z)$
can be obtained by applying the procedure $\Lambda$ above to $u$ in a 
relative form. The relative case of 
lemma~\ref{lem-morphismstoconcordances} below implies 
$\Psi(v)=u$. 
\newline
It is a formality to show that a map $q\co C\to D$ between CW--spaces 
which induces bijections $[X,C]\to [X,D]$ for every $X$ in $\sX$ 
and surjections $\pi_1(C,z)\to \pi_1(D,q(z))$ for every  
$z\in C$ induces bijections $\pi_n(C,z)\to \pi_n(D,q(z))$ for 
$n\ge 0$ and $z\in C$. Such a map is therefore a homotopy 
equivalence. We have just verified that this criterion applies 
with $q=\Psi$, showing that $\Psi\co |K_{\bullet}|\to |\beta\cF|$ 
is a homotopy equivalence. \qed

\begin{lem}
\label{lem-morphismstoconcordances}
Let $(\sY,\varphi\bbul)$ and 
$(\sY',\varphi\bbul')$ be elements 
of $\beta\cF(X)$. Suppose that there exists a map $\kappa\co J\to J$ 
such that $Y'_j\subset Y_{\kappa(j)}$ for all $j\in J$, 
and $\varphi'_{RS}$ is the restriction of
$\varphi_{\kappa(R)\kappa(S)}$
to $Y'_S$, for all finite nonempty $R,S\subset J$ with $R\subset S$. 
Then $(\sY,\varphi\bbul)$ and $(\sY',\varphi\bbul')$ are
concordant.  
If $(\sY,\varphi\bbul)$ and $(\sY',\varphi\bbul')$ are 
in $\beta\cF(X,A;z)$ for some closed $A\subset X$ and some 
$z\in \beta\cF(\pt)$, and if $\kappa(j)=j$ for all $j\in J$ 
such that the closure of $Y_j$ has nonempty intersection with $A$,
then the concordance can be taken relative to $A$. 
\end{lem} 

\proof We assume first that the fixed indexing set 
$J$ is uncountable, rather than just infinite, and concentrate 
on the absolute case, $A=\emptyset$. \newline 
The case where $\kappa=\id_J$ is straightforward.
Hence $(\sY',\varphi\bbul')$ is concordant to $(\sY'',
\varphi\bbul'')$ 
where $Y''_j=Y_{\kappa(j)}$ and $\varphi''_{RS}=
\varphi_{\kappa(R)\kappa(S)}$. It remains to find a concordance 
from $(\sY'',\varphi\bbul'')$ to $(\sY,\varphi\bbul)$. Alternatively, 
to keep notation under control, we may assume from now 
on that $(\sY',\varphi\bbul')
=(\sY'',\varphi\bbul'')$, in other words $Y'_j=Y_{\kappa(j)}$ 
for all $j\in J$. \newline
The sets $\{j\in J\mid Y'_j\ne\emptyset\}$ and 
$\{i\in J\mid Y_i\ne \emptyset\}$ are countable, since 
the coverings $\sY'$ and $\sY$ are locally finite and $X$ admits a 
countable base. Hence there exists a bijection $\lambda\co J\to J$ such that 
$Y_{\lambda(j)}\cap Y_j=\emptyset = Y_{\lambda(j)}\cap
Y_{\kappa(j)}$
for all $j\in J$; for example, $\lambda$ can be chosen so that 
$Y_\lambda(j)=\emptyset$ if $Y_j\ne \emptyset$ or 
$Y_{\kappa(j)}\ne\emptyset$. Now let
\[ W_j= \big(Y_j\times\,]-\infty,1/2[\,\big)
\cup \big(Y_{\lambda(j)}\times\,]1/4,3/4[\,\big) \cup
\big(Y_{\kappa(j)}\times\,]1/2,\infty[\,\big). \]
The $W_j$ for $j\in J$ constitute an open covering $\sW$ of $X\times\RR$. 
For any finite nonempty $S\subset J$, we have a decomposition of $W_S$ 
into disjoint open sets  
\[
\begin{array}{lll}
Y_S\times\,]-\infty,1/2[\,, & 
Y_{\lambda(S)}\times\,]1/4,3/4[\,, &
Y_{\kappa(S)}\times\,]1/2,\infty[\,, \\
Y_{Q\cup\lambda(S\smin Q)}\times\,]1/4,1/2[\,, &
Y_{\lambda(Q)\cup\kappa(S\smin Q)}\times\,]1/2,3/4[\,, &
\end{array}
\]
where $Q$ runs through the nonempty proper subsets of $S$.
Therefore, given finite 
nonempty $R,S\subset J$ with $R\subset S$, there is a unique morphism 
$\psi_{RS}$ in $\cF(W_S)$ whose restrictions 
to the various summands of $W_S$ in the above decomposition
are the pullbacks of 
$\varphi_{RS}$, $\varphi_{\lambda(R)\lambda(S)}$,
$\varphi_{\kappa(R)\kappa(S)}$, etc. etc., under the projections 
to $Y_S$, $Y_{\lambda(S)}$, $Y_{\kappa(S)}$, $
Y_{Q\cup\lambda(S\smin Q)}$ and $Y_{\lambda(Q)\cup\kappa(S\smin Q)}$, 
respectively. (Here the two ``etc.'' are short for $\varphi_{TU}$ where 
$U=Q\cup\lambda(S\smin Q)$ and $T=(R\cap Q)\cup \lambda(R\smin Q)$
in the first case, while $U=\lambda(Q)\cup \kappa(S\smin Q)$ and 
$T=\lambda(R\cap Q)\cup \kappa(R\smin Q)$ in the second case.)
Clearly $(\sW,\psi\bbul)$ is a concordance 
from $(\sY,\varphi\bbul)$ to $(\sY',\varphi\bbul')$. \newline 
Next we look at the relative case, $A\ne\emptyset$, but continue 
to assume that $J$ is uncountable. As in the absolute case we may 
assume that $Y'_j=Y_{\kappa(j)}$ for all $j\in J$. Choose a 
bijection $\lambda\co J\to J$ such that $\lambda(j)=j$ 
whenever $\kappa(j)=j$, and such that $Y_j\cap Y_{\lambda(j)}
=\emptyset = Y_{\kappa(j)}\cap Y_{\lambda(j)}$ for the remaining 
$j$. Again let
\[ W_j= \big(Y_j\times\,]-\infty,1/2[\,\big)
\cup \big(Y_{\lambda(j)}\times\,]1/4,3/4[\,\big) \cup
\big(Y_{\kappa(j)}\times\,]1/2,\infty[\,\big). \]
The $W_j$ for $j\in J$ constitute an open covering $\sW$ of $X\times\RR$. 
For a finite nonempty $S\subset J$ which is contained in the fixed
point set of $\kappa$, we simply have $W_S=Y_S\times\RR$.
For a finite nonempty $S\subset J$ which does not 
contain any fixed points of $\kappa$, we have a decomposition of $W_S$ 
into disjoint open sets   
\[
\begin{array}{lll}
Y_S\times\,]-\infty,1/2[\,, & 
Y_{\lambda(S)}\times\,]1/4,3/4[\,, &
Y_{\kappa(S)}\times\,]1/2,\infty[\,, \\
Y_{Q\cup\lambda(S\smin Q)}\times\,]1/4,1/2[\,, &
Y_{\lambda(Q)\cup\kappa(S\smin Q)}\times\,]1/2,3/4[\,, &
\end{array}
\]
as before, where $Q$ runs through the nonempty proper subsets of $S$. 
For finite nonempty $S\subset J$ which contains some 
fixed points of $\kappa$ and some non--fixed points of $\kappa$,
write $S=S_1\cup S_2$ 
where $S_1=\{j\in S\mid \kappa(j)=j\}$ and $S_2=S\smin S_1$. 
Then $W_{S_2}$ decomposes into disjoint open sets as above, 
whereas $W_{S_1}=Y_{S_1}\times\RR$. Hence $W_S=W_{S_1}\cap W_{S_2}$
still decomposes as a disjoint union of open sets 
\[
\begin{array}{lll}
Y_S\times\,]-\infty,1/2[\,, & 
Y_{\lambda(S)}\times\,]1/4,3/4[\,, &
Y_{\kappa(S)}\times\,]1/2,\infty[\,, \\
Y_{Q\cup\lambda(S\smin Q)}\times\,]1/4,1/2[\,, &
Y_{\lambda(Q)\cup\kappa(S\smin Q)}\times\,]1/2,3/4[\,, &
\end{array}
\]
where $Q$ runs through the nonempty proper subsets of $S_2$ only. 
We can therefore define morphisms 
$\psi_{RS}$ in $\cF(W_S)$ much as in the absolute case and obtain 
a relative concordance $(\sW,\psi\bbul)$ from $(\sY,\varphi\bbul)$ to 
$(\sY',\varphi\bbul')$. \newline
Now we must consider the case(s) where $J$ is countably infinite. 
We can in fact reason as before provided that $X$ is a closed
manifold, because 
in that case the sets $\{j\in J\mid Y'_j\ne\emptyset\}$ 
and $\{i\in J\mid Y_j\ne \emptyset\}$ are finite. While this is not 
exactly what we want, it allows us to make a comparison between the 
case where $J$ is countable and the case where it is uncountable. 
To this end, choose an uncountable set $J^{\sharp}$ 
containing $J$ as a subset. Corresponding to $J$ and $J^{\sharp}$
we have two variants of $\beta\cF$. 
We keep the notation $\beta\cF$ for the $J$--variant, and write 
$\beta^{\sharp}\cF$ for the $J^{\sharp}$--variant. There is a natural 
inclusion $\beta\cF(X)\to \beta^{\sharp}\cF(X)$; namely, any  
$J$--indexed open covering of $X$ can be regarded as a
$J^{\sharp}$--indexed covering of $X$ where all open sets with 
labels in $J^{\sharp}\smin J$ are empty. By all the above, 
$|\beta\cF|\to |\beta^{\sharp}\cF|$ induces an isomorphism of 
homotopy groups or homotopy sets, for any choice of base vertex
in $|\beta\cF|$, the point being that spheres are closed manifolds. 
By proposition~\ref{prp-whatitrepresents}, 
this implies that the inclusion--induced map of concordance sets
\[  \beta\cF[X,A;z] \lra \beta^{\sharp}\cF[X,A;z] \]
is always a bijection, and not just when $X$ is closed. 
We have therefore reduced the case of a countable $J$ to the 
case of an uncountable one, and that has been dealt with. \qed

\section{Realization and homotopy colimits}
\label{sec-barconstruction}

\subsection{Realization and squares}

\begin{lem} 
\label{lem-simplicialfibration}
Let $u_{\bullet}\co E_{\bullet}\lra B_{\bullet}$ be a map 
between incomplete simplicial spaces (or good simplicial 
spaces). Suppose that the squares 
\[
\xymatrix{
E_k \ar[r]^-{u_k} \ar[d]^{d_i} & B_k \ar[d]^{d_i} \\
E_{k-1} \ar[r]^-{u_{k-1}} & B_{k-1} 
}
\]
are all homotopy cartesian ($k\ge i\ge 0$). Then the
following is also homotopy cartesian: 
\[
\xymatrix@C+4pt{
E_0  \ar[r]^-{u_0} \ar[d]^{\textup{incl.}}  &  B_0 \ar[d]^{\textup{incl.}} \\
|E_{\bullet}| \ar[r]^-{|u_{\bullet}|} &  |B_{\bullet}|. 
}
\]
\end{lem}

\begin{lem}
\label{lem-simplicialhomologyfibration}
Let $u_{\bullet}\co E_{\bullet}\lra B_{\bullet}$ be a map 
between incomplete simplicial spaces (or good simplicial spaces). 
Suppose that, in each square 
\[
\xymatrix{
E_k \ar[r]^-{u_k} \ar[d]^{d_i} &  B_k \ar[d]^{d_i} \\
E_{k-1} \ar[r]^{u_{k-1}} & B_{k-1} 
}
\]
the canonical map from any homotopy fiber of $u_k$ to the 
corresponding homotopy fiber of $u_{k-1}$ induces an isomorphism in 
integer homology. Then in the square 
\[
\xymatrix{
E_0 \ar[r]^-{u_0} \ar[d]^{\textup{incl.}} &  B_0
\ar[d]^{\textup{incl.}} \\
|E_{\bullet}| \ar[r]^-{|u_{\bullet}|} & |B_{\bullet}|,
}
\]
the canonical map from any homotopy fiber of $u_0$ to 
the corresponding homotopy fiber of $|u_{\bullet}|$ induces 
an isomorphism in integer homology. 
\end{lem}

\proof[Proofs.] It is shown in
\cite[1.6]{Segal74} and \cite[Prop.4]{McDuffSegal76} that the
geometric realization procedure for simplicial spaces respects
degreewise quasifibrations and homology fibrations under 
reasonable conditions. The two lemmas follow from these 
statements upon converting the maps $u_k$ into fibrations. \qed

\begin{cor} 
\label{cor-hocolimhomologyfibration} Let $\sC$ be a small category 
and let $u\co\cG_1\to \cG_2$ be a natural transformation 
between functors from $\sC$ to spaces. Suppose that, for each morphism
$f\co a\to b$ in $\sC$, the map $f_*$ from any homotopy 
fiber of $u_a$ to the corresponding homotopy fiber of $u_b$ 
induces an isomorphism in integer homology. Then for each object $a$
of $\sC$, the inclusion of any homotopy fiber of $u_a$ in the 
corresponding homotopy fiber of 
$u_*\co\hocolim\,\cG_1\to \hocolim\,\cG_2$ induces an 
isomorphism in integer homology.   
\end{cor} 

\proof Apply lemma~\ref{lem-simplicialhomologyfibration} with 
$E_k:=\coprod \cG_1(D(k))$ and $B_k=\coprod \cG_2(D(k))$, where both 
coproducts run over the set of contravariant functors $D$ from the 
poset $\uli k$ to $\sC$. Then $|E_{\bullet}|$ is $\hocolim\, \cG_1$ 
and $|B_{\bullet}|$ is $\hocolim\, \cG_2$. \qed

\subsection{Homotopy colimits} 
\label{subsec-hocolim}
Any functor $\cD$ from a small (discrete) category $\sC$ to the category of 
spaces has a {\it colimit}, $\colim\, \cD$. This is the quotient 
space of the coproduct 
\[ \coprod_{a\textup{ in }\sC} \, \cD(a) \]
obtained by identifying $x\in \cD(a)$ with $f_*(x)\in \cD(b)$ 
for any morphisms $f\co a\to b$ in $\sC$ and elements $x\in \cD(a)$. It is
well known that the colimit construction is not well behaved
from a homotopy theoretic point of view. Namely,
suppose that $w\co \cD_1\to \cD_2$ is a
natural transformation between functors from $\sC$ to spaces 
and that $w_a:\cD_1(a)\to \cD_2(a)$ is a homotopy equivalence for any 
object $a$ in $\sC$. Then this does not in general imply that the
map induced by $w$ from $\colim\, \cD_1$ to $\colim\, \cD_2$  
is again a homotopy equivalence. (It is easy to make examples 
with $\sC$ equal to the poset of proper subsets of a two--element 
set, so that the colimits become pushouts.)

\medskip
Call a functor $\cD$ from $\sC$ to spaces {\it cofibrant} if, for any 
diagram of functors (from $\sC$ to spaces) and natural transformations
\[
\xymatrix{
\cD \ar[r]^v & \cE  & \ar[l]_w \cF 
}
\]
where $w_a\co \cF(a)\to \cE(a)$ is a homotopy equivalence for all 
$a\in \sC$, there exists a natural transformation 
$v'\co \cD\to \cF$ and a {\it natural} homotopy 
$\cD(a)\times[0,1]\to \cE(a)$ (for 
all $a$) connecting $wv'$ and $v$. 
It is not hard to show the following. If $v\co \cD_1\to \cD_2$ is a
natural transformation between cofibrant functors 
such that $v_a\co \cD_1(a)\to \cD_2(a)$ is a 
homotopy equivalence for each $a\in \sC$, then $v$ has 
a natural homotopy inverse (with {\it natural} homotopies) and therefore
the induced map $\colim\, \cD_1\to\colim\, \cD_2$ is a homotopy equivalence. 

\medskip
This suggests the following procedure for making colimits 
homotopy invariant. Suppose that $\cD$ from $\sC$ to spaces is any
functor. Try to
find a natural transformation 
$\cD'\to \cD$ specializing to homotopy equivalences 
$\cD'(a)\to \cD(a)$ for all $a$ in $\sC$, 
where $\cD'$ is cofibrant. Then define the
\emph{homotopy colimit} of $\cD$ to be 
$\colim\, \cD'$. 
If it can be done, $\hocolim \,\cD$ is at least well defined
up to homotopy equivalence. 

\medskip 
This point of view is carefully presented in \cite{Dror87}. 
Some of the ideas go back to \cite{Milnor62}.
As we will see in a moment, there is a  
construction for $\cD'$ which depends naturally on $\cD$.

\medskip
The standard foundational reference for homotopy colimits and homotopy 
limits is the book \cite{BousfieldKan72} by Bousfield and Kan.
But the first explicit construction of homotopy colimits in general 
appears to be due to Segal \cite{Segal68}.  

Again let $\cD$ be a functor from a discrete 
small category $\sC$ to the category of spaces. Following 
Segal we introduce a topological category denoted $\sC\ssmallint\cD$,
the \emph{transport category} of $\cD$:
\begin{equation*}
\label{eqn-transportcategory}
\ob(\sC\ssmallint\cD) \,= \, \coprod_{a\in \ob(\sC)} \cD(a)\,, \qquad 
\mor(\sC\ssmallint\cD) \,=\, 
\coprod_{f\in \mor(\sC)} \cD(\sigma(f))\,. 
\end{equation*}
Here $\sigma(f)$ denotes the source of a morphism $f$ in $\sC$.
We will write morphisms in $\sC\ssmallint\cD$ as pairs $(f,x)$ where
$f\in \mor(\sC)$ and $x\in \cD(\sigma(f))$. 
The composition $(g,y)\circ(f,x)$ of two 
such morphisms is defined if an only if $g\circ f$ is defined in
$\sC$ and $f_*(x)=y$, in which case $(g,y)\circ(f,x)=(g\circ f,x)$. 
The classifying space $B(\sC\ssmallint\cD)$ is a model for
the homotopy colimit of $\cD$. 

\medskip
To relate $B(\sC\ssmallint\cD)$ to our earlier discussion we define a functor 
$\cD'$ from $\sC$ to spaces as follows. 
For $a\in\ob(\sC)$ let $\sC\!\downarrow\! a$
be the category of $\sC$-objects over $a$, \cite[II.6]{MacLane71}. 
Let 
\[
\begin{array}{c}
\cD'(a) := B\left((\sC\!\downarrow\! a)\ssmallint \cD\right) 
\end{array}
\]
for objects $a$ in $\sC$, where we view $\cD$ as a functor on 
$\sC\!\downarrow\! a$. Then $\cD'$ is cofibrant and the 
canonical map $\cD'(a)\to \cD(a)$ is a homotopy equivalence 
for every $a$ in $\sC$. Moreover,  
\[ 
\begin{array}{c}
B(\sC\ssmallint\cD) \cong \colim\,\cD' . 
\end{array}
\]
   
Note in passing that if $\cD(a)$ is a singleton for each $a$ in $\cC$, 
then the transport category $\sC\ssmallint\cD$ is identified with $\sC$ and 
so $\hocolim\,\cD= B\sC$.  

\begin{prp} Let $w\co \cD_1\to \cD_2$ be a
natural transformation between functors from $\sC$ to spaces. 
Suppose that $w_a:\cD_1(a)\to \cD_2(a)$ is a homotopy equivalence for any 
object $a$ in $\sC$. Then the map 
$\hocolim\, \cD_1 \lra \hocolim\, \cD_2$
induced by $w$ is a homotopy equivalence, where 
$\hocolim\, \cD_i = B(\sC\!\downarrow\!\cD_i)$. 
\end{prp} 

This is just a partial summary of our conclusions above. We proceed 
to a reformulation,~\ref{prp-otherhocoliminvariance} below, in which 
homotopy colimits are not mentioned explicitly.

\begin{dfn} 
\label{dfn-easytransportprojection}
{\rm Let $p\co \sE\to \sC$ be a continuous functor between 
small topological categories, where $\sC$ happens to be discrete. We say 
that $p$ is a \emph{transport projection} if the following 
is a pullback square of spaces: 
\[ 
\xymatrix@C=40pt{
\mor(\sE) \ar[r]^{\textrm{source}} \ar[d]^p & \ob(\sE) \ar[d]^p \\
\mor(\sC) \ar[r]^{\textrm{source}} & \ob(\sC) 
}
\]
}
\end{dfn} 

\begin{prp} 
\label{prp-otherhocoliminvariance} Let $p\co \sE\to \sC$ and $p'\co
\sE'\to \sC$ be transport projections as in 
definition~{\rm \ref{dfn-easytransportprojection}}. Let $u\co \sE\to \sE'$ 
be a continuous functor over $\sC$. Suppose also that, for each 
object $c$ in $\sC$, the restriction $\sE_c\to \sE'_c$ of $u$ 
to the fibers over $c$ is a homotopy equivalence. Then 
$Bu\co B\sE\to B\sE'$ is a homotopy equivalence. 
\end{prp} 

\proof Note that $\sE\cong \sC\!\downarrow\!\cD$ and 
$\sE'\cong \sC\!\downarrow\!\cD$ where $\cD(c)=\sE_c$ and 
$\cD'(c)=\sE'_c$ for an object $c$ in $\sC$. Note also that 
$\sE_c$ and $\sE'_c$ are topological categories in which every 
morphism is an identity, that is, they are just spaces. \qed

\medskip 
Next we mention two useful naturality properties 
of homotopy colimits. 
To make a homotopy colimit, we need a pair $(\sC,\cD)$ 
consisting of a small category $\sC$ and a functor $\cD$ 
from $\sC$ to spaces. By a \emph{morphism} from one such 
pair $(\sC^s,\cD^s)$ to another, $(\sC^t,\cD^t)$, we understand 
a pair $(\cF,\nu)$ consisting of a functor $\cF\co \sC^s\to \sC^t$ 
and a natural transformation $\nu$ from $\cD^s$ to $\cD^t\cF$. 

\begin{rmk}
\label{rmk-naturalhocolim1}
Such a morphism induces a map $(\cF,\nu)_*$ from 
$\hocolim\, \cD^s$  to $\hocolim\,\cD^t$.  
\end{rmk}

Suppose that $(\cF_0,\nu_0)$ and $(\cF_1,\nu_1)$ are morphisms 
from $(\sC^s,\cD^s)$ to $(\sC^t,\cD^t)$. Let $\theta$ be 
a natural transformation from $\cF_0$ to $\cF_1$ 
such that $\nu_1=\cD^t(\theta)\circ\nu_0$. 

\begin{rmk} 
\label{rmk-naturalhocolim2}
Such a $\theta$ induces a homotopy $\theta_*$ from 
$(\cF_0,\nu_0)_*$ to $(\cF_1,\nu_1)_*$. 
\end{rmk}

\proof Let $\sI=\{0,1\}$, viewed as an ordered set with the 
usual order and then as a category. Then $B\sI\cong[0,1]$. 
Let $p\co \sC\times\sI\to \sC$ 
be the projection. The data $(\cF_0,\nu_0)$, $(\cF_1,\nu_1)$
and $\theta$ together define a morphism from 
$(\sC^s\times\sI,\cD^s\circ p)$ to $(\sC^t,\cD^t)$. By 
remark~\ref{rmk-naturalhocolim1}, this induces a map 
from $\hocolim\,(\cD^s\circ p) \cong (\hocolim\, \cD^s)\times B\sI$
to $\hocolim\,\cD^t$. \qed

\bigskip
Let $\sC$ be a small category and let $a\mapsto \cF_a$ 
be a covariant functor from $\sC$ to the category of sheaves on
$\sX$. 

\begin{lem} 
\label{lem-internalhocolim}
$|\hocolim_a\,\cF_a|\,\simeq\, \hocolim_a\, |\cF_a|$. 
\end{lem}

\proof Definition~\ref{dfn-earlysheafhocolim} and 
theorem~\ref{thm-cocyclesheaf} give
$|\hocolim_a\,\cF_a|\,\simeq\, B|\sC\ssmallint\cF|$ and 
propositions~\ref{prp-productsheaf},~\ref{prp-coproductsheaf}
imply $\textstyle B|\sC\ssmallint\cF| \cong 
B(\sC\ssmallint\,|\cF_{\bullet}|\,)$, where $|\cF_{\bullet}|$ 
denotes the functor 
$a\mapsto |\cF_a|$ from $\sC$ to spaces. \qed

\begin{cor}
\label{cor-hoinvofhocolim}
Let $\sC$ be a small category and let 
$a\mapsto \cE_a$ and $\mapsto \cE'_a$ be covariant functors from $\sC$ 
to the category of sheaves on
$\sX$. Let $\nu=\{\nu_a\co \cE_a\to \cE'_a\}$ be a natural 
transformation such that every $\nu_a\co \cE_a\to \cE'_a$ is a weak 
equivalence. 
Then the induced map $\hocolim_a\,\cE_a\to 
\hocolim_a\,\cE'_a$ is a weak equivalence (between 
sheaves on $\sX$). \qed
\end{cor}


\medskip
\begin{thebibliography}
 
\bibitem{BousfieldKan72} {\bf A K Bousfield, D Kan} {\it Homotopy 
limits, completions and localizations}, Springer Lecture Notes in
Math., vol. 304, Springer (1972)
\bibitem{Brown62} {\bf E Brown} {\it Cohomology theories}, Ann. of
Math. 75 (1962), 467--484; correction, Ann. of
Math. 78 (1963), 201
\bibitem{BroeckerJaenich73} {\bf T Br\"{o}cker, K J\"{a}nich} 
{\it Introduction to Differential
Topology}, Engl. edition Cambridge Univ. Press (1982); German 
edition Springer-Verlag (1973)
\bibitem{Dror87} {\bf E Dror} {\it Homotopy and homology 
of diagrams of spaces}, in: {\it Algebraic Topology, Seattle 
Wash. 1985}, Springer Lec. Notes in Math. 1286 (1987), 93--134 
\bibitem{EarleEells69} {\bf C J Earle, J Eells}
{\it A fibre bundle description of Teichmueller 
theory}, J. Differential Geom. 3 (1969), 19--43
\bibitem{EarleSchatz70} {\bf C J Earle, A Schatz} 
{\it Teichmueller theory for surfaces with 
boundary}, J. Differential Geom. 4 (1970), 169--185
\bibitem{GabrielZisman67} {\bf P Gabriel, M Zisman} 
{\it Calculus of fractions and homotopy theory}, Springer--Verlag,
Ergebnisse series (1967)
\bibitem{Galatius04} {\bf S Galatius} {\it Mod $2$ homology of the 
stable spin mapping class group}, University of Aarhus preprint (2004)
\bibitem{Galatius02} {\bf S Galatius} {\it Mod $p$ homology 
of the stable mapping class group}, Topology, to appear
\bibitem{GolubitskyGuillemin73} {\bf M Golubitsky, V Guillemin} 
{\it Stable mappings and their singularities}, Springer Grad. Texts (1973)
\bibitem{GoodwillieKleinWeiss02} {\bf T Goodwillie, J Klein, M Weiss},
{\it A Haefliger style description of the embedding calculus 
tower}, Topology 42 (2003), 509--524
\bibitem{Gromov84} {\bf M Gromov} {\it Partial Differential Relations}, 
Springer-Verlag, Ergebnisse series (1984)
\bibitem{Gromov71} {\bf M Gromov} {\it A  topological technique for the
construction of  solutions of differential equations and
inequalities}, Actes du Congr\`{e}s  International des Mathematiciens, 
Nice 1970,  Gauthier--Villars (1971)
\bibitem{Haefliger71} {\bf A Haefliger} {\it Lectures on the 
theorem of Gromov}, in Proc. of 1969/70 Liverpool Singularities Symp., 
Lecture Notes in Math. vol. 209, Springer (1971) 128-141 
\bibitem{Harer85} {\bf J L Harer} {\it Stability of the homology 
of the mapping class groups of oriented surfaces}, Ann. of Math. 121
(1985) 215--249 
\bibitem{Hirsch59} {\bf M Hirsch} {\it Immersions of manifolds}, 
Trans. Amer. Math. Soc. 93 (1959) 242--276
\bibitem{Hirsch76} {\bf M Hirsch} {\it Differential Topology}, 
Springer--Verlag (1976) 
\bibitem{Ivanov89} {\bf N V Ivanov} {\it Stabilization of the 
homology of the Teichmueller modular groups}, Algebra i Analiz 1
(1989) 120--126; translation in: Leningrad Math. J. 1 (1990) 675--691
\bibitem{Kan57} {\bf D Kan} {\it On c.s.s. complexes},
Amer. J. Math. 79 (1957), 449--476. 
\bibitem{Lang72} {\bf S Lang} {\it Differential manifolds}, 
Addison-Wesley (1972)
\bibitem{MacLane71} {\bf S MacLane} {\it Categories for the 
working mathematician}, Grad. texts in Math., Springer-Verlag (1971)
\bibitem{MadsenTillmann01} {\bf I Madsen, U Tillmann} {\it The stable 
mapping class group and $Q(\CC P^{\infty})$}, Invent. Math. 145 (2001)
509--544
\bibitem{McDuffSegal76} {\bf D McDuff, G Segal} {\it Homology 
fibrations and the ``Group--Completion'' theorem}, 
Inventiones Math. 31 (1976), 279--284
\bibitem{Miller86} {\bf E Miller} {\it The homology of the mapping 
class group}, J. Diff. Geom. 24 (1986), 1--14
\bibitem{Milnor62} {\bf J Milnor} {\it On axiomatic homology theory}, 
Pacific J. Math. 12 (1962) 337--341
\bibitem{Milnor63} {\bf J Milnor} {\it Morse Theory}, Ann. of
Math. Studies 51, Princeton University Press (1963, 1969)
\bibitem{Milnor65} {\bf J Milnor} {\it Lectures on the h-cobordism 
theorem}, Mathematical notes 1, Princeton University Press (1965)
\bibitem{Morita84} {\bf S Morita} {\bf Characteristic classes of
surface bundles}, Bull.Amer.Math.Soc. 11 (1984), 386--388
\bibitem{Morita87}  {\bf S Morita} {\bf Characteristic classes of
surface bundles}, Invent. Math. 90 (1987), 551--577
\bibitem{Mumford83} {\bf D Mumford} {\it Towards an enumerative geometry 
of the moduli space of curves}, Arithmetic and geometry, Vol. II, 
Progr. in Math. 36, Birkh\"{a}user (1983), 271--328
\bibitem{Phillips67} {\bf A Phillips} {\it Submersions of open manifolds},
Topology 6 (1967) 170--206
\bibitem{Quillen71} {\bf D Quillen} {\it Elementary proofs of 
some results of cobordism theory using Steenrod operations}, 
Advances in Math. 7 (1971) 29--56
\bibitem{Quillen73} {\bf D Quillen} {\it Higher algebraic
$K$-theory. I}, in Algebraic $K$-theory, I, Proc. of 1972 
Battelle Memorial Inst. conf., Springer 
Lecture Notes in Math. 341 (1973), 85--147
\bibitem{Ravenel84} {\bf D Ravenel} {\it The Segal conjecture for 
cyclic groups and its consequences}, Amer. J. Math. 106 (1984) 415--446
\bibitem{RourkeSanderson71} {\bf C P Rourke, B Sanderson} 
{\it $\Delta$--sets. I. Homotopy theory}, Quart. J. Math. Oxford
22 (1971) 321--338
\bibitem{Segal74} {\bf G Segal} {\it Categories and cohomology theories},
Topology 13 (1974) 293--312
\bibitem{Segal73} {\bf G Segal} {\it Configuration-spaces and 
infinite loop-spaces}, Invent. Math. 21 (1973) 213--221
\bibitem{Segal68} {\bf G Segal} {\it Classifying spaces and 
spectral sequences}, Inst. Hautes Etudes Sci. Publ Math. 34 (1968)
105--112
\bibitem{Smale59} {\bf S Smale} {\it The classification of immersions of 
spheres in Euclidean spaces}, Ann. of Math. 69 (1959) 327--344
\bibitem{Stong68} {\bf R E Stong} {\it Notes on cobordism theory}, 
Mathematical Notes, Princeton University Press (1968)
\bibitem{Tillmann97} {\bf U Tillmann} {\it On the homotopy of the 
stable mapping class group}, Invent. Math. 130 (1997) 257--275
\bibitem{Vassiliev94} {\bf V Vassiliev} {\it Complements of Discriminants of 
Smooth Maps: Topology and Applications}, Transl. of Math. Monographs
Vol. 98, revised edition, Amer. Math. Soc. (1994)
\bibitem{Vassiliev89} {\bf V Vassiliev} {\it Topology of spaces of 
functions without complicated singularities}, Funktsional Anal. 
i Prilozhen 93 no. 4 (1989), 24--36; Engl. translation in Funct. 
Analysis Appl. 23 (1989) 266--286.
\end{thebibliography}
\end{document}